%% file: main.tex
\documentclass[a4paper]{elsarticle}
\usepackage{indentfirst}
\usepackage{multirow}
\usepackage{tikz}
\usepackage{fullpage}
\usetikzlibrary{arrows,shapes}
\usepackage{fancybox}
\usetikzlibrary{shapes.multipart}
\usepackage{graphicx}
\usepackage{subcaption}
\usepackage{float}
\usepackage{pgfplots}
\pgfplotsset{width=10cm,compat=1.9} 
\RequirePackage{amsmath}
\RequirePackage{amssymb}
\RequirePackage{amsthm}
\RequirePackage{relsize}
\RequirePackage{graphicx}
\newcommand{\bfsigma}{{\boldsymbol{\sigma}}}

\newcommand{\uu}{\textbf{u}}
\newcommand{\vv}{\textbf{v}}
\newcommand{\ww}{\textbf{w}}
\newcommand{\xx}{\textbf{x}}

\newcommand{\nn}{\textbf{n}}
\newcommand{\bff}{\textbf{f}}

\newcommand{\bft}{\textbf{t}}
\newcommand{\bfn}{\textbf{n}}

\definecolor{myblue1}        {RGB}{0,177,234}                        
\definecolor{mygray1}        {RGB}{76,84,93}                                

\journal{}
\begin{document}
\begin{frontmatter}
\title{
Adaptive Immersed Mesh Method (AIMM) for Fluid--Structure Interaction
}
\cortext[cor1]{~Corresponding author}
\address[]{MINES~Paris, PSL~-~Research~University, CEMEF~-~Centre~for~Material~Forming,\\ CNRS~UMR~7635, CS~10207~rue~Claude~Daunesse, 06904~Sophia-Antipolis~Cedex, France.}
\author[]{R.~Nemer}
\ead{ramy.nemer@minesparis.psl.eu}
\author[]{A.~Larcher}
\ead{aurelien.larcher@minesparis.psl.eu}
\author[]{E.~Hachem\corref{cor1}}
\ead{elie.hachem@minesparis.psl.eu}
\begin{abstract}
Our paper proposes an innovative approach for modeling Fluid-Structure Interaction (FSI).Our method combines both traditional monolithic and partitioned approaches, creating a hybrid solution that facilitates FSI. 
At each time iteration, the solid mesh is immersed within a fluid-solid mesh, all while maintaining its independent Lagrangian hyperelastic solver. 
The Eulerian mesh encompasses both the fluid and solid components and accommodates various physical phenomena.
We enhance the interaction between solid and fluid through anisotropic mesh adaptation and the Level-Set methods. 
This enables a more accurate representation of their interaction. Together, these components constitute the Adaptive Immersed Mesh Method (AIMM).
For both solvers, we utilize the Variational Multi-Scale (VMS) method, mitigating potential spurious oscillations common with piecewise linear tetrahedral elements. 
The framework operates in 3D with parallel computing capabilities.
Our method's accuracy, robustness, and capabilities are assessed through a series of 2D numerical problems. 
Furthermore, we present various three-dimensional test cases and compare their results to experimental data.
\end{abstract}
\end{frontmatter}
%
\pagenumbering{arabic} 
\setcounter{page}{1}
\input{01-Introduction}
\input{02-SolidSolver}
\input{03-FluidSolver}
\input{04-Coupling}
\input{05-NumericalTests}
\input{06-Perspectives_andconclusion}
\clearpage
\biboptions{sort&compress}
\bibliographystyle{elsarticle-num}
\bibliography{main}
\end{document}

%% file: 01-Introduction.tex
\section{Introduction}
\indent
Fluid–Structure Interaction (FSI) has gained increasing prominence due to advances in computational power, scientific progress, and the necessity for solving real-world problems through numerical simulations.
This phenomenon involves the interaction of two key components: a fluid and a solid structure.
FSI applications span a wide range of engineering fields \cite{hamdan1995fluid} such as technology \cite{hartwanger20083d}, automotive engineering \cite{hou2014coupled}, aerodynamics \cite{trimarchi1970fluid}, and biomechanics \cite{leung2006fluid}. Despite its significance, experimental studies involving FSI are often time-consuming and complex, motivating the use of numerical simulations\cite{deparis2003acceleration}\cite{chouly2006simulation}.

Structures interacting with fluids tend to deform significantly, impacting their performance. 
This interaction is fundamental and necessitates the development of efficient FSI modeling approaches. 
Various FSI formulations are employed, including Eulerian, Lagrangian, Arbitrary Lagrangian-Eulerian (ALE), Eulerian-Lagrangian, and the Fictitious Domain method \cite{garg2009numerically}\cite{41e42816959649909db90a93179d8dbf}. Each of these formulations has its unique characteristics and benefits.

Eulerian formulations, commonly referred to as the monolithic approach for FSI, treat both fluid and solid as a single entity on a single grid with a fixed reference system. 
Interface tracking methods, such as the level set method, are often used with this approach. 
Additionally, methods like the embedded boundary method \cite{doi:10.1137/040604728} and the immersed boundary method \cite{peskin_2002} have been developed to manage large structural deformations and meshing complexities. 
However, these methods may struggle with the coupling of fluid-structure stresses, which is addressed through specific techniques \cite{wang:hal-00651118}\cite{41e42816959649909db90a93179d8dbf}.

In contrast, Lagrangian formulations focus on studying particle dynamics relative to a known reference frame. 
The reference frame moves with the particle, eliminating the need to account for the advective term in equations.
While this approach simplifies the tracking of the interface, it may lead to mesh quality issues during large deformations. 
Mesh adaptation techniques are employed to mitigate these challenges, and Lagrangian formulations are often preferred in Computational Solid Mechanics (CSM) \cite{persson2009curved}.

Arbitrary Lagrangian-Eulerian (ALE) formulations were initially designed to address FSI and free surface problems \cite{donea1977lagrangian}\cite{HUERTA1988277}\cite{hughes1981lagrangian}. 
However, for elements experiencing substantial distortion, the ALE formulation may require remeshing to prevent element stretching. This limitation can be overcome by using the fixed mesh ALE formulation, which involves coupling the results from the deformed fluid mesh to a fixed Eulerian mesh using methods like the extended finite element method (XFEM) and Lagrange multipliers \cite{GERSTENBERGER20081699}\cite{wall2009advances}\cite{wall2006large}.

The Eulerian-Lagrangian formulation is a combination of the Eulerian approach for fluid and the Lagrangian approach for solid. 
It is well-suited for fluid mechanics, where the focus is on the effect of flow over a region in space. 
In this method, the fluid mesh remains unaffected by the movement of the FSI interface, which is tracked using the level set method. Lagrangian multipliers and penalty methods are employed for coupling, making it a versatile choice  \cite{legay2006eulerian}. 
However, this method has some limitations, such as using a single physical model for the solid and limited flexibility in including various behavior laws \cite{legay2006eulerian}\cite{article}\cite{zilian2008enriched}.

In this paper, we propose a novel approach for solving Fluid-Structure Interaction called the Adaptive Immersed Mesh Method (AIMM). 
AIMM is a hybrid method that combines the advantages of using separate solvers for both the fluid and solid. It employs an Eulerian fluid-solid mesh where the solid is immersed in the fluid-solid mesh through the level set method, allowing for the tracking of the FSI interface over time.
Stress and velocity boundary conditions are imposed at the FSI interface, and mesh adaptation techniques are applied to capture the boundary layer and enhance precision at the fluid-solid interface. Both the fluid and solid solvers are stabilized using the Variational Multi-Scale (VMS) method. The fluid is solved naturally on an Eulerian grid, accounting for the advective term, while the solid is modeled as an elastic structure with an independent Lagrangian grid. 
This approach is designed to handle three-dimensional simulations with parallel computing capabilities. 
In the subsequent sections, we provide details of the AIMM coupling method, the stabilized solid solver formulation, the stabilized fluid solver, the mesh adaptation technique, and present various numerical examples and benchmarks to assess the method's accuracy, robustness, and capabilities \cite{article2}.

The structure of the paper is as follows: Section 2 details the stabilized solid solver formulation, Section 3 discusses the stabilized fluid solver and mesh adaptation technique, Section 4 outlines the AIMM coupling method, and Section 5 presents numerical examples and benchmarks. The paper concludes with perspectives and conclusions in Section 6.

%% file: 02-SolidSolver.tex
\section{Solid solver}
\subsection{Lagrangian Solid Dynamics and Hyperelastic Model}
In this section, we will explore the fundamental concepts of Lagrangian Solid Dynamics and the Hyperelastic model.
Lagrangian Solid Dynamics is a framework used to model the displacement and density variations within a solid structure. 
We refer to the current and initial domains as $\Omega_s$ and $\Omega_{s_0}$, respectively.
 Both domains are open sets in $d$ dimensional real space, denoted as ${\rm I\!R}^d$, with Lipshitz boundaries.
 The boundary of the domain, denoted as $\Gamma_s$, is defined using specific boundary conditions, separating it into the Dirichlet boundary ($\partial\Omega_{s\textbf{u}}$), which defines the displacement, and the Neumann boundary ($\partial\Omega_{s\textbf{t}}$), designating the traction.

The dynamics of a solid structure are described through invertible and smooth mappings. 

\begin{equation}
\boldsymbol{\phi} : = \Omega_{s_0} \rightarrow \Omega_s = \boldsymbol{\phi} (\Omega_{s_0}), 
\end{equation}

\begin{equation}
\boldsymbol{\phi} : = \Gamma_{s_0} \rightarrow \Gamma_s = \boldsymbol{\phi} (\Gamma_{s_0}), 
\end{equation}

\begin{equation}
\mathbf{X} \rightarrow \mathbf{x} = \boldsymbol{\phi}(\mathbf{X},t) \: \forall \mathbf{X} \in \Omega_{s_0}.
\end{equation}

These mappings distinguish between the material coordinate $\mathbf{X}$, and the coordinate $\mathbf{x}$ in the total Lagrangian and updated Lagrangian frameworks, respectively. 
The displacement of a solid particle is represented by $\mathbf{u}= \mathbf{x} - \mathbf{X}$.
Additionally, the Jacobian determinant ($ J = \text{det} \mathbf{F}$) and the deformation gradient ($ \mathbf{F} =  \nabla_{\mathbf{X}} \boldsymbol{\phi}$) are key measures that play a crucial role in these dynamics.

The governing equations for solid dynamics include the equation of motion, and a constraint on the density:

\begin{equation}
\rho_s \ddot{\uu} - \nabla_{\mathbf{x}} \cdot \: \bfsigma = \bff \: \text{in} \: \Omega_s,
\end{equation}
\begin{equation}
\rho_s J = \rho_{s_0} \: \text{on} \: \partial \Omega_{s_u}.
\end{equation}

In these equations, $\rho_s$ represents the density, $ \bff $ is the source term, and $\bfsigma$ is the symmetric Cauchy stress tensor.
The equations are applicable to isotropic, compressible, and incompressible materials, thanks to the decomposition of stress into its volumetric and deviatoric constituents.

\begin{equation}
\bfsigma = p_s \mathbf{I} + dev[\bfsigma] .
\end{equation}

The Hyperelastic model is a fundamental concept used to describe material behavior.
It involves a Helmholtz free energy function, denoted as $\Psi (\textbf{C})$,where $\textbf{C}$ is the right Cauchy-Green strain tensor ($\textbf{C} = \textbf{F}^T\textbf{F}$).
The second Piola--Kirchhoff stress tensor given by $\textbf{S} =  J \textbf{F}^{-1} \bfsigma \textbf{F}^{-T}$ is also obtained by deriving the Helmholtz free energy functional  $\boldsymbol{\Psi} (\textbf{C})$ with respect to $\textbf{C}$:

\begin{equation}
\textbf{S}=2\partial_{\textbf{C}} \Psi (\textbf{C}).
\end{equation}

This energy function is decomposed into volumetric and deviatoric parts.

\begin{equation}
\Psi ( \textbf{C} ) = U(J) + W( \bar{\textbf{C}} ).
\end{equation}

Where  $J = \sqrt{det \textbf{C}}$, and $\bar{\textbf{C}} = J^{-\frac{2}{3}} \textbf{C}$ is the volumetric/deviatoric part of $\textbf{C}$.

In specific volumetric models, such as the Neo-Hookean and Simo--Taylor models:

\begin{equation}
U( J ) = \frac{1}{4} \kappa ( J^2 - 1 ) - \frac{1}{2} \kappa \text{ln} J,
\end{equation} 

\begin{equation}
W ( \bar{\textbf{C}} ) = \frac{1}{2} \mu_s ( \text{tr}\bar{\textbf{C}} - 3 )  =\frac{1}{2} \mu_s ( \bar{\mathbf{I}_1} - 3 ).
\end{equation}

where $\kappa$ and $\mu_s$ are material properties,and ${ \mathbf{I}_1}=\text{tr}\bar{\textbf{C}}$.
The stress can also be split to its deviatoric and volumetric part:

\begin{equation}
p_s = 2 J^{-1} \textbf{F} \frac{\partial U(J)}{\partial \textbf{C}} \textbf{F}^T = U'(J) = \frac{1}{2} \kappa (J + J^{-1}),
\end{equation}

\begin{equation}
\textbf{dev}[\bfsigma] = 2 J^{-1} \textbf{F} \frac{\partial W(\bar{\textbf{C}})}{\partial \textbf{C}} \textbf{F}^T = \mu_s J^{-\frac{5}{3}} \textbf{dev} [ \textbf{FF}^T].
\end{equation}

We can also write:

\begin{equation}
\mathbf{F} \mathbf{F^T} =  \nabla_{\mathbf{X}} \uu +  \nabla_{\mathbf{X}}^T \uu + \nabla_{\mathbf{X}} \uu  \nabla_{\mathbf{X}}^T \uu + \mathbf{I}.
\end{equation}

A linearization of the above equation is considered, by starting with:

\begin{equation}
\nabla_{\mathbf{X}} \uu = (  \mathbf{I} - \nabla_{\mathbf{x}} \uu ) ^{-1} -  \mathbf{I}.
\end{equation}

Assuming a small variation in displacement denoted by $\delta \uu$, and
recalling that for very small displacement, $(  \mathbf{I} - \nabla \uu ) ^{-1} =  \mathbf{I} + \nabla \uu$, we end up with the following, as explained in \cite{NEMER2021113923}

\begin{equation}
\begin{split}
\mathbf{F} \mathbf{F^T} =  & ( \mathbf{I} - \nabla_{\mathbf{x}} \uu ) ^{-1} - \mathbf{I} + ( ( \mathbf{I} - \nabla_{\mathbf{x}} \uu ) ^{-1} - \mathbf{I})^T +  (( \mathbf{I} - \nabla_{\mathbf{x}} \uu ) ^{-1} - \mathbf{I})( ( \mathbf{I} - \nabla_{\mathbf{x}} \uu ) ^{-1} - \mathbf{I})^T  \\
& + \mathbf{I} + 2\epsilon (\delta \uu) + \nabla_{\mathbf{x}} \delta \uu (\nabla_{\mathbf{x}} \delta \uu)^T + \nabla_{\mathbf{x}} \delta \uu (\nabla_{\mathbf{x}}  \uu)^T + (\nabla_{\mathbf{x}}  \uu) ( \nabla_{\mathbf{x}} \delta \uu)^T.
\end{split}
\end{equation}

The final system of equations to be solved now is given by:

\begin{equation}
\rho_s \ddot{\uu} - \nabla_{\mathbf{x}} p_s - \nabla_{\mathbf{x}} \cdot \: dev[\bfsigma] = \bff \: \text{in} \: \Omega_s,
\end{equation}

\begin{equation}
  \nabla_{\mathbf{x}} \cdot \uu - \frac{1}{\kappa} p_s  = g \: \text{in} \: \Omega_s,
\end{equation}

\begin{equation}
\uu=\textbf{l} \: \text{on} \: \partial \Omega_{s_u},
\end{equation}

\begin{equation}
\bfsigma\bfn=\bft \: \text{on} \: \partial \Omega_{s_t},
\end{equation}

\begin{equation}
\rho_s J = \rho_{s_0}.
\end{equation}

The weak discrete form of the above equation in Galerkin finite elements is given by:

\begin{equation}
(\rho_s \frac{\partial^2 \uu_h}{\partial t^2},\ww_h) + a'(\uu_h, \ww_h)+(p_{s_h}, \nabla_{\mathbf{x}} .\ww_h)  =L(\ww_h) \: \forall \: \ww_h \in W_{h,0},
\end{equation}
\begin{equation}
(\nabla_{\mathbf{x}} \cdot \uu_h, q_h) - (\frac{1}{\kappa} p_{s_h}, q_h)=(g, q_h) \: \forall \: q_h \in Q_h,
\end{equation} 
where $a'$ is given by:
\begin{equation}
a'(\uu_h, \ww_h)= \int_{\Omega_s} \mu_s \: dev[\boldsymbol{\sigma}] : \nabla_{\mathbf{x}}^s \ww.
\end{equation}

\subsection{Variational Multi Scale stabilization}

In this section, we discuss Variational Multi-Scale (VMS) stabilization, a technique used to ensure the stability and accuracy of numerical simulations. 
VMS stabilization is particularly valuable in problems involving the simulation of fluids and solids and is aimed at improving the performance of finite element methods by addressing issues such as the inf-sup condition or Babuska--Brezzi\cite{babuvska1971error}.

The inf-sup stability  condition coerces the interpolation relationship between the variables.
This leads to different interpolation order for $\uu$ and $p_s$.
Same order interpolation exhibits weak numerical performance since it does not respect the inf-sup condition.
Different types of stabilization exist, that help alleviate this problem.
P1/P1 elements are used in our case for the displacement and pressure variables, with a Variational Multi-Scale Method (VMS) stabilization.
This allows us have same order interpolation for both variables.
In \cite{HUGHES198685}, equal order elements for velocity and pressure were utilized to solve the Stokes equations. 
The author also proved the convergence and stability of this method.
This work inspired the extension of the VMS method to solve the Navier--Stokes equations \cite{HACHEM20108643}.
The linear elastic equations were also tackled in \cite{HUGHES198785, 10.1007/BF01395881} using VMS.
It provides a natural stabilization using an orthogonal decomposition of the solution spaces.
The function spaces are first decomposed into their coarse and fine scale components, which yields:

\begin{equation}
W_0 = W_{h,0} + W_0',
\end{equation}
\begin{equation}
W = W_{h} + W',
\end{equation}
\begin{equation}
Q = Q_h + Q'.
\end{equation}

Following \cite{hughes1998variational}, the displacement and pressure are decomposed:

\begin{equation}
\uu=\uu_h+\uu',
\end{equation}
\begin{equation}
p_s=p_{s_h}+p_s'.
\end{equation}

The same decomposition is also applied for the test functions:

\begin{equation}
\ww=\ww_h+\ww',
\end{equation}
\begin{equation}
q=q_h+q'.
\end{equation}

The transient mixed finite element approximation of equations (21)(22)(23):
\\
Coarse Scale
\begin{equation}
(\rho_s \frac{\partial^2 (\uu + \uu')}{\partial t^2},\ww_h) + a'((\uu_h+\uu'), \ww_h)+(p_{s_h}+p_s', \nabla_{\mathbf{x}} \cdot \ww_h)  =L(\ww_h) \: \forall \: \ww_h \in W_{h,0} \subset [H_0^1]^d,
\end{equation}
\begin{equation}
(\nabla_{\mathbf{x}} \cdot (\uu_h+\uu'), q_h) - (\frac{1}{\kappa} (p_{s_h}+p_s'), q_h)=(g, q_h) \: \forall \: q_h \in Q_h \subset L_{\int=0}^2,
\end{equation} 
Fine Scale
\begin{equation}
(\rho_s \frac{\partial^2 (\uu + \uu')}{\partial t^2},\ww') + a'((\uu_h+\uu'), \ww')+(p_{s_h}+p_s', \nabla_{\mathbf{x}} .\ww')  =L(\ww') \: \forall \: \ww' \in W',
\end{equation}
\begin{equation}
(\nabla_{\mathbf{x}} \cdot (\uu_h+\uu'), q') - (\frac{1}{\kappa} (p_{s_h}+p_s'), q')=(g, q') \: \forall \: q' \in Q'.
\end{equation}

The fine-scale problem is first solved and modeled implicitly in terms of the time-dependent coarse-scale problem, thus capturing the sub-scales' behavior.
An elaboration of this choice is found in \cite{dubois1999dynamic}.
Sub-scales are not explicitly tracked in time but are considered quasi time-dependent since they respond to the time-dependent large-scale residual. 
For more information on time tracked sub-scales, please refer to \cite{https://doi.org/10.1002/fld.1481}.

The fine scale approximations are given by:

\begin{equation}
\uu' = (\tau_{\uu} P_{\uu}' (R_{\uu})),
\end{equation}
\begin{equation}
p_s' = (\tau_c P_{c}' (R_c)).
\end{equation}

Where $R_{\uu}$ and $R_c$ are the finite element residuals, $P_{\uu}'$ and $P_{c}'$ are the projection operators, and $\tau_{\uu}$ and $\tau_c$ are tuning parameters.
Note that in this current work, both $P_{\uu}'$, and $P_{c}'$ are taken as the identity matrix $\mathbf{I}$.

Afterwards, the coarse scale equations are given by:

\begin{equation}
(\rho_s \frac{\partial^2 \uu}{\partial t^2},\ww_h) + a'(\uu_h, \ww_h)+(p_{s_h}, \nabla_{\mathbf{x}} \cdot\ww_h)+ (p_s', \nabla_{\mathbf{x}} \cdot\ww_h) =L(\ww_h) \: \forall \: \ww_h \in W_{h,0},
\end{equation}

\begin{equation}
(\nabla_{\mathbf{x}} \cdot \uu_h, q_h) - (\frac{1}{\kappa} p_{s_h}, q_h) - (\frac{1}{\kappa} p_s', q_h) - (u',\nabla_{\mathbf{x}} q_h) =(g, q_h) \: \forall \: q_h \in Q_h.
\end{equation} 

The finite element residuals are given by

\begin{equation}
R_{\uu}=\bff-\rho_s\ddot{\uu_h}+\nabla_{\mathbf{x}} p_{s_h}+\nabla_{\mathbf{x}}\cdot\:dev[\bfsigma],
\end{equation} 
\begin{equation}
R_c=g-\nabla_{\mathbf{x}}\cdot\uu_h+\frac{1}{\kappa}p_{s_h}.
\end{equation}

Modeling the fine scales as in (35)(36), we finally get

\begin{equation}
(\rho_s \frac{\partial^2 (\uu)}{\partial t^2},\ww_h) + a'(\uu_h), \ww_h)+(p_{s_h}, \nabla \cdot\ww_h)+ \sum_{K \in T_h}(\tau_{c}(g - \nabla_{\mathbf{x}} \cdot \uu_h + \frac{1}{\kappa} p_{s_h}), \nabla .\ww_h) =L(\ww_h) \: \forall \: \ww_h \in W_{h,0},
\end{equation}

\begin{equation}
\begin{split}
(\nabla . (\uu_h), q_h) - (\frac{1}{\kappa} (p_{s_h}), q_h)+ \sum_{K \in T_h} (\frac{\tau_{c}}{\kappa} (\nabla_{\mathbf{x}} \cdot \uu_h - \frac{1}{\kappa} p_{s_h} - g), q_h)
\\ 
+ \sum_{K \in T_h} (\tau_{\uu} ( \rho_s \ddot{\uu_h} - \nabla_{\mathbf{x}} p_{s_h} - \nabla_{\mathbf{x}} \cdot\:dev[\bfsigma] -\bff),\nabla q_h) = 
(g, q_h) \: \forall \: q_h \in Q_h.
\end{split}
\end{equation}

Stabilization parameters, such as $\tau_{\uu}$ and $\tau_{c}$, are computed within each element to determine the impact of sub-scales. 
These parameters contribute additional terms to the finite element equations, addressing problems like spurious pressure oscillations and improving overall accuracy and stability in simulations.
The general definition of the stabilization parameters, $\tau_{\uu}$ and $\tau_{c}$, considers factors like the density, time step size, element size, and material properties. 
These parameters play a critical role in shaping the stability and performance of the finite element simulations.
The general definition of the stabilization parameters \cite{HACHEM2016238}\cite{CODINA20024295}, computed within each element gives:

\begin{equation}
\tau_{\uu} = (( \frac{ \rho_s}{(c_0\Delta t)^2} )^2+( \frac{ 2\mu_s}{c_1 h^2_K} )^2)^{-\frac{1}{2}},
\end{equation}

\begin{equation}
\tau_{c} = ((2c_2\mu_s )^2)^{-\frac{1}{2}}.
\end{equation}

Where $h_k$ is the characteristic length of the element, and $c_0$,$c_1$, and $c_2$ are constants to be determined.
These parameters are computed within each element to account for local variations and provide stabilization.

In summary, Variational Multi-Scale (VMS) stabilization enhances the stability and accuracy of finite element simulations, particularly in problems involving fluids and solids. 
It achieves this by addressing issues like the inf-sup condition and modeling fine scales to capture sub-scale behavior, contributing to the overall performance and reliability of simulations. 
Stabilization parameters, such as $\tau_{\mathbf{u}}$ and $\tau_c$, are computed based on various factors to ensure the effectiveness of the method.

%% file: 03-FluidSolver.tex
\section{Fluid solver}
\subsection{Newtonian incompressible equations}

Let $\Omega_f \subset \mathbb{R}^d$ be the spatial domain at time $t \in [0, T]$, where $d$ is the spatial dimension, and $\Gamma_f$ is the boundary of $\Omega_f$. 
The mixed formulation for the transient incompressible Navier-Stokes equations, expressed in terms of velocity and pressure, is as follows:

\begin{equation}
    (\rho_f \frac{\partial \vv}{\partial t} + (\vv \cdot \nabla) \vv) - \nabla \cdot \bfsigma = \bff \: \text{in} \: \Omega_f,
\end{equation}
\begin{equation}
    \nabla \cdot \vv = 0 \: \text{in} \: \Omega_f.
\end{equation}

In these equations, $\rho_f$ and $\mathbf{v}$ represent the fluid density and velocity, while $\mathbf{f}$ is the source term. The stress tensor $\boldsymbol{\sigma}$ is defined as:

\begin{equation}
    \bfsigma = 2\mu_f \epsilon(\vv) - p_f \mathbf{I},
\end{equation}

with the strain rate tensor $\epsilon(\mathbf{v})$ defined as:

\begin{equation}
    \epsilon(\vv) = \frac{1}{2}(\nabla \vv + \nabla^T \vv).
\end{equation}

The boundary conditions for this problem are as follows:

\begin{equation}
    \vv = \textbf{m} \: \text{on} \: \partial \Omega_{f_u},
\end{equation}
\begin{equation}
   \bfsigma\bfn=\bft \: \text{on} \: \partial \Omega_{f_t}.
\end{equation}

Here, $\textbf{m}$ is a known imposed value, and $\Omega_{f_u}$ and $\Omega_{f_t}$ represent the boundaries where Dirichlet and Neumann boundary conditions are applied, respectively.
The weak discrete form of these equations in Galerkin finite elements is given by:

\begin{equation}
(\rho_f\frac{\partial \vv_h}{\partial t},\ww_h) + (\rho_f(\vv_h \cdot \nabla) \vv_h, \ww_h) + (2\mu \epsilon (\vv_h):\epsilon(\ww_h)) - (p_{f_h},\nabla \cdot \ww_h) = (\bff,\ww_h) \: \forall \: \ww_h \in W_{h,0} \subset [H_0^1]^d,
\end{equation}
\begin{equation}
(\nabla \cdot \vv_h, q_h) = 0 \: \forall \: q_h \in Q_h \subset L_{\int=0}^2.
\end{equation} 

\subsection{Variational Multi Scale (VMS) stabilization}

Applying the principles of VMS, as previously employed in the solid solver, leads to the following formulation \cite{HUGHES198685, HACHEM20108643}:

\begin{equation}
\begin{split}
& (\rho_f\frac{\partial \vv_h}{\partial t},\ww_h) + (\rho_f(\vv_h \cdot \nabla) \vv_h, \ww_h) + (2\mu \epsilon (\vv_h):\epsilon(\ww_h)) - (p_{f_h},\nabla \cdot \ww_h) + (\nabla \cdot \vv_h, q_h) 
\\ 
&  - (\bff,\ww_h) + \sum_{K \in T_h} (\tau_{\uu,K} ((\rho_f \frac{\partial \vv_h}{\partial t} + (\vv_h \cdot \nabla) \vv_h) + \nabla p_{f_h} -\bff),\rho_f \vv_h \nabla \ww_h)_K
\\ 
& + \sum_{K \in T_h} (\tau_{\uu,K} ((\rho_f\frac{\partial \vv_h}{\partial t} + (\vv_h \cdot \nabla) \vv_h) + \nabla p_{f_h} -\bff),\nabla q_h)_K
\\
& + \sum_{K \in T_h} (\tau_{c,K} \nabla \cdot \vv_h,\nabla \cdot \ww_h)_K=0  \: \forall \: \ww_h \in W_{h,0},
\: \forall \: q_h \in Q_h.
\end{split}
\end{equation} 

Compared to the standard Galerkin formulation, this VMS formulation includes additional integrals that are evaluated element-wise. 
These terms account for the effects of sub-scales and are developed consistently to address instabilities in advective-dominated regimes. 
Additionally, an extra term modeling small-scale pressure helps stabilize high Reynolds number flows. 
An implicit linearization of the advective term is incorporated using a Newton-Raphson linearization method. 
For further details, refer to \cite{HACHEM20108643}.

%% file: 04-Coupling.tex
\section{Coupling}
In the context of Fluid-Structure Interaction (FSI), the choice of coupling strategy is pivotal for the accurate simulation of the complex interplay between fluid and solid domains. 
Two primary approaches exist in the FSI literature: the monolithic approach and the partitioned approach. 
Understanding the fundamental differences between these approaches is essential for selecting the most suitable method for a given FSI problem. 
The monolithic approach involves solving the governing equations for both the fluid and solid components on a single mesh. 
This unified framework directly reflects the strong coupling between the fluid and solid domains, as they share the same computational grid. Consequently, information is exchanged between fluid and solid regions seamlessly within a single numerical framework. 
In contrast, the partitioned approach decouples the fluid and solid simulations, typically employing separate meshes for each domain. This decoupling allows for a degree of independence in solving the fluid and solid equations. 
However, the partitioned approach can further be categorized into two sub-types:

\begin{itemize}
    \item Two-Way Strong Coupling: In this configuration, information is exchanged bidirectionally between the fluid and solid domains, ensuring that changes in one domain influence the other. This approach maintains a high level of coupling, but it comes with computational challenges due to the simultaneous interaction of both domains.
    \item One-Way Weak Coupling: Here, the coupling is unidirectional, with one domain influencing the other while remaining relatively unaffected by the feedback. This approach simplifies the computational burden but may lead to some loss of accuracy in modeling the mutual impact of fluid and solid components.
\end{itemize}
    
\begin{figure}[H]
\centering
   \includegraphics[width=1\linewidth]{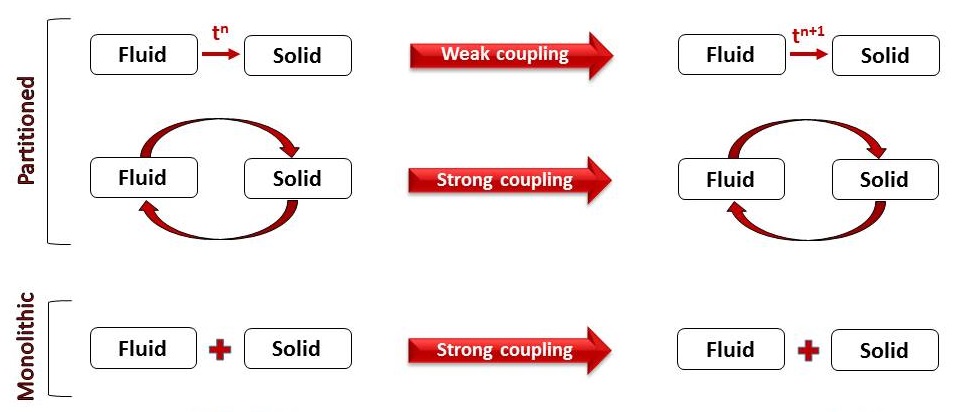}
\caption{Different approaches for FSI.}
\label{fig:FigureBis1}
\end{figure}

Figure \ref{fig:FigureBis1} visually illustrates the differences between these coupling approaches. The choice between monolithic and partitioned FSI strategies should align with the specific requirements and characteristics of the problem at hand.

The Adaptive Immersed Mesh Method (AIMM) introduces a novel hybrid approach, designed to harness the advantages of both monolithic and partitioned methods while mitigating their inherent limitations. The AIMM technique not only maintains the flexibility of separate solvers for the fluid and solid domains but also ensures robust, stable, and strong coupling, akin to the monolithic approach. 
To achieve this, AIMM immerses the solid mesh within a fluid-solid mesh at each time step, creating a shared computational domain where both fluid and solid entities coexist. 
The signed distance function (level set) is then employed to track the fluid-structure interface, providing a foundation for seamless interaction. 
Furthermore, a Moving Mesh Method (MMM) ensures that the interface remains accurately represented within the solid mesh.

To achieve the high precision necessary for FSI simulations, AIMM incorporates a mesh refinement technique, which selectively enhances mesh resolution in an anisotropic manner along the fluid-structure interface. 
This innovative approach combines the robustness of monolithic FSI with the advantages of separate domain simulations. 
Figure \ref{fig:FigureBis2} provides an illustrative representation of the AIMM hybrid method, while Figure \ref{fig:Figure22} offers an overview of the AIMM coupling framework for FSI. 
The following sub-sections will explain the different components of the AIMM coupling method:

\begin{figure}[H]
\centering
   \includegraphics[width=1\linewidth]{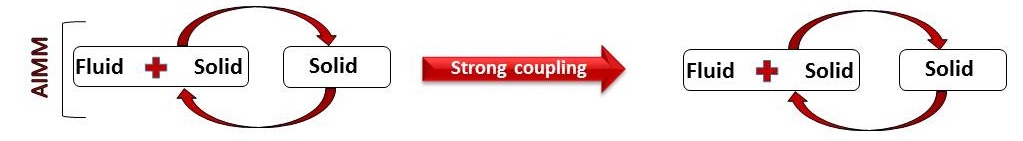}
\caption{Different approaches for FSI.}
\label{fig:FigureBis2}
\end{figure}

\subsection{Level set approach}

The AIMM method employs the signed distance function of the FSI interface, denoted as $\Gamma_{FSI}$, to define the contact surface between the fluid and solid components. This approach uses the level set function, $\alpha_{FSI}$, which represents the signed distance from any point in the fluid-solid mesh to $\Gamma_{FSI}$. Thus, the interface $\Gamma_{FSI}$ is defined as the iso-zero of the level set function:

\begin{equation}
    \left\{
    \begin{array}{ll}
        \alpha_{FSI}(\xx)= \pm d(\xx, \Gamma_{FSI}), \: \xx \: \in \: \Omega \\
        \Gamma_{FSI} = \left\{ \xx, \alpha_{FSI}(\xx)=0 \right\}
    \end{array}
    \right.
\end{equation}

Readers interested in further algorithm details can refer to \cite{bruchon:emse-00475556}. Additionally, alternative smoother functions may be considered for regions away from $\Gamma_{FSI}$ \cite{article1}.

\begin{figure}[H]
\centering
   \includegraphics[width=0.7\linewidth]{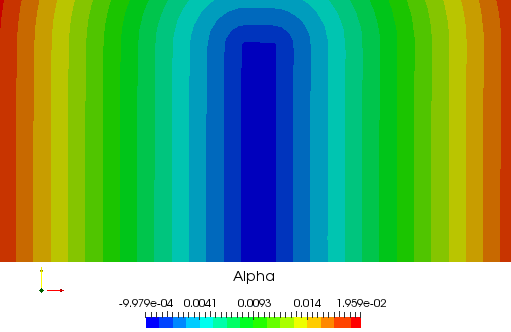}
\caption{The signed distance function $\alpha$ example of an immersed circle.}
\label{fig:LevelSet}
\end{figure}

\begin{figure}[H]
\centering
   \includegraphics[width=0.7\linewidth]{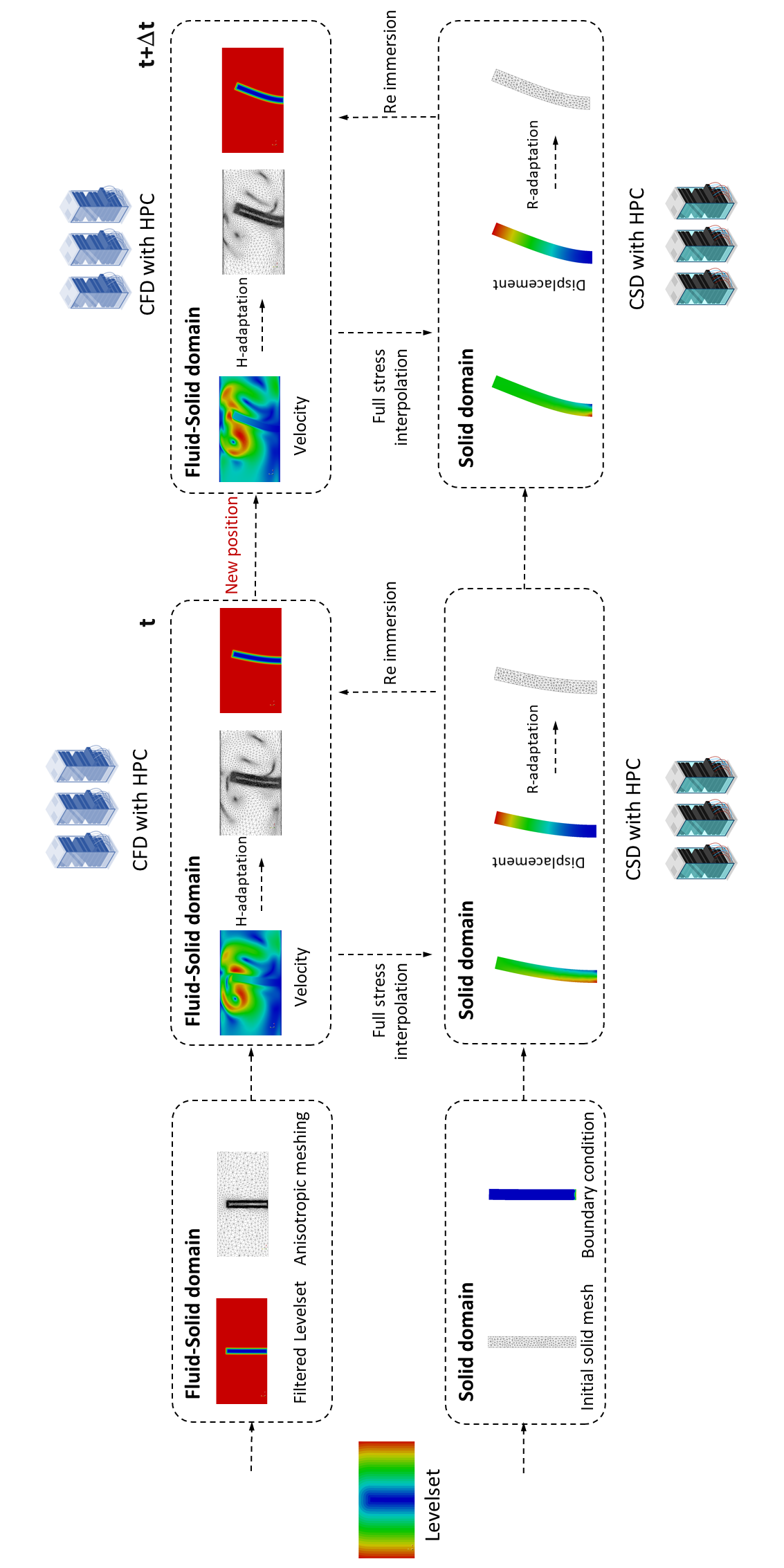}
\caption{Two way coupling loop.}
\label{fig:Figure22}
\end{figure}

\subsection{Physical continuity}

Ensuring two-way coupling between the fluid and solid on $\Gamma_{FSI}$ requires maintaining velocity and stress continuity at the interface:

\begin{equation}
    \vv_s=\vv_f \: \text{on} \: \Gamma_{FSI},
\end{equation}
\begin{equation}
    \bfsigma_s \nn = \bfsigma_f \nn \: \text{on} \: \Gamma_{FSI}.
\end{equation}

\subsection{Edge-based mesh adaptation}

The AIMM methodology incorporates anisotropic mesh adaptation for unstructured meshes, a key factor in achieving high accuracy and reducing computational time \cite{formaggia2003anisotropic, hoffman2003adaptive, https://doi.org/10.1002/nme.4481}.
Anisotropic meshes prove superior in capturing smaller features than octree-based methods, while simultaneously reducing the number of integration points \cite{Khalloufi201644, Legrain2018}.
This approach concentrates elements in regions of interest characterized by significant variations in variables or their gradients. 
For FSI, this enables accurate modeling of the fluid-solid interface where rapid variable changes can occur.

The adaptive process relies on a monitor function, either scalar or vector, to control the shape, size, and orientation of mesh elements. 
This monitor function estimates the solution error and distributes it equitably over each mesh element. 
The mesh adaptation algorithm generates a mesh and computes a numerical solution at each time step, evaluating an interpolation error estimate. 
Subsequently, a minimization problem aims to reduce the interpolation error in the L1-norm, independently of the current problem \cite{Coupez11}. 
An optimal metric is derived to monitor solution development and minimize the interpolation error. 
A new mesh is generated, aligning with the metric field. 
This approach simplifies the computation of the metric and its associated edge-based error.

\subsubsection{Edge-based error estimation}
\begin{figure}[h!]
\centering
\includegraphics[scale=0.20]{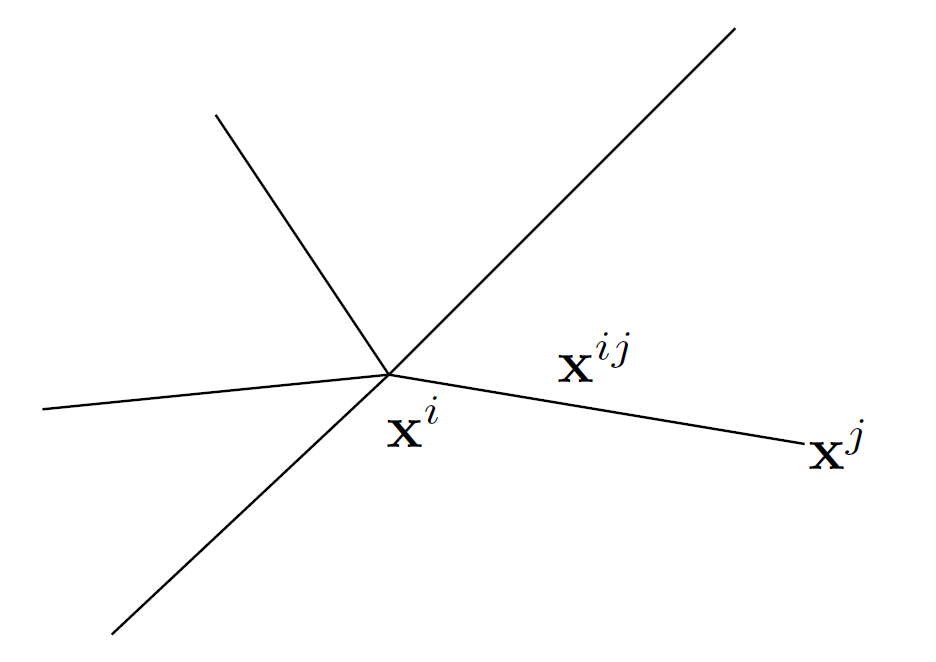}
\caption{Patch associated with node $x^i$}
\label{fpatch}
\end{figure}

Let $u_h$ be a a first order finite element approximation, acquired through the Lagrange interpolation operator of a regular function $u \in C^2 (\Omega)$.
At each vertex $i$ of the mesh, we have $U_i = u(x^i) = u_h (x^i)$ (where $x^i$ are the coordinates of the vertex $i$).
Let $\Gamma (i)$ be the "patch" associated to a vertex $x_i$ of the mesh defined as the set of nodes which share one edge with $x_i$, and let us denote by $x_{ij}$ the edge connecting $x_i$ to $x_j$ as in Figure \ref{fpatch}.

\noindent The continuity of the gradient $\nabla u^h \cdot x^{ij}$ on the edge $x_{ij}$ enables us to write

\begin{eqnarray}
U^j = U^i + \nabla u^h \cdot x^{ij},
\end{eqnarray}

This leads to

\begin{eqnarray}
\nabla u_h \cdot x^{ij} = U^j - U^i \,.
\end{eqnarray}

Following the works from \cite{Coupez11} we can write the following error estimator

\begin{eqnarray}
\mid \mid \nabla u^h \cdot x^{ij} - \nabla u (x^i) \cdot x^{ij}\mid \mid \leq \max_{y \in \mid x^i , x^j \mid} \mid x^{ij} \cdot H_{u} (y) \cdot x^{ij} \mid ,
\end{eqnarray}

with $H_{u}$ being the Hessian of $u$.
At the node $x^i$, we seek the recovered gradient $g^i$ of $u_h$

\begin{eqnarray}
\nabla g_h \cdot x^{ij} = g^j - g^i\,.
\end{eqnarray}

We are interested in the projection of the Hessian based on the gradient at the edge extremities, thus

\begin{eqnarray}
(\nabla g_h \cdot x^{ij})\cdot  x^{ij}  = (g^j - g^i)\cdot  x^{ij},\\
(H_{u}\cdot x^{ij})\cdot  x^{ij}  = g^{ij}\cdot  x^{ij},
\end{eqnarray}

with $g^{ij} = g^j - g^i$.
It can be shown in \cite{Coupez11} that the quantity $\mid g^{ij} \cdot x^{ij} \mid$ gives a second order accurate approximation of the second derivative of $u$ along the edge $x^{ij}$.
Motivated by the fact that, for first order finite elements on anisotropic meshes, edge residuals dominate a posteriori errors\cite{kunert2000edge}, it is therefore suitable to define an error indicator function associated to the edge $x^{ij}$ as

\begin{eqnarray}
e^{ij} = \mid g^{ij} \cdot x^{ij} \mid \,.
\end{eqnarray}

And this error, is the exact interpolation error along the edge and allows to evaluate the global $L_1$ error.
However, the gradient is not know at the vertices, thus a recovery procedure has to be considered.

\subsubsection{Gradient recovery procedure}

The gradient recovery procedure relies on the following optimization problem

\begin{eqnarray}
G^i = arg \min_{G} \left(\sum_{j \in \Gamma (i)} \mid (G - \nabla u_h) \cdot  x^{ij} \mid^2 \right),
\end{eqnarray}

where $G^i$ is the recovered gradient.
Denoting by $\otimes$ the tensor product between two vectors, let us introduce $X^i$ the length distribution tensor at node $i$

\begin{eqnarray}
X^i = \frac{1}{\mid \Gamma (i) \mid} \left( \sum_{j \in \Gamma (i)} x^{ij} \otimes x^{ij} \right),
\end{eqnarray}

this gives us an average representation of the distribution of edges in the patch.
Let us express the recovered gradient $G^i$ in terms of the length distribution tensor

\begin{eqnarray}
G^i = (X^i)^{-1} \sum_{j \in \Gamma (i) } U^{ij} x^{ij}\,.
\end{eqnarray}

Therefore, the estimated error $e_{ij}$ is thus written as

\begin{eqnarray}
\label{eq}
e_{ij} = G^{ij}.x^{ij}\,.
\end{eqnarray}

\subsubsection{Metric construction}

A stretching factor $s^{ij}$ defined as the ratio between the length of the edges $x^{ij}$ before and after the adaptation procedure is introduced in order to correlate the error indicator defined in (\ref{eq}) to the associated metric \cite{Coupez11}.
We end up with the following expression for the metric

\begin{eqnarray}
\widetilde {M}^{i} = (\widetilde{X}^i)^{-1},
\end{eqnarray}

where $\widetilde{X}^i$ is defined as

\begin{eqnarray}
\widetilde{X}^i = \frac{1}{\mid \Gamma (i) \mid} \left( \sum_{j \in \Gamma (i)} s^{ij} \otimes s^{ij} \right)\,.
\end{eqnarray}

\begin{figure}[H]
\centering
   \includegraphics[width=1\linewidth]{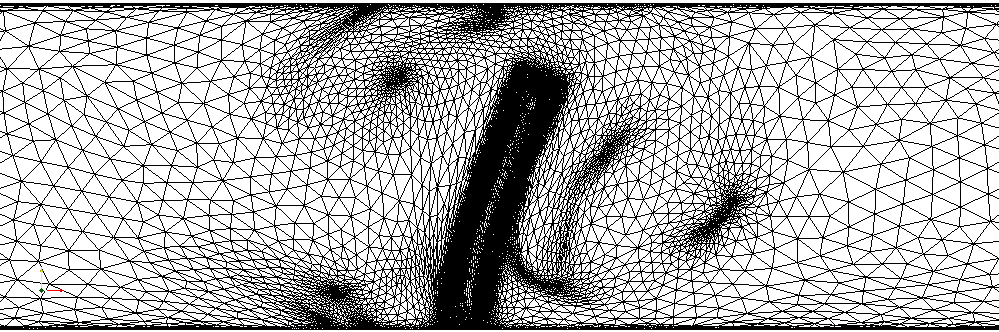}
\caption{An example of the mesh adaptation on multiple criteria.}
\label{fig:MeshAdap2}
\end{figure}

The stretching factor $s^{ij}$ of the edge $ij$ is chosen so that the total number of nodes in the mesh is kept fixed and is defined as

\begin{eqnarray}
s^{ij} = \left(\frac{e_{ij}}{e(N)}\right),
\end{eqnarray}

where $e(N)$ the total error.
An example of the mesh adaptation is shown in Figure \ref{fig:MeshAdap2}.

%% file: 05-NumericalTests.tex
\section{Numerical Validation}
In the following subsections, we present a series of test cases designed to rigorously validate the numerical methods proposed in this study. 
These test cases provide a comprehensive assessment of the model's accuracy and capabilities. 
Beginning with 2D scenarios, we progressively transition into 3D simulations to assess the model's performance under varying conditions. 
Through these test cases, we aim to thoroughly evaluate the effectiveness of our approach in simulating and analyzing complex phenomena.

\subsection{Bending beam 1}
This test examines a 2D semi-stationary beam bending problem, originally proposed in \cite{Baiges}\cite{https://doi.org/10.1002/nme.6321}. 
In this problem, a clamped plate is positioned parallel to the direction of the flow, and the aerodynamic forces exerted by the fluid induce bending in the plate. 
As a result, the plate bending must be accurately accounted for in the fluid flow simulations. 
Over time, the plate reaches a stationary position, where it no longer exhibits oscillations and stabilizes in its steady-state equilibrium.

Two variations of this problem are considered, where the density of the solid material differs, significantly impacting the dynamics of the solid-structure response. The fluid and solid properties, as summarized in Table \ref{table:1}, play a crucial role in determining the system's behavior.

\begin{table}[ht!]
\centering
\begin{tabular}{|p{3cm}|p{3cm}||p{3cm}|p{3cm}|}
 \hline
 \multicolumn{2}{|c|}{Fluid} & \multicolumn{2}{|c|}{Solid}\\
 \hline
$\rho_f$ & 2.0 & $\rho_s$ & 1.0 or 10.0\\
$\mu_f$ & 0.2 & $\mu_s$ & 5000\\
&&$\lambda_s$& 2000\\
Model & Newtonian & Model &  Neo-Hookean\\
\hline
\end{tabular} 
\caption{Fluid and solid properties for the bending beam 1 problem.}
\label{table:1}
\end{table}

Here, $\rho_f$ represents the fluid density, $\mu_f$ is the dynamic viscosity of the fluid, and $\rho_s$ and $\mu_s$ denote the solid density and shear modulus, respectively. Additionally, $\lambda_s$ represents the Lamé parameter for the solid, and the fluid is modeled as Newtonian, while the solid follows a Neo-Hookean model.

\begin{figure}[H]
\centering
\begin{subfigure}{1\textwidth}
   \includegraphics[width=1\linewidth]{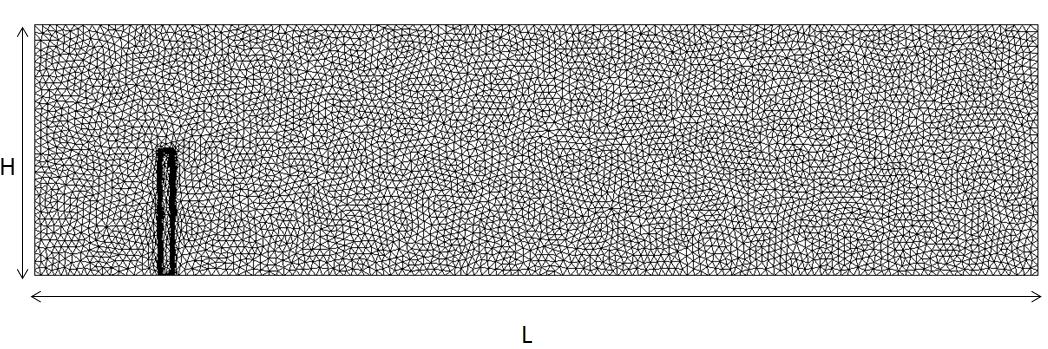}
    \caption{Fluid-Solid Mesh} \label{fig:1a}
  \end{subfigure}%
  \\
  \begin{subfigure}{0.4\textwidth}
  \centering
   \includegraphics[width=1\linewidth]{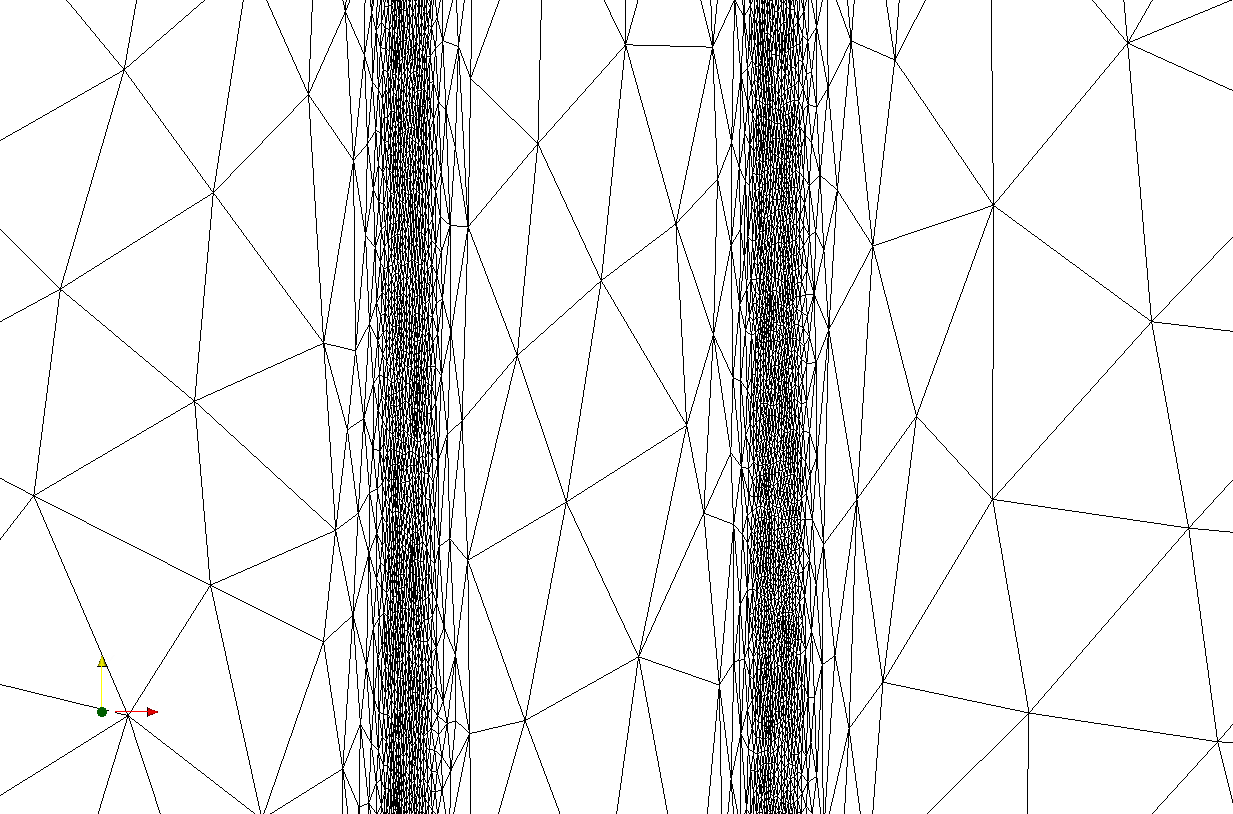}
    \caption{Anisotropic mesh at the interface} \label{fig:1b}
  \end{subfigure}
  \begin{subfigure}{0.15\textwidth}
  \centering
   \includegraphics[width=1\linewidth]{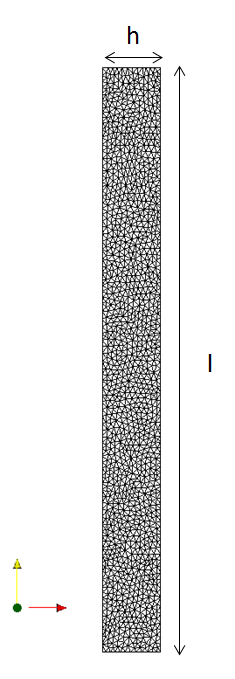}
    \caption{Solid Mesh} \label{fig:1c}
  \end{subfigure}
\caption{Problem set up and mesh}
\label{fig:ProblemDescription}
\end{figure}

The geometrical setup of the problem and the associated meshes are visualized in Figure \ref{fig:ProblemDescription}. 
The parameters used in this setup are as follows: $L = 80$, $H = 20$, $l = 10$, and $h = 1$. 
Figure \ref{fig:1a} displays the fluid-solid mesh, which exhibits higher mesh resolution at the interface, where the fluid and solid domains interact. 
For a closer inspection of this interface region, Figure \ref{fig:1b} provides a magnified view of the fluid-solid mesh, highlighting the presence of anisotropic stretched elements. Lastly, Figure \ref{fig:1c} showcases the mesh designed for the solid domain.

To simulate the flow in this problem, an inlet velocity of $1$ in the $x$ direction is imposed at the inlet boundary. 
The outlet is kept free, allowing for fluid outflow. 
The top and bottom boundaries of the domain are free to move tangentially. 
Specifically, the solid's velocity is enforced at the fluid-structure interface, and zero-displacement boundary conditions are applied to the bottom of the solid plate. 
Moreover, fluid stress (traction forces) is imposed on the fluid-structure interface within the solid domain.

The two following figures, Figure \ref{fig:Codina} and Figure \ref{fig:Codina2}, display the x-displacement of the tip of the solid plate versus time for two different scenarios: one with a solid density of 10 and the other with a solid density of 1. 
These figures clearly illustrate how variations in solid density significantly influence the dynamics of the solid structure, leading to differences in the time required to reach a steady-state solution.
These insightful simulations and their interpretations contribute to our understanding of how variations in solid properties influence the dynamics of fluid-structure interaction in the Bending Beam 1 problem.

\begin{figure}[H]
        \centering
        \begin{tikzpicture}[trim axis left, trim axis right]
                \begin{axis}[
                scale=0.8,
                transform shape,
                width=.7\textwidth,
                height=.5\textwidth,
                label style={font=\scriptsize},
                tick label style={font=\scriptsize},
                legend style={font=\scriptsize},
                clip=true,
                grid=both,
                legend cell align=left,
                legend pos=south east,
                grid=major,
                xlabel={time (s)},
                ylabel={x displacement of tip}  
                ]                 
                \addplot[ 
                draw=black,
                smooth
                ]
                table[
                x index=0,
                y index=1
                ]
                {data/Codina.dat}; 
                \end{axis}
        \end{tikzpicture}
        \caption{x displacement of tip versus time for Bending Beam 1 and a density of 10.}
		\label{fig:Codina}
\end{figure}
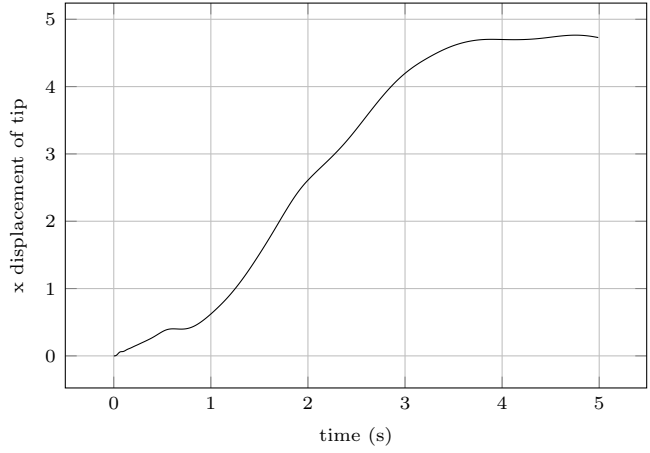

\begin{figure}[H]
        \centering
        \begin{tikzpicture}[trim axis left, trim axis right]
                \begin{axis}[
                scale=0.8,
                transform shape,
                width=.7\textwidth,
                height=.5\textwidth,
                label style={font=\scriptsize},
                tick label style={font=\scriptsize},
                legend style={font=\scriptsize},
                clip=true,
                grid=both,
                legend cell align=left,
                legend pos=south east,
                grid=major,
                xlabel={time (s)},
                ylabel={x displacement of tip}  
                ]                 
                \addplot[ 
                draw=black,
                smooth
                ]
                table[
                x index=0,
                y index=1
                ]
                {data/Codina2.dat}; 
                \end{axis}
        \end{tikzpicture}
        \caption{x displacement of tip versus time for Bending Beam 1 and a density of 1.}
		\label{fig:Codina2}
\end{figure}
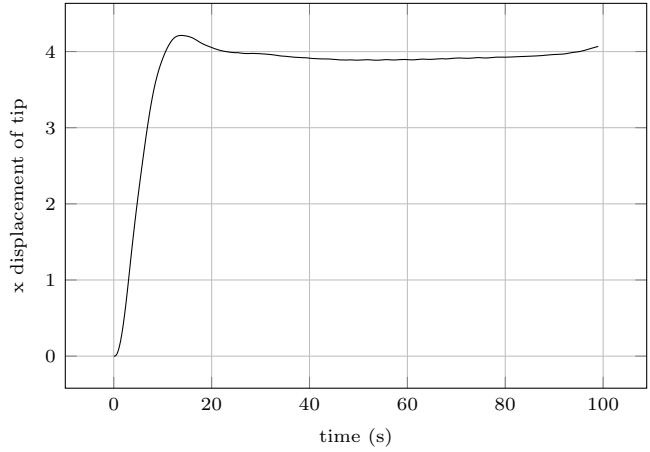

\subsection{Bending beam 2}

\begin{figure}[H]
\centering
   \includegraphics[width=1\linewidth]{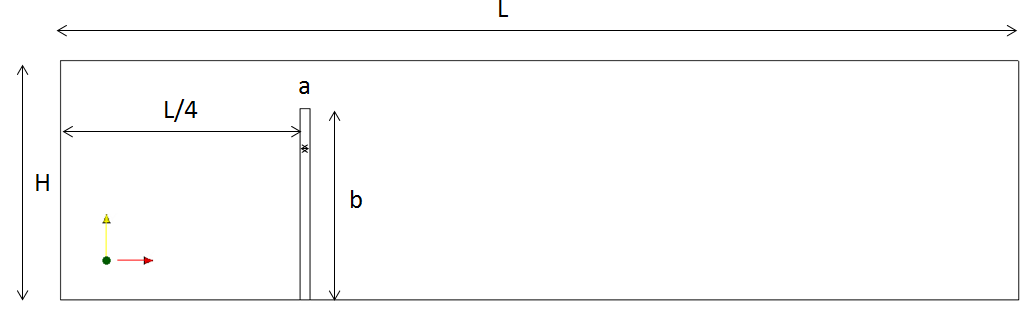}
\caption{Problem set up for bending beam 2}
\label{fig:BB2GEO}
\end{figure}

This section discusses another variation of the problem as described in \cite{HAN2020106179}\cite{article3}, where the aspect ratio of the beam is half that of the first iteration. 
The geometrical setup for this case is depicted in Figure \ref{fig:BB2GEO}. 
The fluid tunnel's height is $H=1$ cm, and its length is $L=4$ cm. The beam inside the tunnel has a thickness $a=0.04cm$ and a length of $b=0.8$ cm.
No-slip boundary conditions are imposed on the bottom part of the tunnel, while at the top, the flow is constrained to move only tangentially (symmetry). 
The inlet velocity is defined as $v_1(t)= 1.5 (-y^2+2y)$ cm/s, and $v_2(t)=0$. 
Zero gauge pressure is imposed at the outlet. 
The bottom part of the beam has zero Dirichlet boundary conditions. 
The fluid and solid properties are summarized in Table \ref{table:2} for two variations of the case.

\begin{table}[ht!]
\centering
\begin{tabular}{|p{2cm}|p{4cm}||p{2cm}|p{4cm}|}
 \hline
 \multicolumn{2}{|c|}{Fluid} & \multicolumn{2}{|c|}{Solid}\\
 \hline
$\rho_f$ & 1.0 g/$\text{cm}^3$ & $\rho_s$ & 7.8 g/$\text{cm}^3$\\
$\mu_f$ & 0.1 g/($\text{cm.s})$ & $\mu_s$ & 10$^5$ or 2x10$^{12}$ g/($\text{cm.s}^2)$ \\
&&$\nu_s$& 0.3\\
Model & Newtonian & Model &  St. Venant--Kirchhoff\\
\hline
\end{tabular} 
\caption{Fluid and solid properties for the bending beam 2 problem.}
\label{table:2}
\end{table}

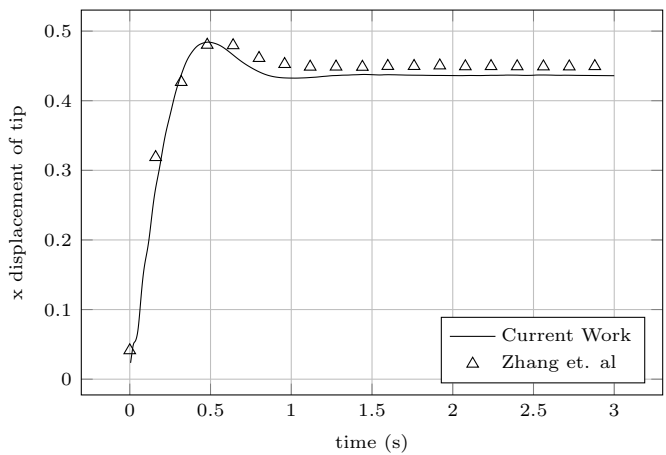
\begin{figure}[H]
        \centering
        \begin{tikzpicture}[trim axis left, trim axis right]
                \begin{axis}[
                scale=0.8,
                transform shape,
                width=.7\textwidth,
                height=.5\textwidth,
                label style={font=\scriptsize},
                tick label style={font=\scriptsize},
                legend style={font=\scriptsize},
                clip=true,
                grid=both,
                legend cell align=left,
                legend pos=south east,
                grid=major,
                xlabel={time (s)},
                ylabel={x displacement of tip}  
                ]                 
                
                \legend{Current Work, Zhang et. al}
                       
                \addplot[ 
                draw=black,
                smooth
                ]
                table[
                x index=0,
                y index=1
                ]
                {data/BB2.dat};

                \addplot[ 
                only marks,
    		    style={solid, fill=gray},
    		    mark=triangle,
   			    mark size=2.5pt
                ]
                 table[
                x index=0,
                y index=1
                ]
                {data/BB2B.dat};    
               
                \end{axis}
        \end{tikzpicture}
        \caption{Bending beam 2 x displacement of tip versus time and comparison with that from  Zhang et. al.}
		\label{fig:BB2D}
\end{figure}

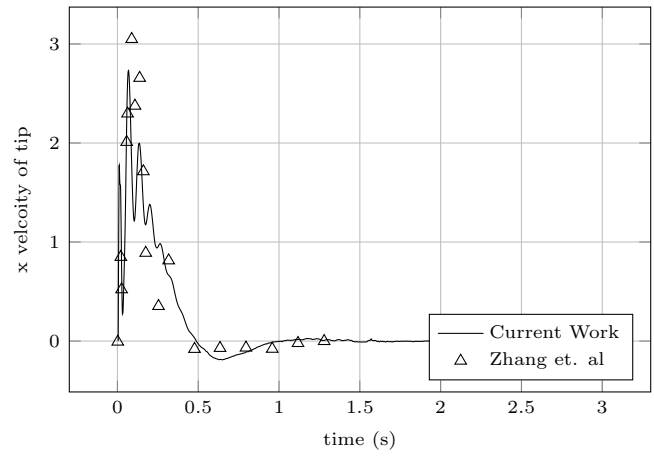
\begin{figure}[H]
        \centering
        \begin{tikzpicture}[trim axis left, trim axis right]
                \begin{axis}[
                scale=0.8,
                transform shape,
                width=.7\textwidth,
                height=.5\textwidth,
                label style={font=\scriptsize},
                tick label style={font=\scriptsize},
                legend style={font=\scriptsize},
                clip=true,
                grid=both,
                legend cell align=left,
                legend pos=south east,
                grid=major,
                xlabel={time (s)},
                ylabel={x velcoity of tip}  
                ]                 
                
                \legend{Current Work, Zhang et. al}
                       
                \addplot[ 
                draw=black,
                smooth
                ]
                table[
                x index=0,
                y index=1
                ]
                {data/BB2V.dat};

                \addplot[ 
                only marks,
    		    style={solid, fill=gray},
    		    mark=triangle,
   			    mark size=2.5pt
                ]
                 table[
                x index=0,
                y index=1
                ]
                {data/BB2BV.dat};    
               
                \end{axis}
        \end{tikzpicture}
        \caption{Bending beam 2 x velocity of tip versus time and comparison with that from  Zhang et. al.}
		\label{fig:BB2V}
\end{figure}

\begin{figure}[H]
\centering
\begin{subfigure}{0.7\textwidth}
   \includegraphics[width=1\linewidth]{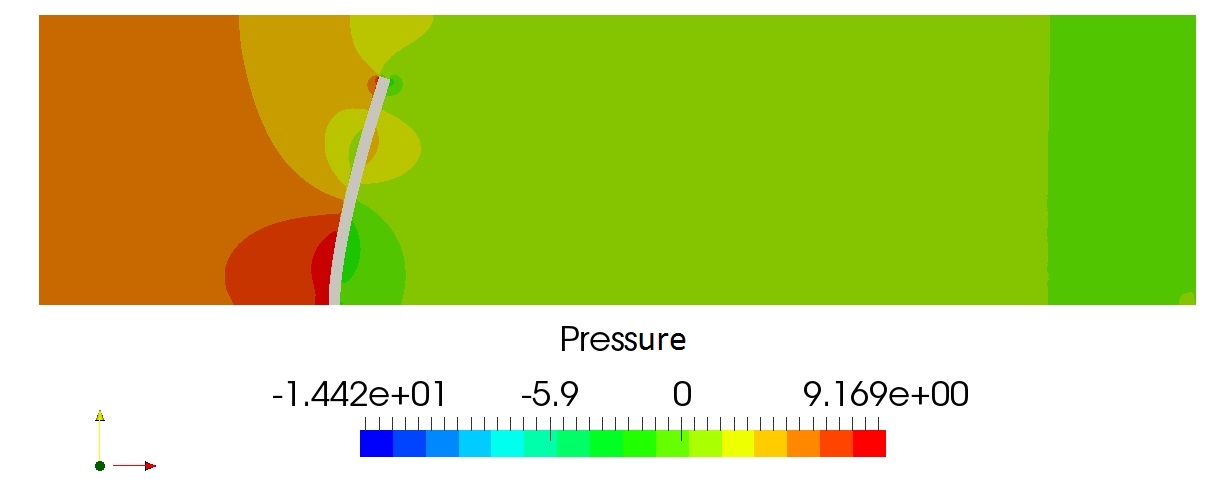}
  \end{subfigure}%
  \hspace*{\fill}   
  \begin{subfigure}{0.3\textwidth}
   \includegraphics[width=1\linewidth]{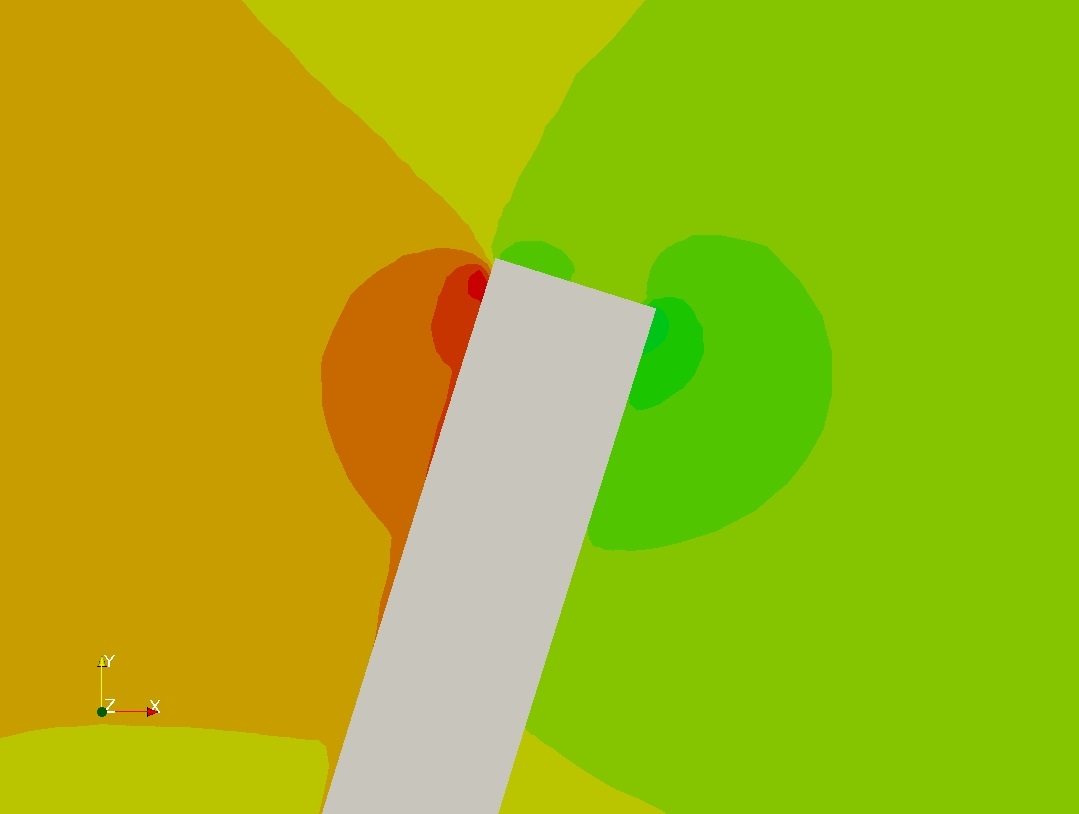}
  \end{subfigure}%
  \\
\centering
\begin{subfigure}{0.7\textwidth}
   \includegraphics[width=1\linewidth]{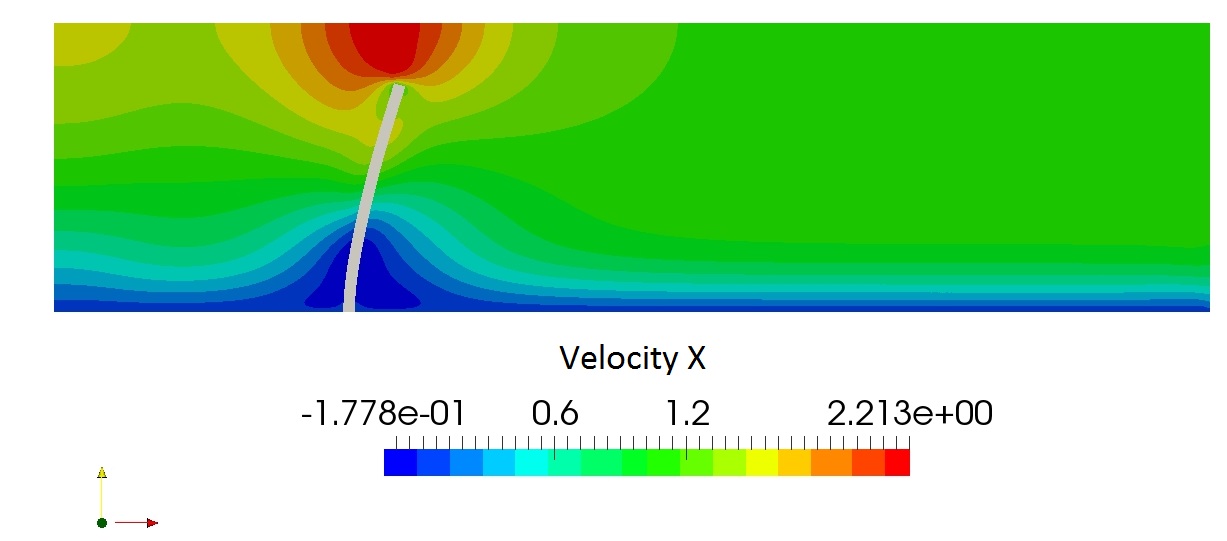}
  \end{subfigure}%
  \hspace*{\fill}   
  \begin{subfigure}{0.3\textwidth}
   \includegraphics[width=1\linewidth]{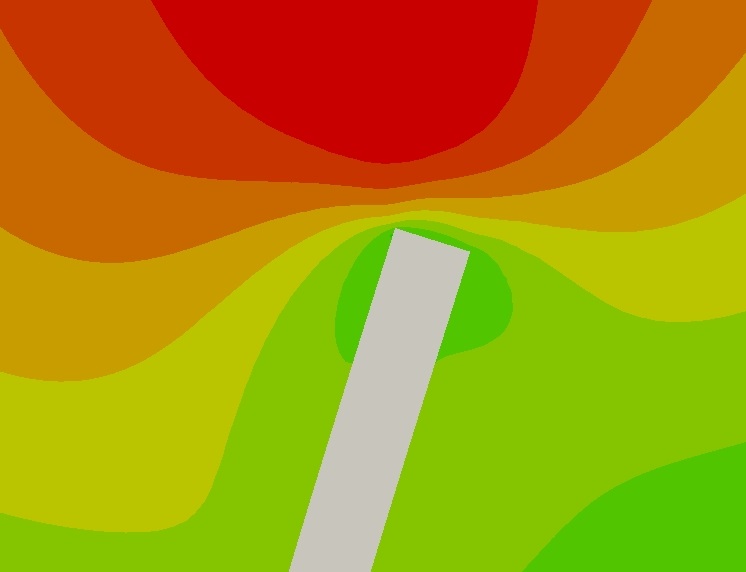}
  \end{subfigure}%
\caption{Velocity and pressure contours at t=0.1 s for bending beam 2}
\label{fig:t=0.1}
\end{figure}

\begin{figure}[H]
\centering
\begin{subfigure}{0.7\textwidth}
   \includegraphics[width=1\linewidth]{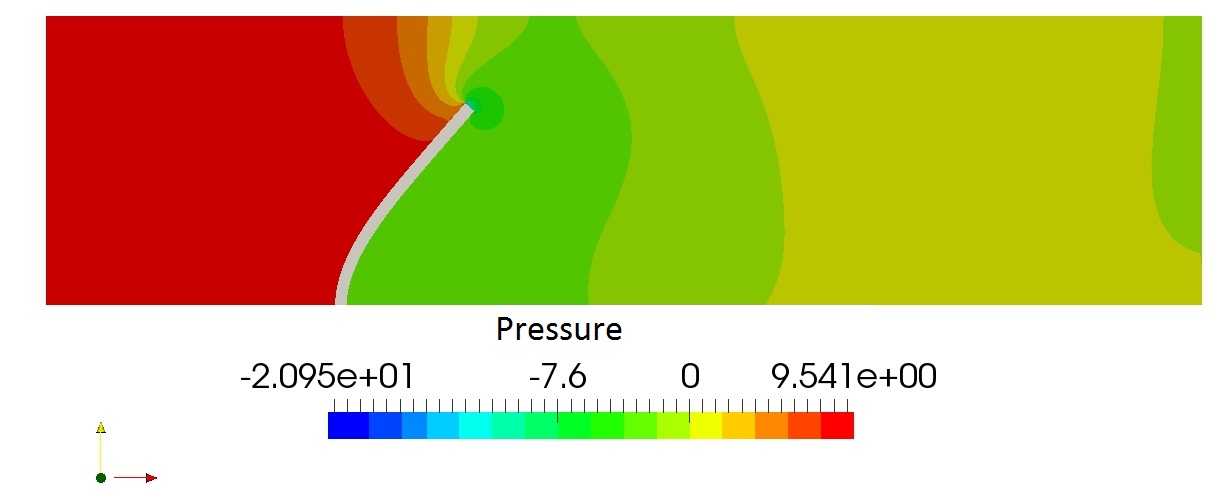}
  \end{subfigure}%
  \hspace*{\fill}   
  \begin{subfigure}{0.3\textwidth}
   \includegraphics[width=1\linewidth]{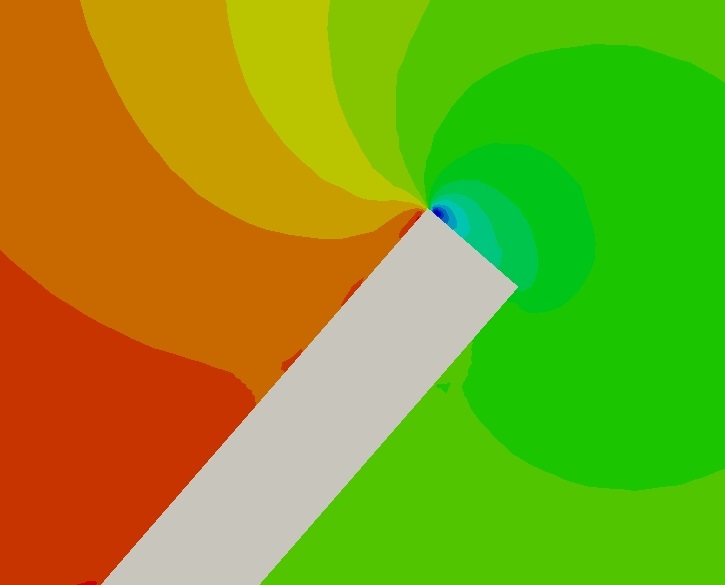}
  \end{subfigure}%
  \\
\centering
\begin{subfigure}{0.7\textwidth}
   \includegraphics[width=1\linewidth]{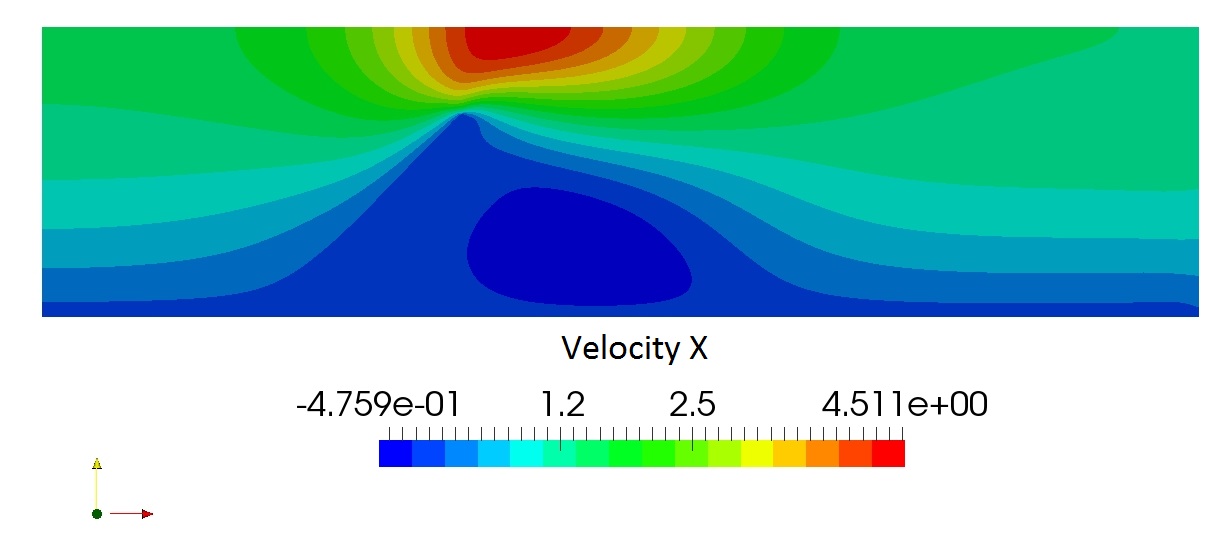}
  \end{subfigure}%
  \hspace*{\fill}   
  \begin{subfigure}{0.3\textwidth}
   \includegraphics[width=1\linewidth]{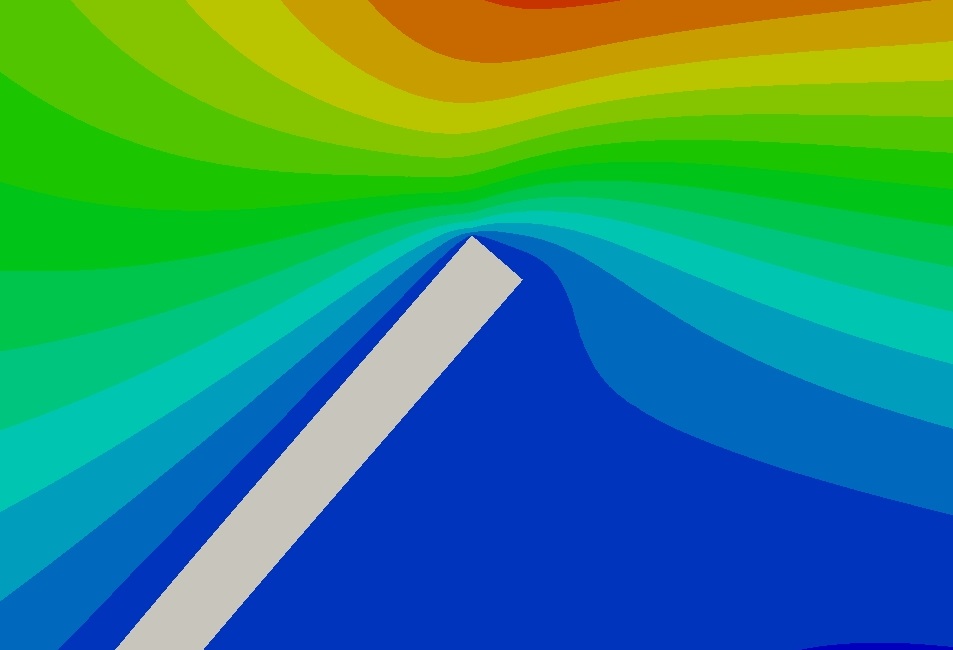}
  \end{subfigure}%
\caption{Velocity and pressure contours at t=0.8 s for bending beam 2}
\label{fig:t=0.8}
\end{figure}

\begin{figure}[H]
\centering
\begin{subfigure}{0.7\textwidth}
   \includegraphics[width=1\linewidth]{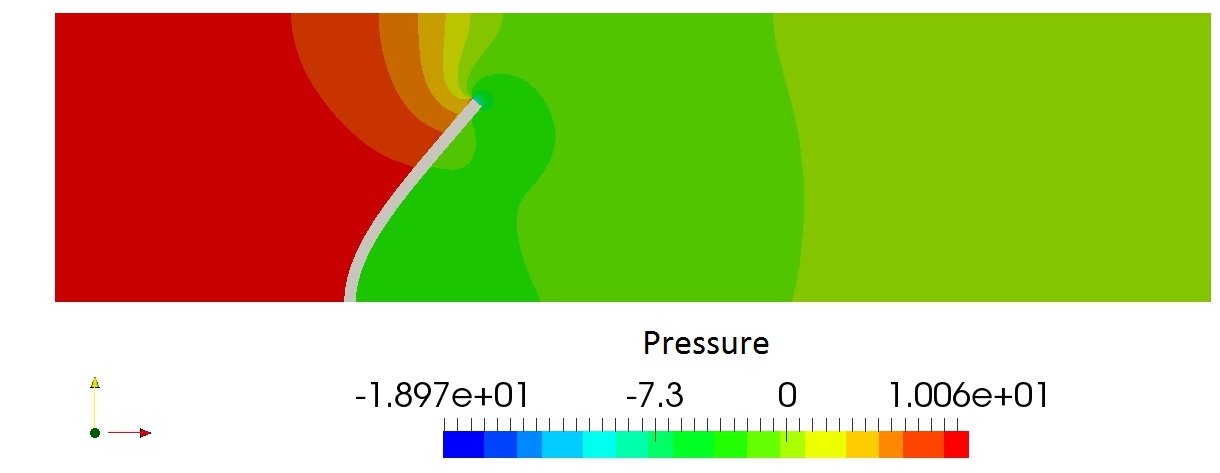}
  \end{subfigure}%
  \hspace*{\fill}   
  \begin{subfigure}{0.3\textwidth}
   \includegraphics[width=1\linewidth]{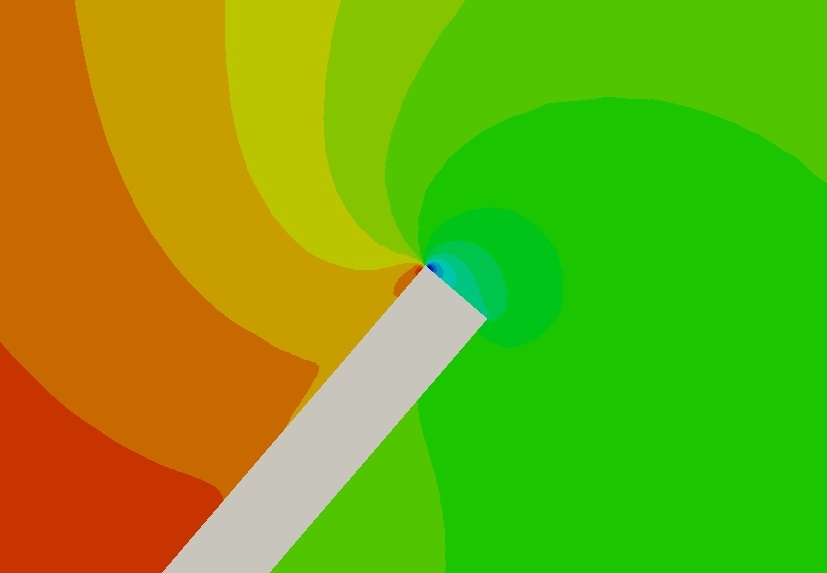}
  \end{subfigure}%
  \\
\centering
\begin{subfigure}{0.7\textwidth}
   \includegraphics[width=1\linewidth]{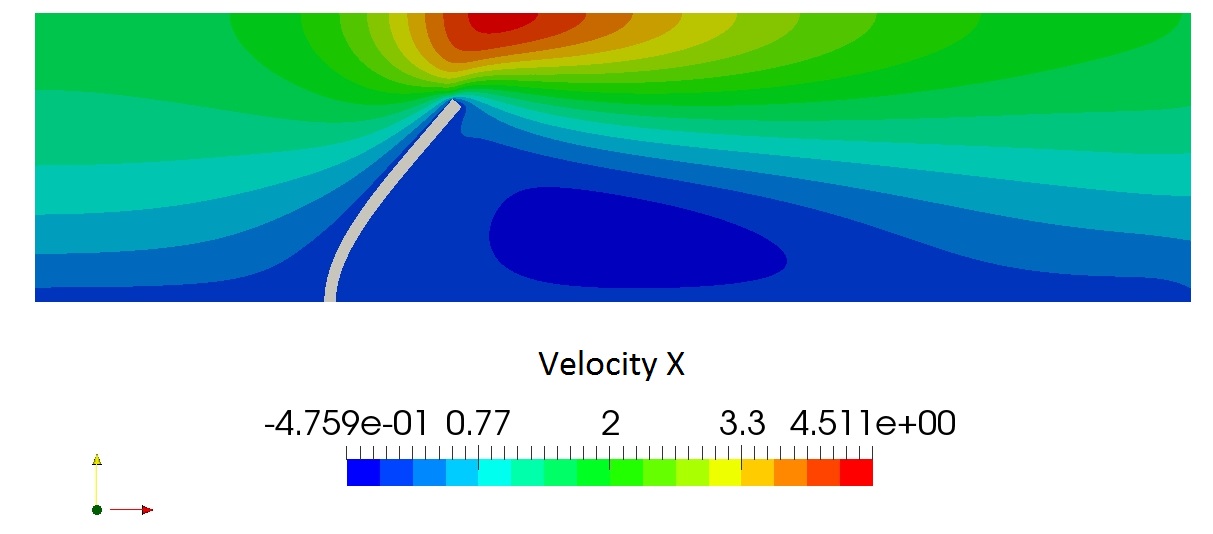}
  \end{subfigure}%
  \hspace*{\fill}   
  \begin{subfigure}{0.3\textwidth}
   \includegraphics[width=1\linewidth]{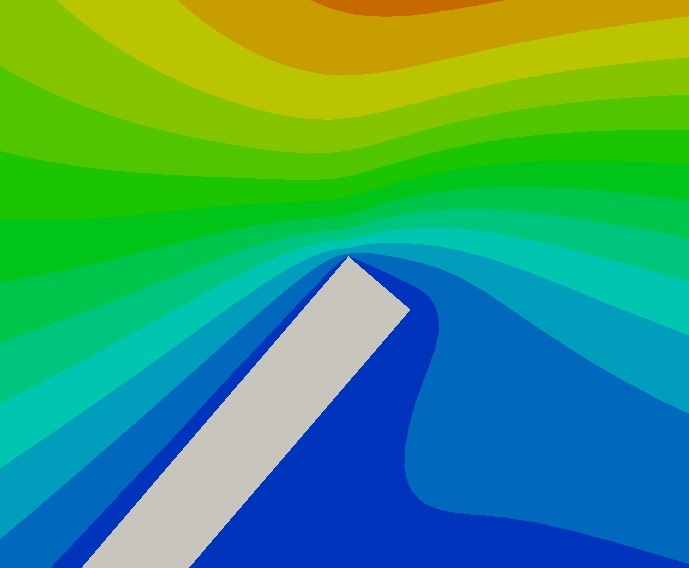}
  \end{subfigure}%
\caption{Velocity and pressure contours at t=3 s for bending beam 2}
\label{fig:t=3}
\end{figure}

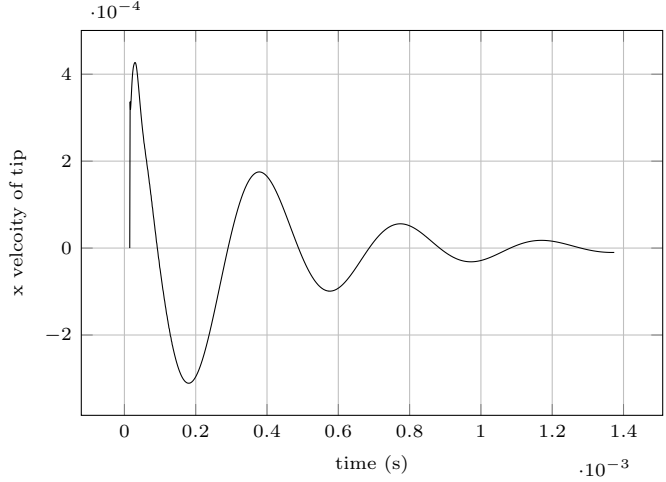
\begin{figure}[H]
        \centering
        \begin{tikzpicture}[trim axis left, trim axis right]
                \begin{axis}[
                scale=0.8,
                transform shape,
                width=.7\textwidth,
                height=.5\textwidth,
                label style={font=\scriptsize},
                tick label style={font=\scriptsize},
                legend style={font=\scriptsize},
                clip=true,
                grid=both,
                legend cell align=left,
                legend pos=south east,
                grid=major,
                xlabel={time (s)},
                ylabel={x velcoity of tip}  
                ]                 
                \addplot[ 
                draw=black,
                smooth
                ]
                table[
                x index=0,
                y index=1
                ]
                {data/BB22V.dat}; 
                \end{axis}
        \end{tikzpicture}
        \caption{Bending beam 2 x velocity of tip versus time for the second variation.}
		\label{fig:BB22V}
\end{figure}

Initially, both the fluid and the solid are at rest. 
In both variations of the case, the system will reach a steady-state solution where the solid will no longer oscillate.
In the first case, where the Young's modulus is relatively low, significant deformations are observed in the solid structure due to its low stiffness. 
The displacement and velocity versus time are shown in Figures \ref{fig:BB2D} and \ref{fig:BB2V}, respectively.
In the second case, a material with a much higher Young's modulus is considered. 
As a result, a much lower time step is required, approximately 1000 times lower. 
This is because the high stiffness of the beam leads to very small deformations. 
The velocity of the tip of the beam in the x-direction is displayed in Figure \ref{fig:BB22V} for a time step $\Delta t= 0.000001 s$.
The vibration of the beam is induced by the fluid traction forces. These vibrations quickly dissipate due to the viscous effect of the fluid. 
Despite no damping being considered in the solid, the velocity decays until it reaches 0.
Pressure oscillations occur in the vicinity of the beam, primarily because of the high gradient in the velocity caused by sharp corners. 
However, these pressure oscillations are dissipated further from the beam. 
Notably, no oscillations are observed in the velocity field itself.
The figures in this section include a comparison with the work by Zhang et. al, offering valuable insights into the accuracy and agreement of the present study with earlier research.
Velocity and pressure contours of the fluid on the fluid-solid mesh are depicted in Figures \ref{fig:t=0.1}, \ref{fig:t=0.8}, and \ref{fig:t=3}.
Overall, the results presented here illustrate the complex interplay between fluid and solid dynamics in the context of Bending Beam 2. 
The choice of material properties, in particular the Young's modulus, significantly impacts the structural behavior and the time-stepping requirements, which is a key aspect of this problem.

\subsection{2D flow induced vibration of an elastic plate}

\begin{figure}[!h]
\centering
   \includegraphics[width=1\linewidth]{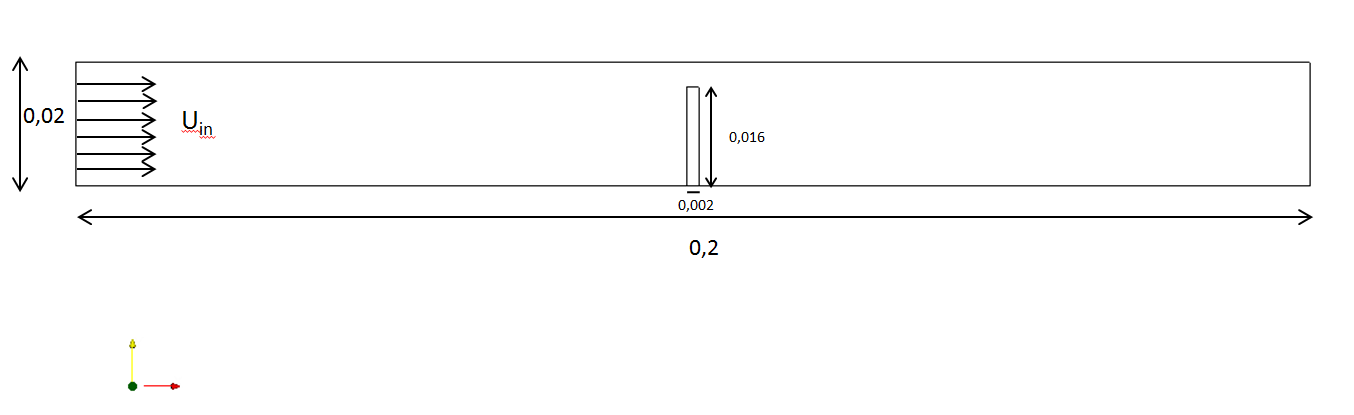}
\caption{Problem set up for the 2D flow-induced vibration of an elastic plate.}
\label{fig:Capture}
\end{figure}

In this case, which was first conducted by \cite{xia2008unstructured}, a plate is placed at the center of the channel as shown in figure \ref{fig:Capture}.
The length and height of the channel are $0.2 m$, and $0.02 m$ respectively.
The thickness and height of the elastic plate are $0.002 m$, and $0.016 m$ respectively.
The fluid and solid properties are tabulated in Table \ref{table:0}.
No-slip Dirichlet boundary conditions are applied at the top and bottom of the 2D channel.
A zero-gauge pressure outlet is imposed.
At the inlet, we impose a sinusoidal velocity function given by $U_{in}=0.015 sin(2\pi t)$, with an equivalent period of $1 s$.
The peak velocity magnitude is equal to $0.015 m/s$, which is equivalent to a Reynolds number of $300$.

\begin{table}[!h]
\centering
\begin{tabular}{|p{3cm}|p{3cm}||p{3cm}|p{3cm}|}
 \hline
 \multicolumn{2}{|c|}{Fluid properties} & \multicolumn{2}{|c|}{Solid properties}\\
 \hline
$\rho_f$ & 1000 Kg/m$^3$ & $\rho_s$ & 1000 Kg/m$^3$\\
$\mu_f$ & 0.001 Kg/ms & $\mu_s$ & 1677.85 Pa\\
&&$E$& 5000 Pa\\
&&$\nu$& 0.49 \\
Model & Newtonian & Model &  Neo-Hookean\\
\hline
\end{tabular} 
\caption{Fluid and solid properties for the 2D flow induced vibration of an elastic plate.}
\label{table:0}
\end{table}

\begin{figure}[!h]
\centering
    \begin{subfigure}{0.9\linewidth}
        \includegraphics[width=\linewidth]{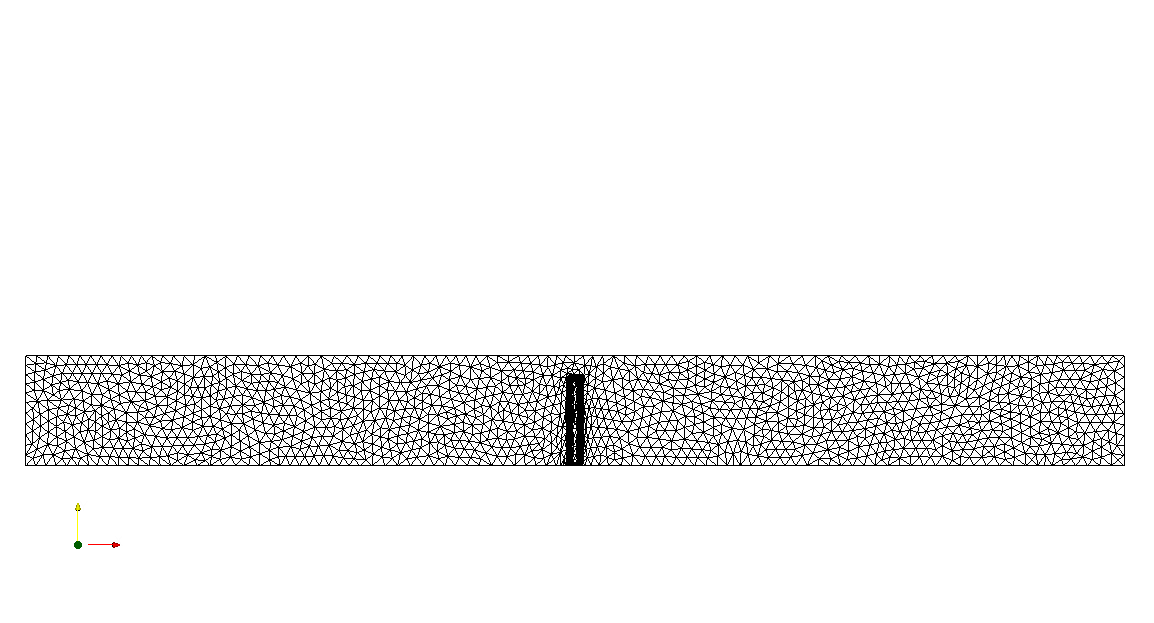}
        \caption{Fluid-Solid Mesh}
        \label{fig:50a}
    \end{subfigure}
    \\
    \begin{subfigure}{0.4\linewidth}
        \includegraphics[width=\linewidth]{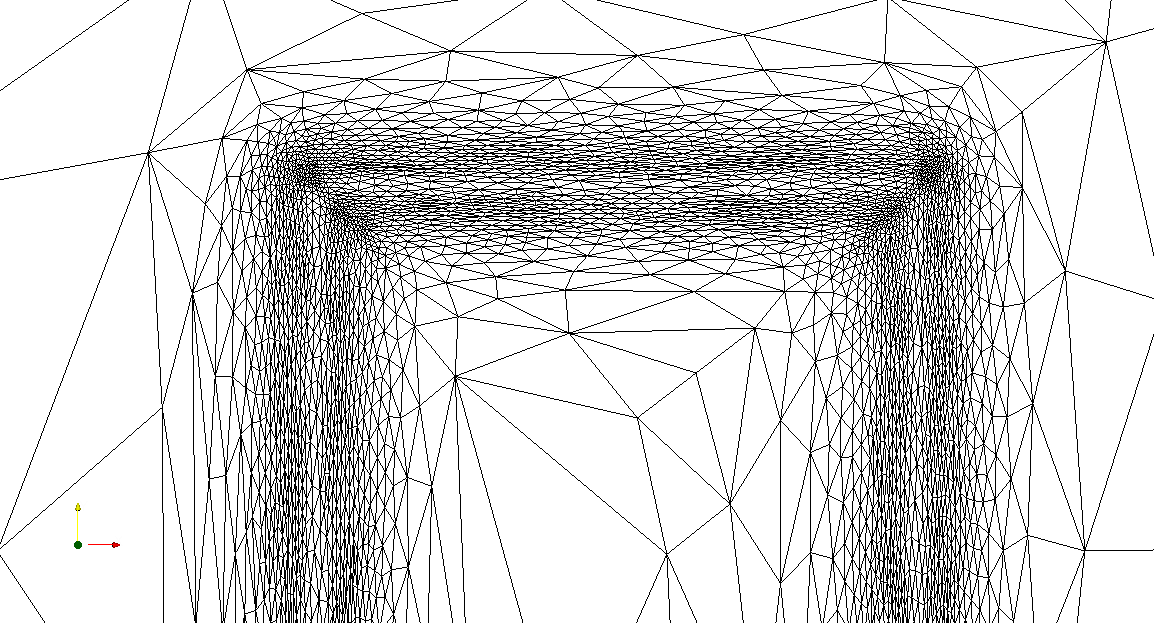}
        \caption{Anisotropic mesh at the interface}
        \label{fig:50b}
    \end{subfigure}
    \begin{subfigure}{0.2\linewidth}
        \includegraphics[width=\linewidth]{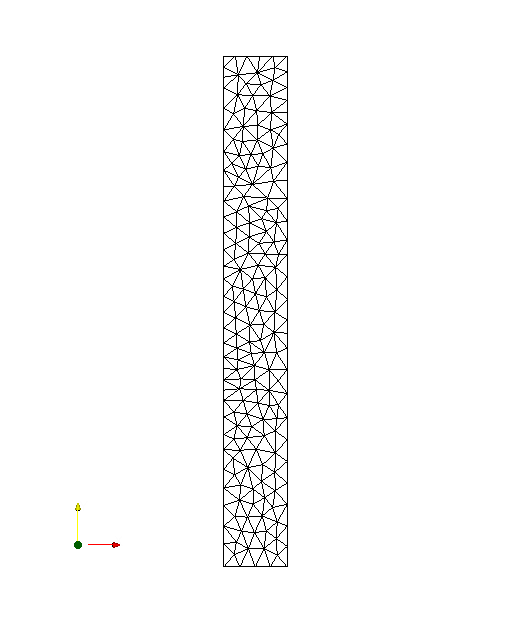}
        \caption{Solid Mesh}
        \label{fig:50c}
    \end{subfigure} 
\caption{Fluid-solid and solid meshes for the 2D flow induced vibration of an elastic plate. }
\label{fig:0Mesh}
\end{figure}

As stated earlier, we use an anisotropic mesh adaptation for the fluid-solid mesh.
It can be seen from Figure \ref{fig:50a}, that the elements are localized at the interface.
A magnified picture at the interface is shown in Figure \ref{fig:50b} that shows the anisotropic properties of the elements.
The number of elements is capped at $30 000$.
The solid mesh is shown in Figure \ref{fig:50c}, where the number of elements is equivalent to 324 elements.
The time step for the simulation is set at $0.005 s$.

\begin{figure}[p]
\centering
    \begin{subfigure}{0.8\linewidth}
        \includegraphics[width=\linewidth]{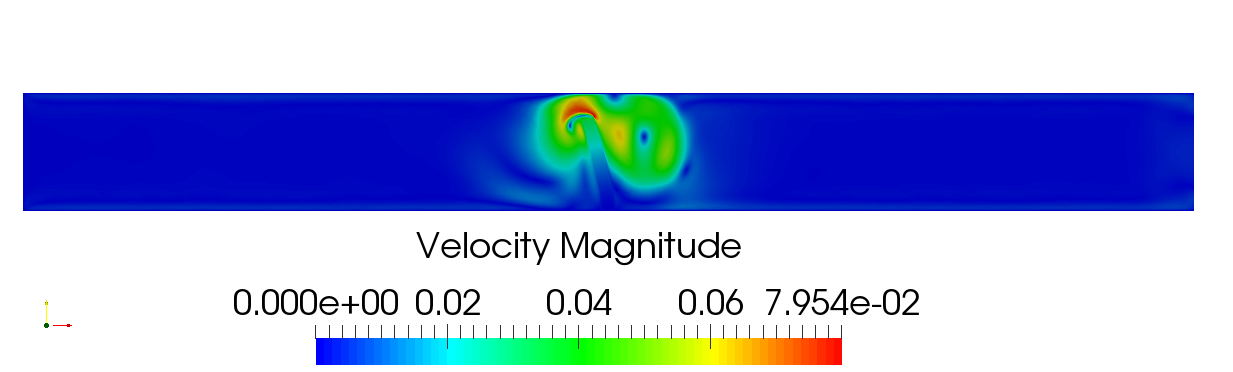}
        \caption{t=10}
        \label{fig:00a}
    \end{subfigure} 
\\
    \begin{subfigure}{0.8\linewidth}
        \includegraphics[width=\linewidth]{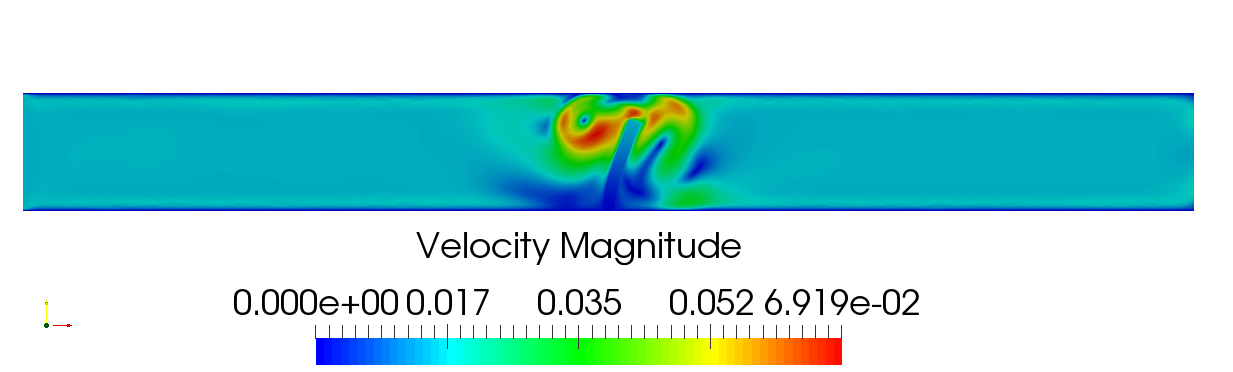}
        \caption{t=10.25}
        \label{fig:00b}
    \end{subfigure} 
\\
    \begin{subfigure}{0.8\linewidth}
        \includegraphics[width=\linewidth]{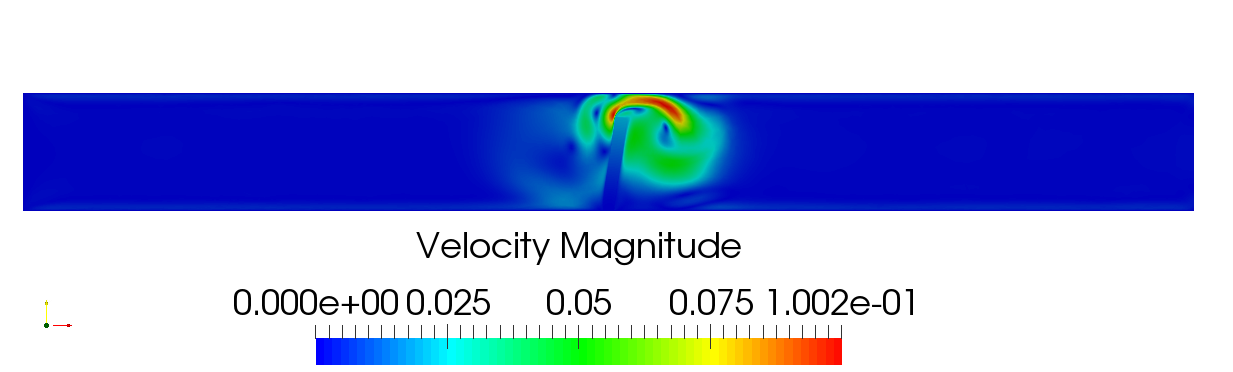}
        \caption{t=10.5}
        \label{fig:00c}
    \end{subfigure} 
\\
    \begin{subfigure}{0.8\linewidth}
        \includegraphics[width=\linewidth]{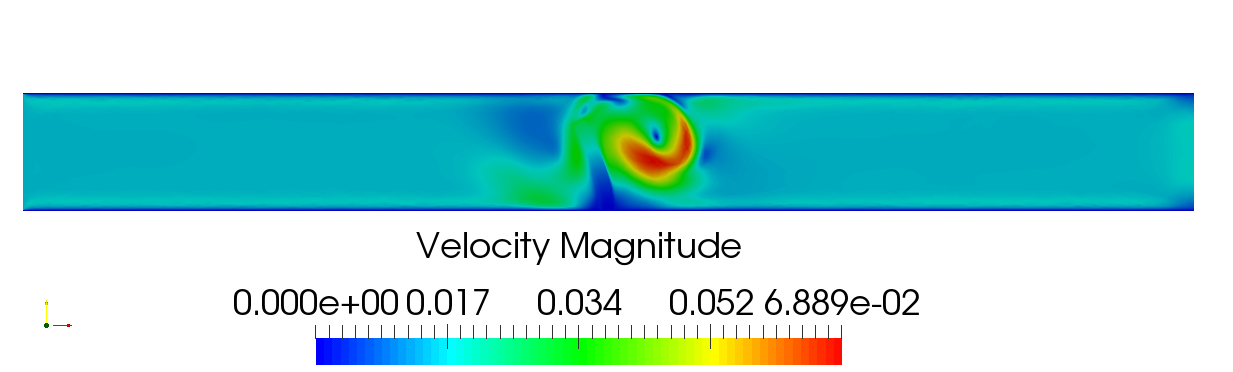}
        \caption{t=10.75}
        \label{fig:00d}
    \end{subfigure} 
\\
    \begin{subfigure}{0.8\linewidth}
        \includegraphics[width=\linewidth]{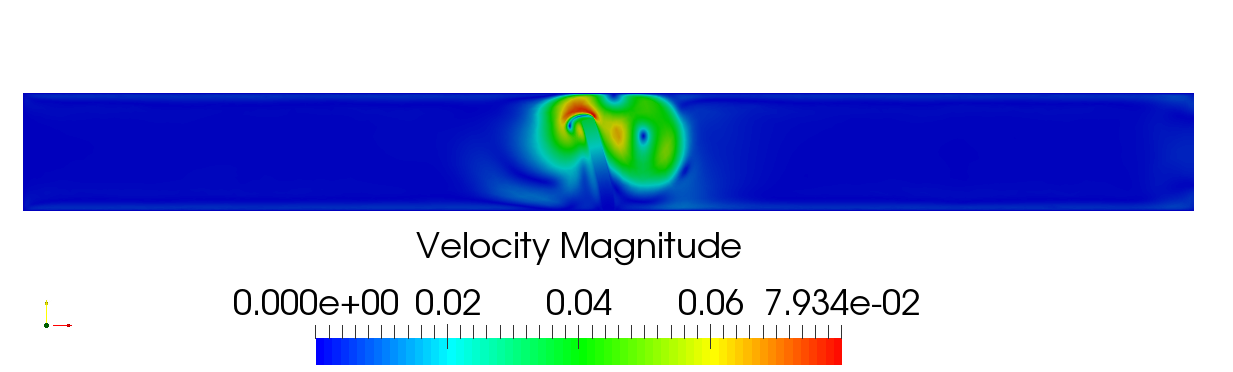}
        \caption{t=11}
        \label{fig:00e}
    \end{subfigure} 
\caption{Velocity magnitude field of the fluid and solid on the fluid-solid mesh at different positions in time for 2D flow induced vibration of an elastic plate.}
\label{fig:2Dt}
\end{figure}

\begin{figure}[!h]
\centering
    \begin{subfigure}{0.6\linewidth}
        \includegraphics[width=\linewidth]{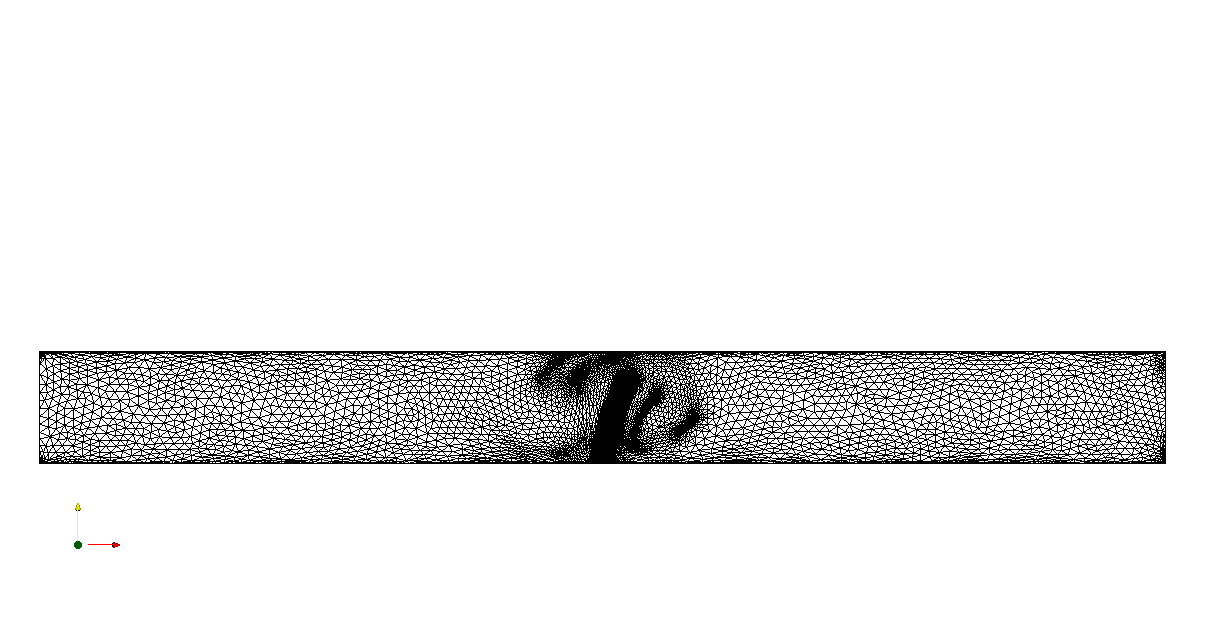}
        \caption{}
        \label{fig:000a}
    \end{subfigure} 
    \begin{subfigure}{0.4\linewidth}
        \includegraphics[width=\linewidth]{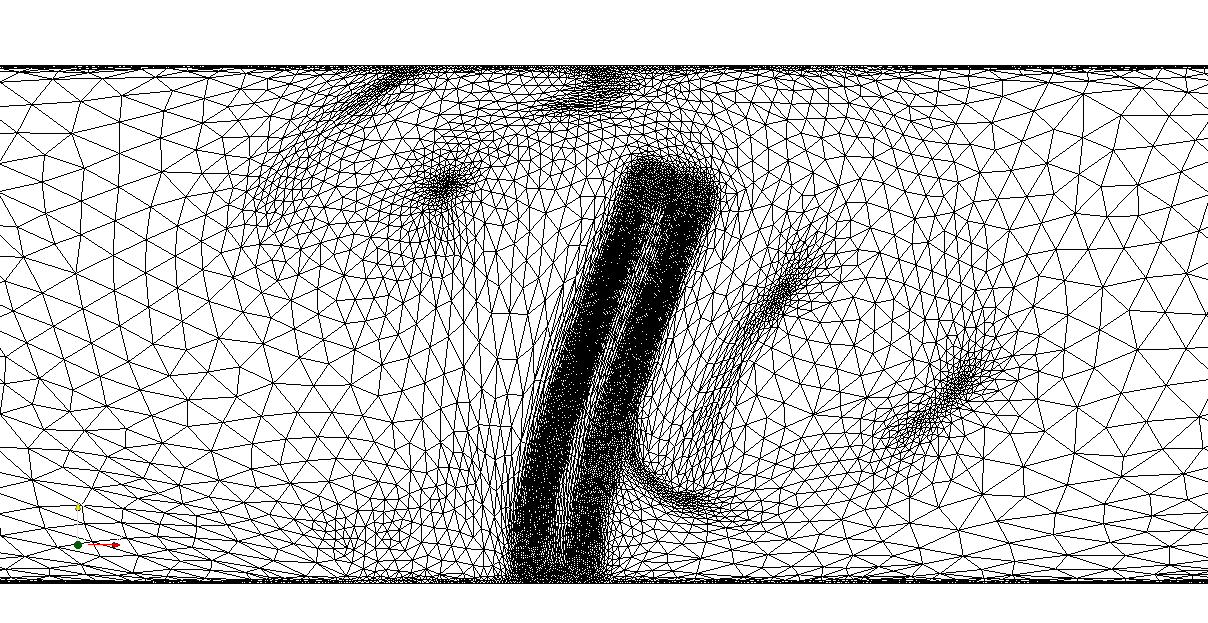}
        \caption{}
        \label{fig:000b}
    \end{subfigure} 
\caption{Full \ref{fig:000a} and magnified \ref{fig:000b} Fluid-Solid Mesh with mesh adaptation on different criteria at time  t.}
\label{fig:00Mesh}
\end{figure}

\begin{figure}[!h]
\centering
    \begin{subfigure}{0.3\linewidth}
        \includegraphics[width=\linewidth]{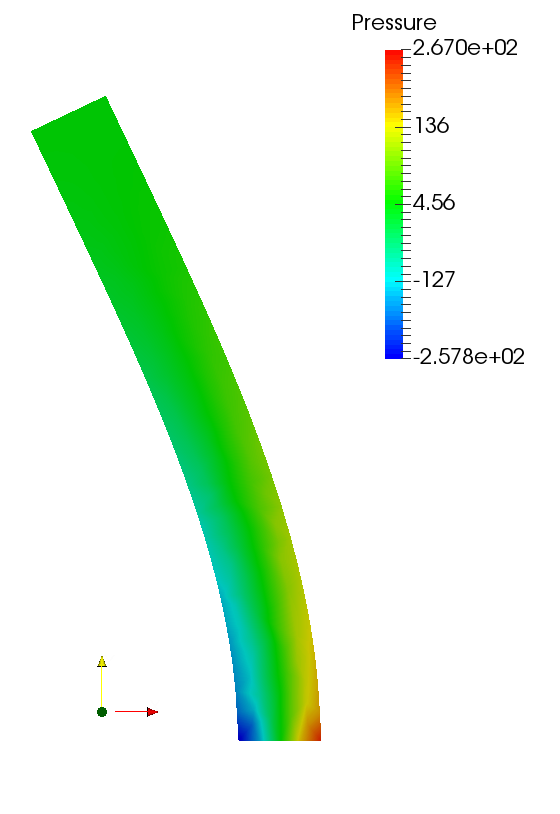}
        \caption{t=10.75}
        \label{fig:0000a}
    \end{subfigure}
    \begin{subfigure}{0.3\linewidth}
        \includegraphics[width=\linewidth]{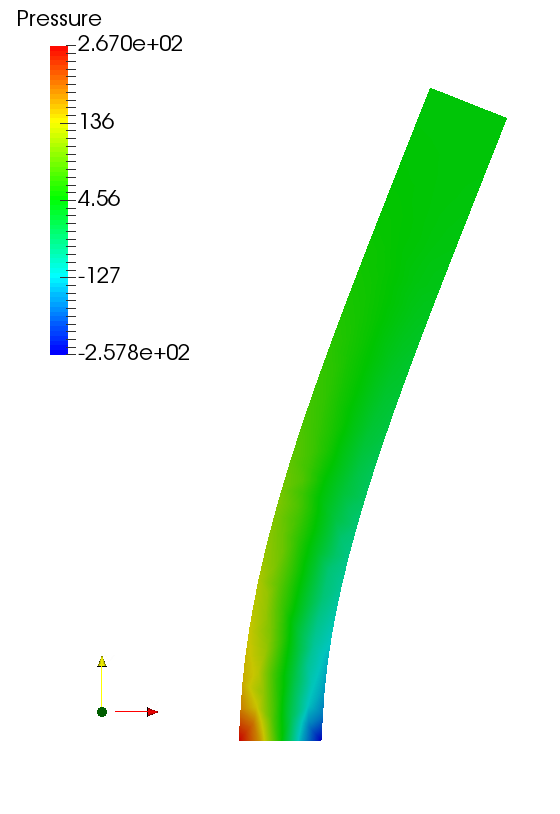}
        \caption{t=11.25}
        \label{fig:0000b}
    \end{subfigure}
\caption{Solid pressure contours at different times for 2D flow induced vibration of an elastic plate.}
\label{fig:Pressure}
\end{figure}

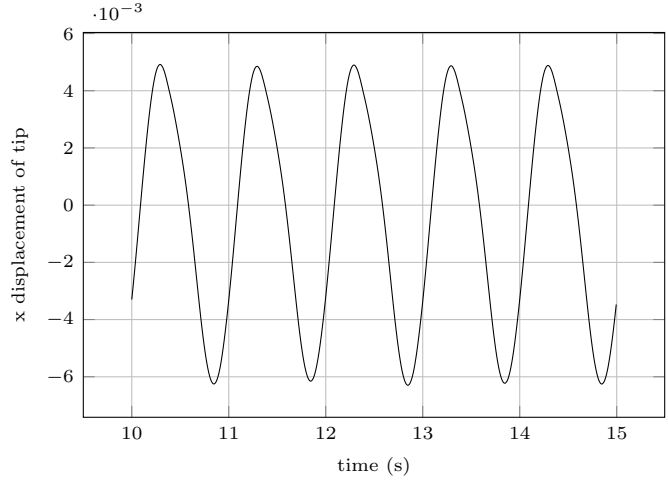
\begin{figure}[!h]
        \centering
        \begin{tikzpicture}[trim axis left, trim axis right]
                \begin{axis}[
                scale=0.8,
                transform shape,
                width=.7\textwidth,
                height=.5\textwidth,
                label style={font=\scriptsize},
                tick label style={font=\scriptsize},
                legend style={font=\scriptsize},
                clip=true,
                grid=both,
                legend cell align=left,
                legend pos=south east,
                grid=major,
                xlabel={time (s)},
                ylabel={x displacement of tip}  
                ]                 
                \addplot[ 
                draw=black,
                smooth
                ]
                table[
                x index=0,
                y index=1
                ]
                {data/2DUx.dat}; 
                \end{axis}
        \end{tikzpicture}
        \caption{x displacement of top right node versus time for 2D flow induced vibration of an elastic plate.}
		\label{fig:2DUx}
\end{figure}

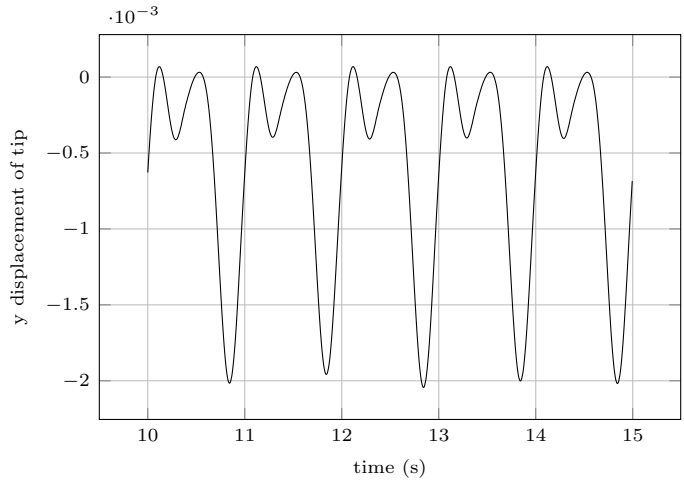
\begin{figure}[!h]
        \centering
        \begin{tikzpicture}[trim axis left, trim axis right]
                \begin{axis}[
                scale=0.8,
                transform shape,
                width=.7\textwidth,
                height=.5\textwidth,
                label style={font=\scriptsize},
                tick label style={font=\scriptsize},
                legend style={font=\scriptsize},
                clip=true,
                grid=both,
                legend cell align=left,
                legend pos=south east,
                grid=major,
                xlabel={time (s)},
                ylabel={y displacement of tip}  
                ]                 
                \addplot[ 
                draw=black,
                smooth
                ]
                table[
                x index=0,
                y index=1
                ]
                {data/2DUy.dat}; 
                \end{axis}
        \end{tikzpicture}
        \caption{y displacement of top right node versus time for 2D flow induced vibration of an elastic plate.}
		\label{fig:2DUy}
\end{figure}

Given the sinusoidal nature of the flow, the elastic plate will swing back and forth from its original position.
We ran the simulation for $15T$.
We are interested in the deflection of the plate at $T/4$ of the beginning of the period where the inlet velocity is at its maximum.
The velocity magnitude field for both the fluid and the solid on the fluid-solid mesh for different positions in time are shown in Figure \ref{fig:2Dt}.
The mesh adaptation highlighting the interface, and the vortices of the flow is shown in Figure \ref{fig:00Mesh}.
Pressure contours inside the solid are shown in Figure \ref{fig:0000b}, for two different positions in time, that showcases the tension and compression that the plate exhibits when undergoing significant bending.
The x and y displacements of the top right node are plotted versus time in Figures \ref{fig:2DUx} and \ref{fig:2DUy} respectively.

\subsection{Turek's FSI benchmark}

In this section, we explore the FSI2 and FSI3 variations of Turek's Fluid-Structure Interaction (FSI) benchmark \cite{schafer1996benchmark}\cite{turek2006proposal}. 
This benchmark has been widely adopted in computational fluid dynamics and serves as an essential validation tool for FSI simulations. 
We will look into the benchmark setup, boundary conditions, parameters, and the results of FSI2 and FSI3 cases.

\begin{figure}[!h]
\centering
\begingroup  
\setlength{\tabcolsep}{0pt}
\renewcommand{\arraystretch}{0}
\begin{tabular}{c c}
    \begin{subfigure}{0.7\linewidth}
        \includegraphics[width=\linewidth]{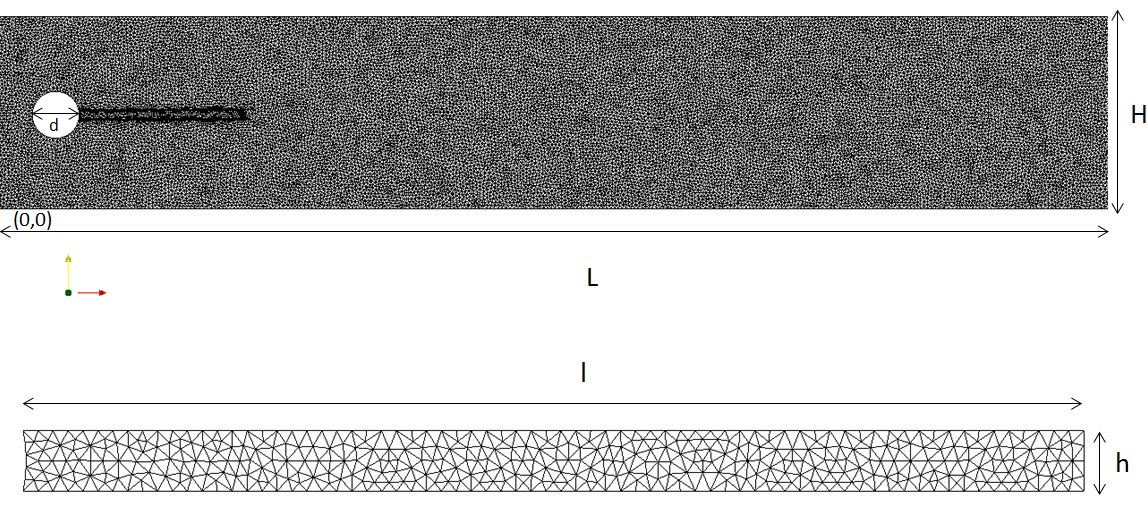}
        \caption{}
        \label{fig:42a}
    \end{subfigure}
    \begin{subfigure}{0.3\linewidth}
        \includegraphics[width=\linewidth]{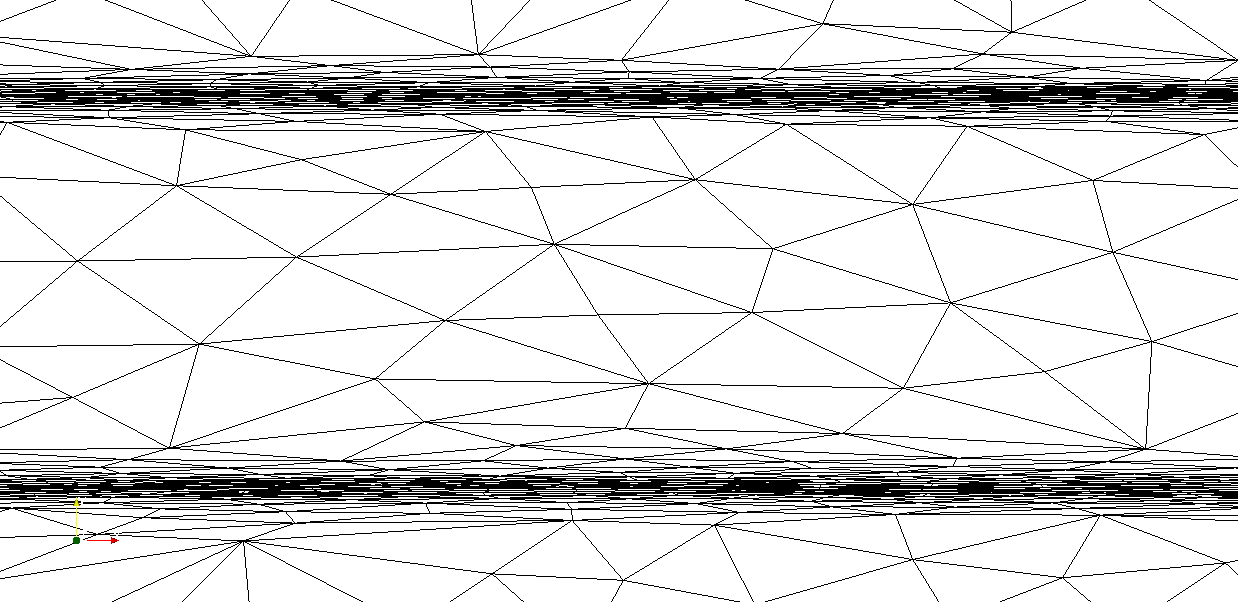}
        \caption{}
        \label{fig:42b}
    \end{subfigure}
\end{tabular}
\endgroup
\caption{Problem set up and mesh}
\label{fig:FSI2ProblemDescription}
\end{figure}

The geometrical setup of FSI2 and FSI3 is consistent with the well-established benchmarks in CFD. 
As shown in Figure \ref{fig:FSI2ProblemDescription}, a 2D channel is considered with a cylindrical solid placed inside. 
Key parameters include the dimensions of the channel, cylinder, solid, and the control point where data is collected.
The bottom left of the domain is at $(0,0)$. The height is $H=0.41$, and the length is $L=2.5$.
The center point of the cylinder is positioned at $(0.2,0.2)$, with a diameter $d=0.1$.
The bottom right of the solid is positioned at $(0.6,0.19)$, and the left part is fully fixated on the cylinder.
The control point for which the values are collected is positioned at $(0.6,0.2)$.
The solid has a length of $l=0.35$, and a height of $h=0.02$.
A non-symmetry is intentionally prevalent in the y-direction, to avoid the dependency of the beginning of the oscillations on the calculations.
Figure \ref{fig:42a} shows the fluid-solid, and solid mesh, and \ref{fig:42b} shows a close-up of the anisotropic stretched elements at the interface.

The benchmark involves setting up essential boundary conditions to replicate real-world flow behavior. 
These include no-slip conditions on the channel's top and bottom walls, a parabolic inlet velocity profile, and a smooth increase in velocity at the inlet:

\begin{equation}
    \vv_f=1.5 \overline{U} \frac{y(H-y)}{(\frac{H}{2})^2}=1.5 \overline{U} \frac{4}{0.1681}y(0.41-y).
\end{equation}

\begin{equation}
    \vv_f(t)=\vv_f(1-e^{-(2.5(t))^2})
\end{equation}

This will ensure a mean velocity of $\overline{U}$ at the inlet, and a maximum velocity of $1.5\overline{U}$.
The outlet features a zero-gauge pressure condition.
Both FSI2 and FSI3 cases involve varying fluid and solid parameters. 
Table \ref{table:3} provides a summary of the parameters.

\begin{table}[!h]
\centering
\begin{tabular}{|p{6cm}|p{4cm}||p{4cm}|}
 \hline
parameter & FSI2 & FSI3 \\
$\rho_s$  [Kg/$\text{m}^3$] & 10000 & 1000\\
$\nu_s$ & 0.4 & 0.4\\
$\mu_s$ [Kg/($\text{m.s}^2)$] & 500000  & 2000000\\
\hline
$\rho_f$  [Kg/$\text{m}^3$] & 1000 & 1000\\
$\mu_f$  [Kg/($\text{m.s})$] & 1 & 1\\
\hline
$\overline{U}$ [m/s] & 1 & 2\\
\hline
$Re=\frac{\rho_f \overline{U} d}{\mu_f}$  & 100 & 200\\
$Ae=\frac{E}{\rho_f \overline{U}^2}$  & 1.4 x 10$^3$ & 1.4 x 10$^3$\\

\hline
\end{tabular} 
\caption{Fluid and solid parameters for the Turek's FSI benchmarks variations.}
\label{table:3}
\end{table}

\begin{figure}[!h]
\centering
\includegraphics[width=.7\linewidth]{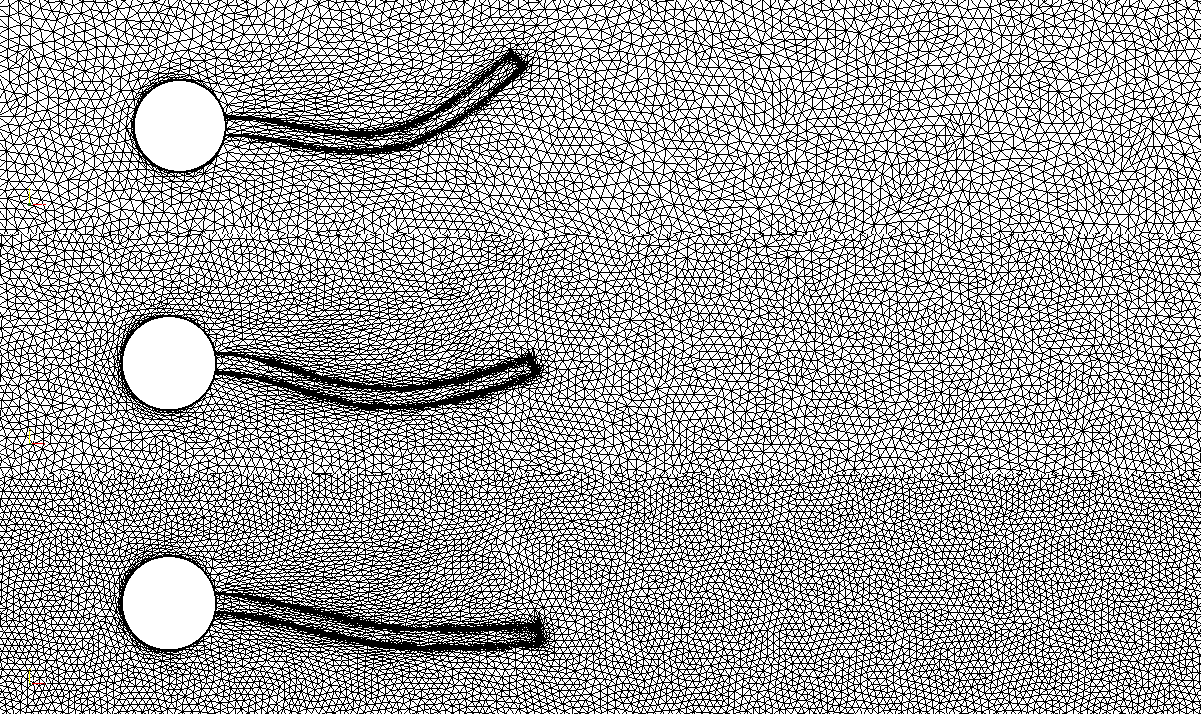}
\caption{Fluid mesh convergence study consisting of three different meshes of  20 000, 30 000, and 40 000 elements respectively.}
\label{fig:FSI2FluidMeshConvergence}
\end{figure}

\begin{figure}[!h]
\centering
\includegraphics[width=.7\linewidth]{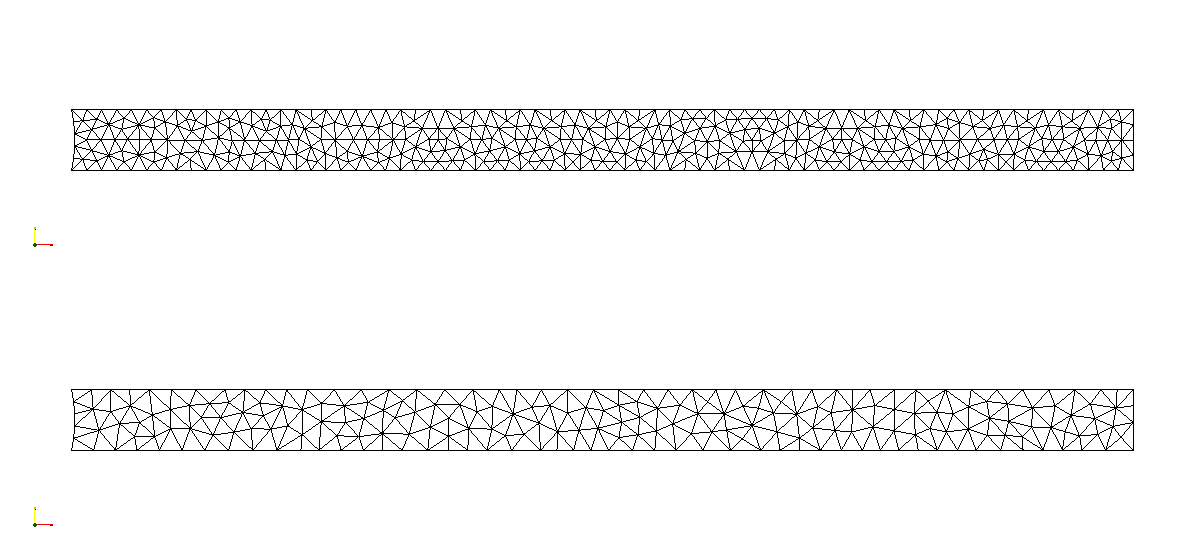}
\caption{Solid mesh convergence study consisting of two different meshes of 354, and 747 elements.}
\label{fig:FSI2SolidMeshConvergence}
\end{figure}

Before analyzing the results, a mesh convergence study is conducted, examining the impact of fluid-solid mesh and solid mesh on the problem's solution. 
This study helps identify optimal mesh settings for the simulations. 
It's essential to have an adequate mesh to capture the flow and deformation behavior accurately. 
Figure \ref{fig:FSI2FluidMeshConvergence} demonstrates the fluid mesh convergence study, where we have meshes of 20 000, 30 000, and 40 000 elements respectively, while Figure \ref{fig:FSI2SolidMeshConvergence} shows the solid mesh convergence study, where we have meshes of 354 and 747 elements.

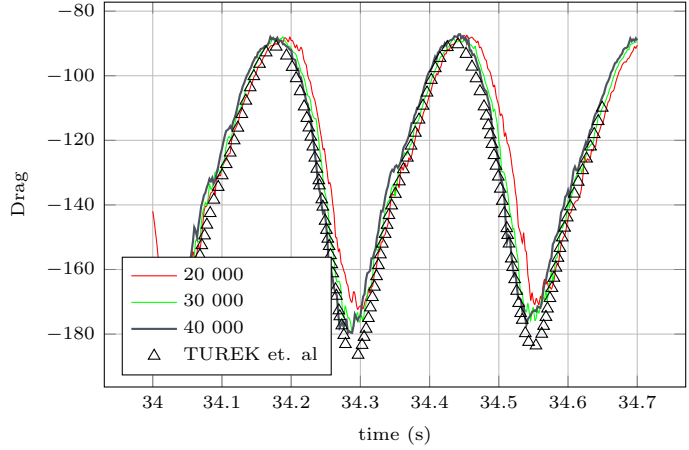
\begin{figure}[!h]
        \centering
        \begin{tikzpicture}[trim axis left, trim axis right]
                \begin{axis}[
                scale=0.8,
                transform shape,
                width=.7\textwidth,
                height=.5\textwidth,
                label style={font=\scriptsize},
                tick label style={font=\scriptsize},
                legend style={font=\scriptsize},
                clip=true,
                grid=both,
                legend cell align=left,
                legend pos=south west,
                grid=major,
                xlabel={time (s)},
                ylabel={Drag}  
                ]                 
                
                \legend{20 000, 30 000, 40 000, TUREK et. al}
                       
                 \addplot[ 
                draw=red,
                smooth
                ]
                table[
                x index=0,
                y index=1
                ]
                {data/FSI2Drag20000.dat};
                
                \addplot[ 
                draw=green,
                smooth
                ]
                table[
                x index=0,
                y index=1
                ]
                {data/FSI2Drag30000.dat};  
                
                \addplot[ 
                mark=none,thick,mygray1
                ]
                table[
                x index=0,
                y index=1
                ]
                {data/FSI2Drag40000.dat};             
                
                \addplot[ 
                only marks,
    		    style={solid, fill=gray},
    		    mark=triangle,
   			    mark size=2.5pt
                ]
                 table[
                x index=0,
                y index=1
                ]
                {data/FSI2DragBenchmark.dat};    
               
                \end{axis}
        \end{tikzpicture}
        \caption{FSI 2 drag over the cylinder and membrane versus time for different number of elements and comparison with that from Turek et. al.}
		\label{fig:FSI2Drag}
\end{figure}

\begin{figure}[!h]
        \centering
        \begin{tikzpicture}[trim axis left, trim axis right]
                \begin{axis}[
                scale=0.8,
                transform shape,
                width=.7\textwidth,
                height=.5\textwidth,
                label style={font=\scriptsize},
                tick label style={font=\scriptsize},
                legend style={font=\scriptsize},
                clip=true,
                grid=both,
                legend cell align=left,
                legend pos=south west,
                grid=major,
                xlabel={time (s)},
                ylabel={Lift}  
                ]                 
                
                \legend{20 000, 30 000, 40 000, TUREK et. al}
                       
                 \addplot[ 
                draw=red,
                smooth
                ]
                table[
                x index=0,
                y index=1
                ]
                {data/FSI2Lift20000.dat};
                
                \addplot[ 
                draw=green,
                smooth
                ]
                table[
                x index=0,
                y index=1
                ]
                {data/FSI2Lift30000.dat};  
                
                \addplot[ 
                mark=none,thick,mygray1
                ]
                table[
                x index=0,
                y index=1
                ]
                {data/FSI2Lift40000.dat};             
                
                \addplot[ 
                only marks,
    		    style={solid, fill=gray},
    		    mark=triangle,
   			    mark size=2.5pt
                ]
                 table[
                x index=0,
                y index=1
                ]
                {data/FSI2LiftBenchmark.dat};    
               
                \end{axis}
        \end{tikzpicture}
        \caption{FSI 2 lift over the cylinder and membrane versus time for different number of elements and comparison with that from Turek et. al.}
		\label{fig:FSI2Lift}
\end{figure}
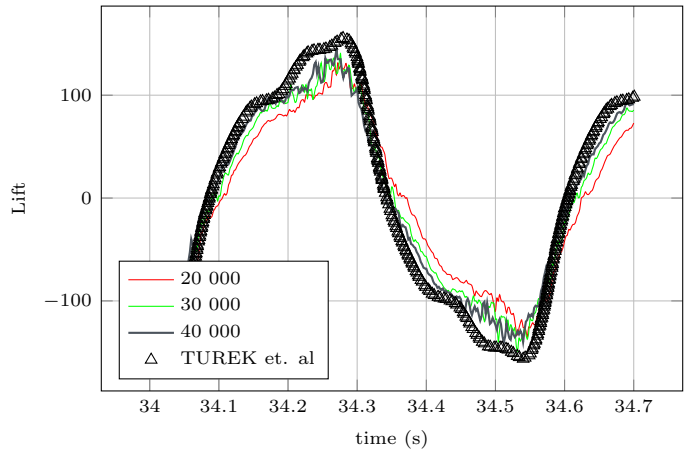

The fluid-solid mesh also plays an important role in the accuracy and the precision of the drag and lift a shown in figures \ref{fig:FSI2Drag}, and \ref{fig:FSI2Lift}.
It is important to note however that the elements are distributed in a way that the displacement of the solid converges rapidly to the final solution for rather coarse meshes, because of the way the elements are distributed close to the interface.
It can be seen from the graphs that the frequency of the oscillations and the amplitude correlate with that of the benchmark.

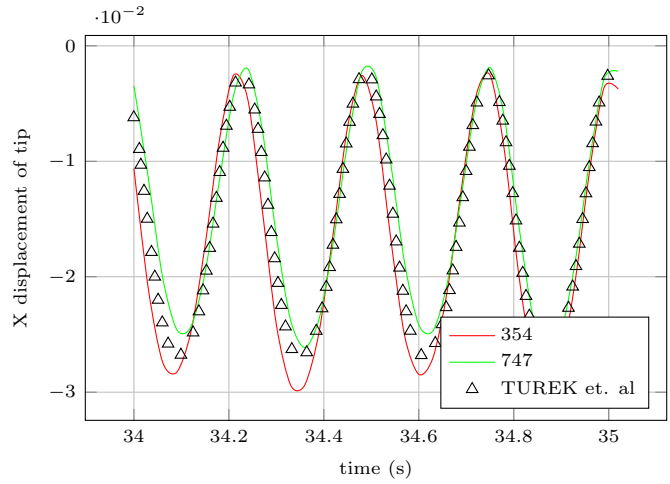
\begin{figure}[!h]
        \centering
        \begin{tikzpicture}[trim axis left, trim axis right]
                \begin{axis}[
                scale=0.8,
                transform shape,
                width=.7\textwidth,
                height=.5\textwidth,
                label style={font=\scriptsize},
                tick label style={font=\scriptsize},
                legend style={font=\scriptsize},
                clip=true,
                grid=both,
                legend cell align=left,
                legend pos=south east,
                grid=major,
                xlabel={time (s)},
                ylabel={X displacement of tip}  
                ]                 
                
                \legend{354, 747, TUREK et. al}
                
                \addplot[ 
                draw=red,
                smooth
                ]
                table[
                x index=0,
                y index=1
                ]
                {data/FSI2UxHS1.dat};
                
                \addplot[ 
                draw=green,
                smooth
                ]
                table[
                x index=0,
                y index=1
                ]
                {data/FSI2UxHS2.dat}; 
                
                \addplot[ 
                only marks,
    		    style={solid, fill=gray},
    		    mark=triangle,
   			    mark size=2.5pt
                ]
                 table[
                x index=0,
                y index=1
                ]
                {data/FSI2UxBenchmark.dat};    
               
                \end{axis}
        \end{tikzpicture}
        \caption{FSI 2 X displacement of tip versus time and comparison with that from Turek et. al for different solid meshes.}
		\label{fig:FSI2Uxsolidmesh}
\end{figure}

The solid mesh was found to correlate the most with the X displacement of the beam, and not so much with the Y displacement of the beam, which converges to the benchmark values for rather coarse meshes.
Figure \ref{fig:FSI2Uxsolidmesh} show the dependency of the X displacement of the beam on the number of elements of the solid mesh.

\begin{figure}[!h]
        \centering
        \begin{tikzpicture}[trim axis left, trim axis right]
                \begin{axis}[
                scale=0.8,
                transform shape,
                width=.7\textwidth,
                height=.5\textwidth,
                label style={font=\scriptsize},
                tick label style={font=\scriptsize},
                legend style={font=\scriptsize},
                clip=true,
                grid=both,
                legend cell align=left,
                legend pos=south east,
                grid=major,
                xlabel={time (s)},
                ylabel={Y displacement of tip}  
                ]                 
                
                \legend{delta t 1, delta t 2, delta t 3 TUREK et. al}
                       
                \addplot[ 
                draw=red,
                smooth
                ]
                table[
                x index=0,
                y index=1
                ]
                {data/FSI2Uydeltat1.dat};
                
                \addplot[ 
                draw=green,
                smooth
                ]
                table[
                x index=0,
                y index=1
                ]
                {data/FSI2Uydeltat2.dat};
                
                 \addplot[ 
                draw=black,
                smooth
                ]
                table[
                x index=0,
                y index=1
                ]
                {data/FSI2Uy.dat};

                \addplot[ 
                only marks,
    		    style={solid, fill=gray},
    		    mark=triangle,
   			    mark size=2.5pt
                ]
                 table[
                x index=0,
                y index=1
                ]
                {data/FSI2UyBenchmark.dat};    
               
                \end{axis}
        \end{tikzpicture}
        \caption{FSI 2 Y displacement of tip versus time and comparison with that from Turek et. al for different time steps.}
		\label{fig:FSI2Uydeltat}
\end{figure}
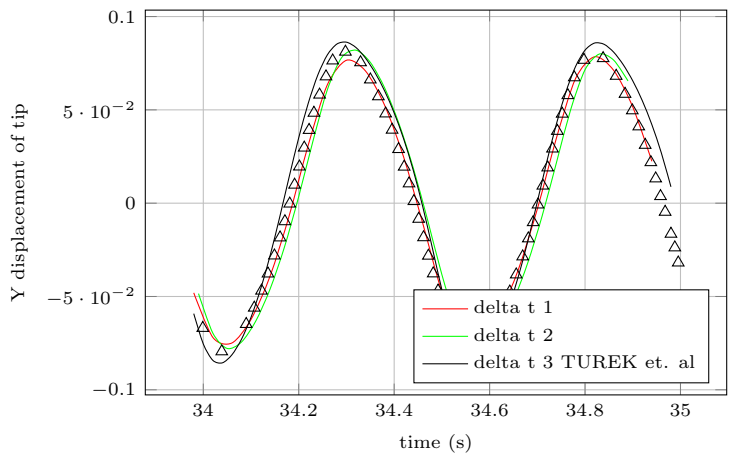

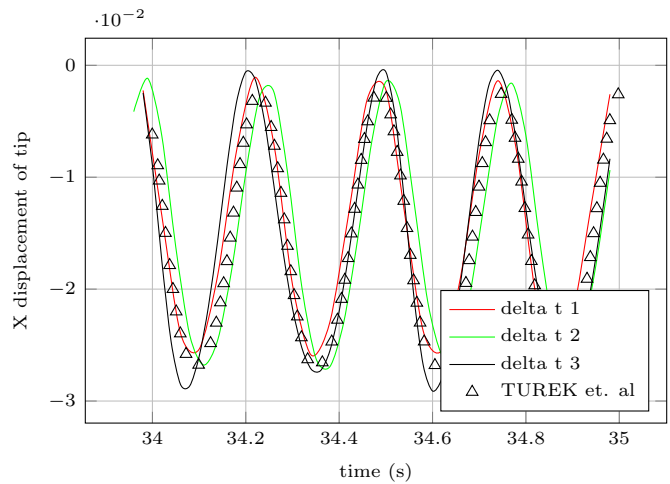
\begin{figure}[!h]
        \centering
        \begin{tikzpicture}[trim axis left, trim axis right]
                \begin{axis}[
                scale=0.8,
                transform shape,
                width=.7\textwidth,
                height=.5\textwidth,
                label style={font=\scriptsize},
                tick label style={font=\scriptsize},
                legend style={font=\scriptsize},
                clip=true,
                grid=both,
                legend cell align=left,
                legend pos=south east,
                grid=major,
                xlabel={time (s)},
                ylabel={X displacement of tip}  
                ]                 
                
                \legend{delta t 1, delta t 2, delta t 3, TUREK et. al}
                
                \addplot[ 
                draw=red,
                smooth
                ]
                table[
                x index=0,
                y index=1
                ]
                {data/FSI2Uxdeltat1.dat};
                
                \addplot[ 
                draw=green,
                smooth
                ]
                table[
                x index=0,
                y index=1
                ]
                {data/FSI2Uxdeltat2.dat}; 
                       
                \addplot[ 
                draw=black,
                smooth
                ]
                table[
                x index=0,
                y index=1
                ]
                {data/FSI2Ux.dat};

                \addplot[ 
                only marks,
    		    style={solid, fill=gray},
    		    mark=triangle,
   			    mark size=2.5pt
                ]
                 table[
                x index=0,
                y index=1
                ]
                {data/FSI2UxBenchmark.dat};    
               
                \end{axis}
        \end{tikzpicture}
        \caption{FSI 2 X displacement of tip versus time and comparison with that from Turek et. al for different time steps.}
		\label{fig:FSI2Uxdeltat}
\end{figure}

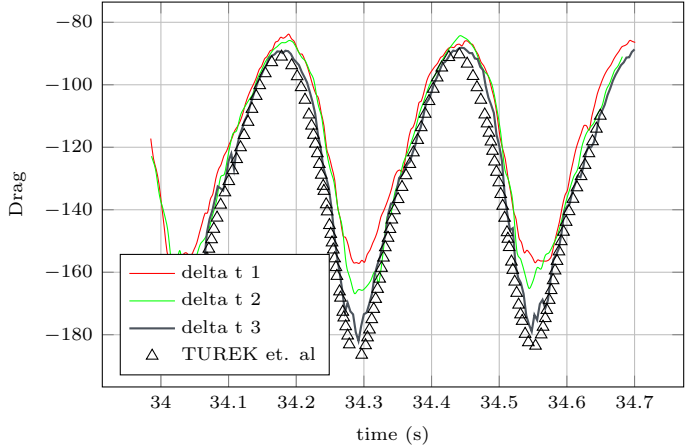
\begin{figure}[!h]
        \centering
        \begin{tikzpicture}[trim axis left, trim axis right]
                \begin{axis}[
                scale=0.8,
                transform shape,
                width=.7\textwidth,
                height=.5\textwidth,
                label style={font=\scriptsize},
                tick label style={font=\scriptsize},
                legend style={font=\scriptsize},
                clip=true,
                grid=both,
                legend cell align=left,
                legend pos=south west,
                grid=major,
                xlabel={time (s)},
                ylabel={Drag}  
                ]                 
                
                \legend{delta t 1, delta t 2, delta t 3, TUREK et. al}
                       
                 \addplot[ 
                draw=red,
                smooth
                ]
                table[
                x index=0,
                y index=1
                ]
                {data/FSI2Dragdeltat1.dat};
                
                \addplot[ 
                draw=green,
                smooth
                ]
                table[
                x index=0,
                y index=1
                ]
                {data/FSI2Dragdeltat2.dat};  
                
                \addplot[ 
                mark=none,thick,mygray1
                ]
                table[
                x index=0,
                y index=1
                ]
                {data/FSI2Dragdeltat3.dat};             
                
                \addplot[ 
                only marks,
    		    style={solid, fill=gray},
    		    mark=triangle,
   			    mark size=2.5pt
                ]
                 table[
                x index=0,
                y index=1
                ]
                {data/FSI2DragBenchmark.dat};    
               
                \end{axis}
        \end{tikzpicture}
        \caption{FSI 2 drag over the cylinder and membrane versus time for different time steps and comparison with that from Turek et. al.}
		\label{fig:FSI2Dragdeltat}
\end{figure}

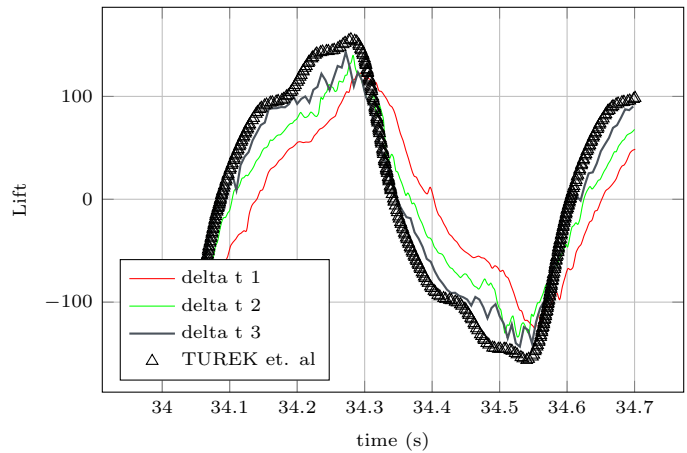
\begin{figure}[!h]
        \centering
        \begin{tikzpicture}[trim axis left, trim axis right]
                \begin{axis}[
                scale=0.8,
                transform shape,
                width=.7\textwidth,
                height=.5\textwidth,
                label style={font=\scriptsize},
                tick label style={font=\scriptsize},
                legend style={font=\scriptsize},
                clip=true,
                grid=both,
                legend cell align=left,
                legend pos=south west,
                grid=major,
                xlabel={time (s)},
                ylabel={Lift}  
                ]                 
                
                \legend{delta t 1, delta t 2, delta t 3, TUREK et. al}
                       
                 \addplot[ 
                draw=red,
                smooth
                ]
                table[
                x index=0,
                y index=1
                ]
                {data/FSI2Liftdeltat1.dat};
                
                \addplot[ 
                draw=green,
                smooth
                ]
                table[
                x index=0,
                y index=1
                ]
                {data/FSI2Liftdeltat2.dat};  
                
                \addplot[ 
                mark=none,thick,mygray1
                ]
                table[
                x index=0,
                y index=1
                ]
                {data/FSI2Liftdeltat3.dat};             
                
                \addplot[ 
                only marks,
    		    style={solid, fill=gray},
    		    mark=triangle,
   			    mark size=2.5pt
                ]
                 table[
                x index=0,
                y index=1
                ]
                {data/FSI2LiftBenchmark.dat};    
               
                \end{axis}
        \end{tikzpicture}
        \caption{FSI 2 lift over the cylinder and membrane versus time for different time steps and comparison with that from Turek et. al.}
		\label{fig:FSI2Liftdeltat}
\end{figure}

The choice of time step in simulations is crucial, as it influences the results' accuracy and stability. 
In this benchmark, different time steps are explored, and their impact on displacement, drag, and lift is evaluated. 
Figures \ref{fig:FSI2Uydeltat}, \ref{fig:FSI2Uxdeltat}, \ref{fig:FSI2Dragdeltat}, and \ref{fig:FSI2Liftdeltat} depict how the time step affects these parameters, with values equal to 0.004, 0.003, and 0.002 seconds.

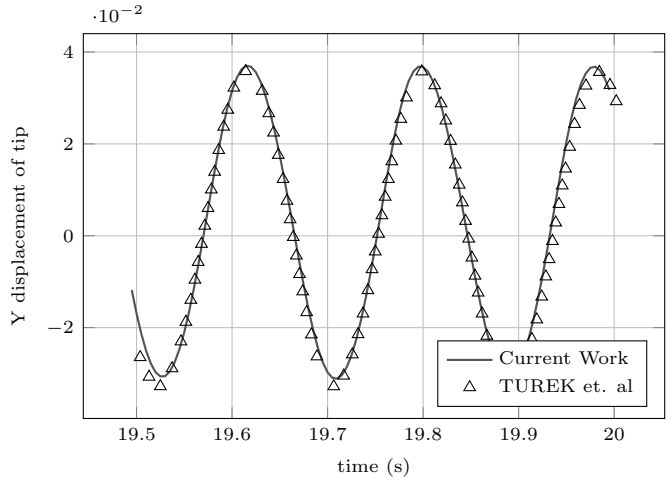
\begin{figure}[!h]
        \centering
        \begin{tikzpicture}[trim axis left, trim axis right]
                \begin{axis}[
                scale=0.8,
                transform shape,
                width=.7\textwidth,
                height=.5\textwidth,
                label style={font=\scriptsize},
                tick label style={font=\scriptsize},
                legend style={font=\scriptsize},
                clip=true,
                grid=both,
                legend cell align=left,
                legend pos=south east,
                grid=major,
                xlabel={time (s)},
                ylabel={Y displacement of tip}  
                ]                 
                
                \legend{Current Work, TUREK et. al}
                       
                \addplot[ 
                mark=none,thick,mygray1
                ]
                table[
                x index=0,
                y index=1
                ]
                {data/FSI3Uy.dat};

                \addplot[ 
                only marks,
    		    style={solid, fill=gray},
    		    mark=triangle,
   			    mark size=2.5pt
                ]
                 table[
                x index=0,
                y index=1
                ]
                {data/FSI3UyBenchmark.dat};    
               
                \end{axis}
        \end{tikzpicture}
        \caption{FSI 3 Y displacement of tip versus time and comparison with that from Turek et. al.}
		\label{fig:FSI3Uy}
\end{figure}

\begin{figure}[!h]
        \centering
        \begin{tikzpicture}[trim axis left, trim axis right]
                \begin{axis}[
                scale=0.8,
                transform shape,
                width=.7\textwidth,
                height=.5\textwidth,
                label style={font=\scriptsize},
                tick label style={font=\scriptsize},
                legend style={font=\scriptsize},
                clip=true,
                grid=both,
                legend cell align=left,
                legend pos=south east,
                grid=major,
                xlabel={time (s)},
                ylabel={X displacement of tip}  
                ]                 
                
                \legend{Current Work, TUREK et. al}
                       
                \addplot[ 
                draw=black,
                smooth
                ]
                table[
                x index=0,
                y index=1
                ]
                {data/FSI3Ux.dat};

                \addplot[ 
                only marks,
    		    style={solid, fill=gray},
    		    mark=triangle,
   			    mark size=2.5pt
                ]
                 table[
                x index=0,
                y index=1
                ]
                {data/FSI3UxBenchmark.dat};    
               
                \end{axis}
        \end{tikzpicture}
        \caption{FSI 3 X displacement of tip versus time and comparison with that from Turek et. al.}
		\label{fig:FSI3Ux}
\end{figure}
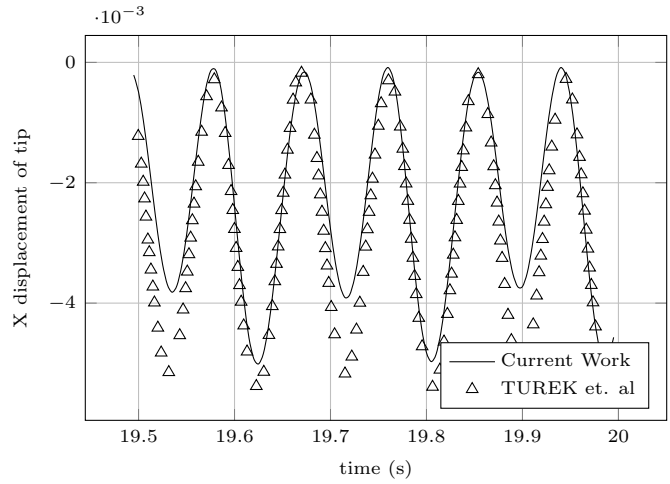

FSI3 is run until 20 seconds with a time step of $\Delta t=0.0005$.
The Y, and X displacements for FSI3 are plotted versus time in Figure  \ref{fig:FSI3Uy}, and \ref{fig:FSI3Ux} respectively.
It can be seen from the graph that the frequency of the oscillations and the amplitude correlate with that of the benchmark.

\begin{figure}[!h]
\centering
    \begin{subfigure}{0.7\linewidth}
        \includegraphics[width=\linewidth]{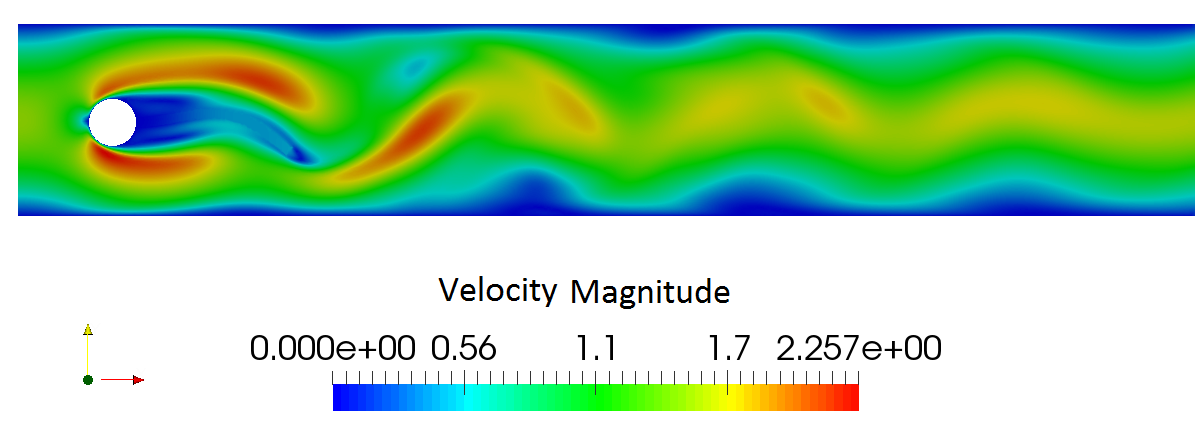}
        \label{fig:43a}
    \end{subfigure}
\\
    \begin{subfigure}{0.7\linewidth}
        \includegraphics[width=\linewidth]{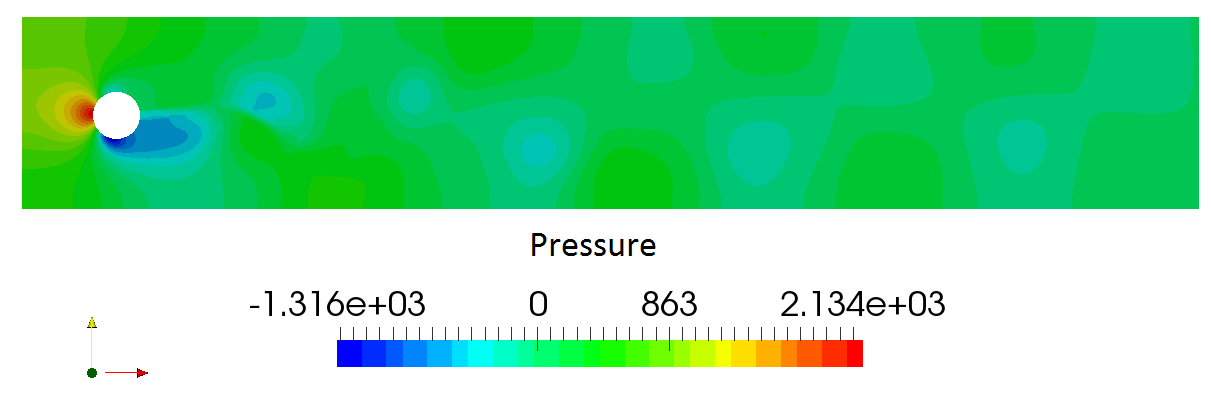}
        \label{fig:43b}
    \end{subfigure}
\\
    \begin{subfigure}{0.7\linewidth}
        \includegraphics[width=\linewidth]{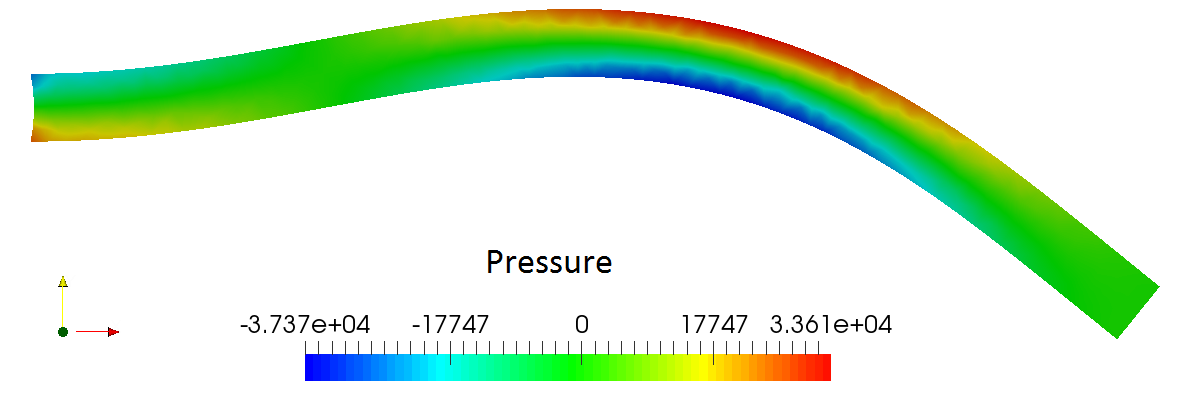}
        \label{fig:43c}
    \end{subfigure}
\caption{Velocity and pressure contour for the fluid-solid and solid domains for FSI2.}
\label{fig:FSI2Contours}
\end{figure}

\begin{figure}[!h]
\centering
    \begin{subfigure}{0.7\linewidth}
        \includegraphics[width=\linewidth]{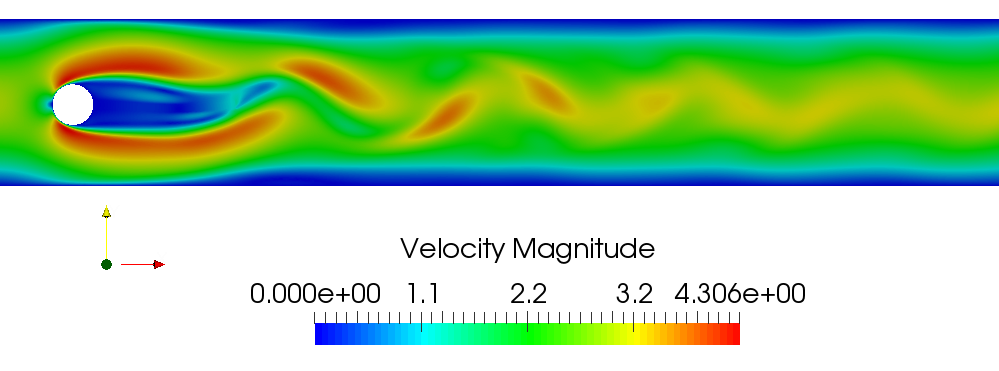}
        \label{fig:44a}
    \end{subfigure} 
\\
    \begin{subfigure}{0.7\linewidth}
        \includegraphics[width=\linewidth]{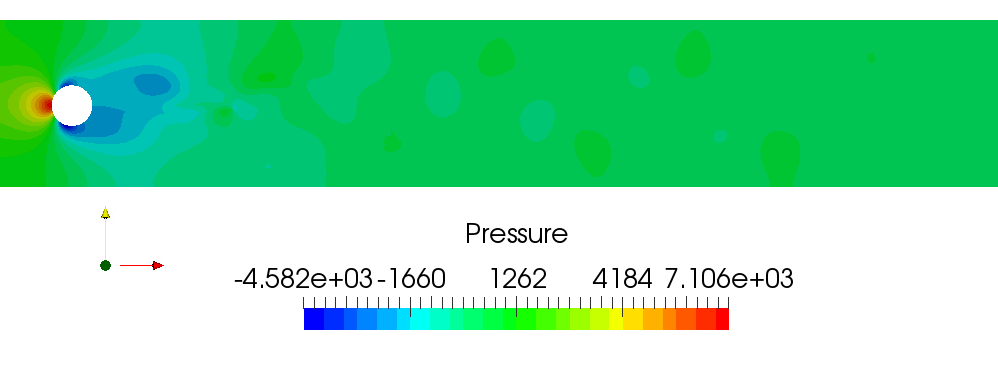}
        \label{fig:44b}
    \end{subfigure}
\\
    \begin{subfigure}{0.7\linewidth}
        \includegraphics[width=\linewidth]{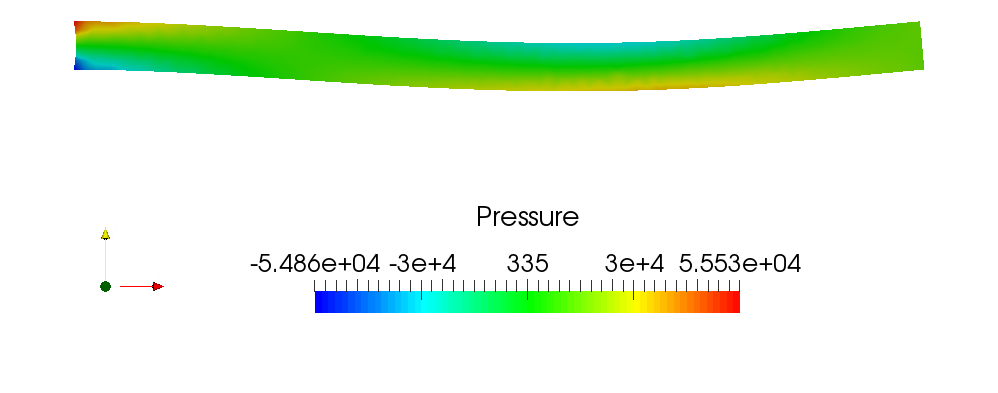}
        \label{fig:44c}
    \end{subfigure}
\caption{Velocity and pressure contour for the fluid-solid and solid domains for FSI3.}
\label{fig:FSI3Contours}
\end{figure}

The velocity and pressure contours for both FSI2 and FSI3 are displayed in Figures \ref{fig:FSI2Contours} and \ref{fig:FSI3Contours}, respectively.
In FSI2 and FSI3, the simulations are run for specific durations and with different parameters. 
The results are presented for both cases, including the Y and X displacements, drag, and lift versus time. 
These results are compared to the benchmark values to assess the fidelity of the simulations.
Both cases exhibit oscillatory behavior, and it's evident from the graphs that the frequency and amplitude of these oscillations correlate well with the benchmark data.
The FSI2 and FSI3 variations of Turek's benchmark offer valuable insights into the accuracy and stability of FSI simulations under different parameter configurations. 
The convergence studies provide guidance on mesh and time step selection, and the overall results showcase the solver's ability to reproduce the benchmark's complex flow phenomena.
The ability to capture the complex flow characteristics and validate against the benchmark data is crucial for FSI simulations, making the Turek FSI benchmark a significant reference in the field of computational fluid dynamics. 
It enables the assessment of solver performance and optimization of simulation parameters.

\subsection{Pillar in a Laminar cross flow}

In our investigation, we perform a comprehensive three-dimensional (3D) fully coupled two-way Fluid-Structure Interaction (FSI) simulation. 
Our simulation is aimed at validating the framework by replicating an experimental setup detailed in \cite{axtmann2016investigation}. This experimental arrangement features a transparent basin with dimensions of 3 meters in length, 2.5 meters in width, and 0.4 meters in height, which is filled with glycerin.
The experimental basin includes a versatile support structure that can attain velocities of up to 1 meter per second. 
To capture crucial data, particularly focusing on the bending behavior of the beam, a high-speed camera is employed.
The choice of materials is a crucial aspect of our simulation. Notably, we use glycerin as the working fluid and silicone for the construction of the flexible cylinder. 
The diameter of the cylinder is 20 millimeters, with a height of 200 millimeters, resulting in an aspect ratio of 10. 
The decision to employ glycerin and silicone is based on their closely matched densities, allowing us to reasonably neglect buoyancy forces in our simulations.
The material properties of glycerin and silicone are presented in Table \ref{table:51} for reference. 
Glycerin exhibits Newtonian behavior, with a density ($\rho_f$) of 1220 kg/m$^3$ and a dynamic viscosity ($\mu_f$) of 1 kg/ms. Silicone, on the other hand, is modeled using a Neo-Hookean model with a density ($\rho_s$) of 1030 kg/m$^3$, an elastic modulus ($E$) of 1.23 MPa, and Poisson's ratio ($\nu$) of 0.3.

The material properties of both fluid and solid are summarized in Table \ref{table:51}.

\begin{table}[!h]
\centering
\begin{tabular}{|p{3cm}|p{3cm}||p{3cm}|p{3cm}|}
 \hline
 \multicolumn{2}{|c|}{Fluid (Glycerin)} & \multicolumn{2}{|c|}{Solid(Silicone)}\\
 \hline
$\rho_f$ & 1220 Kg/m$^3$ & $\rho_s$ & 1030 Kg/m$^3$\\
$\mu_f$ & 1 Kg/ms & $\mu_s$ & 0.473 MPa\\
&&$E$& 1.23 MPa\\
&&$\nu$& 0.3 \\
Model & Newtonian & Model &  Neo-Hookean\\
\hline
\end{tabular} 
\caption{Fluid and solid properties for the pillar in a Laminar cross flow problem.}
\label{table:51}
\end{table}

\begin{figure}[!h]
\centering
\begin{minipage}{0.5\textwidth}
   \includegraphics[width=1\linewidth]{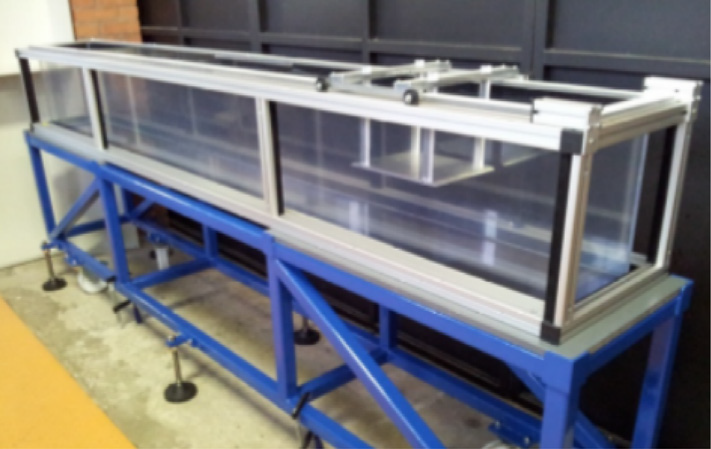}
  \end{minipage}%
  \hspace*{\fill}   
  \begin{minipage}{0.5\textwidth}
   \includegraphics[width=1\linewidth]{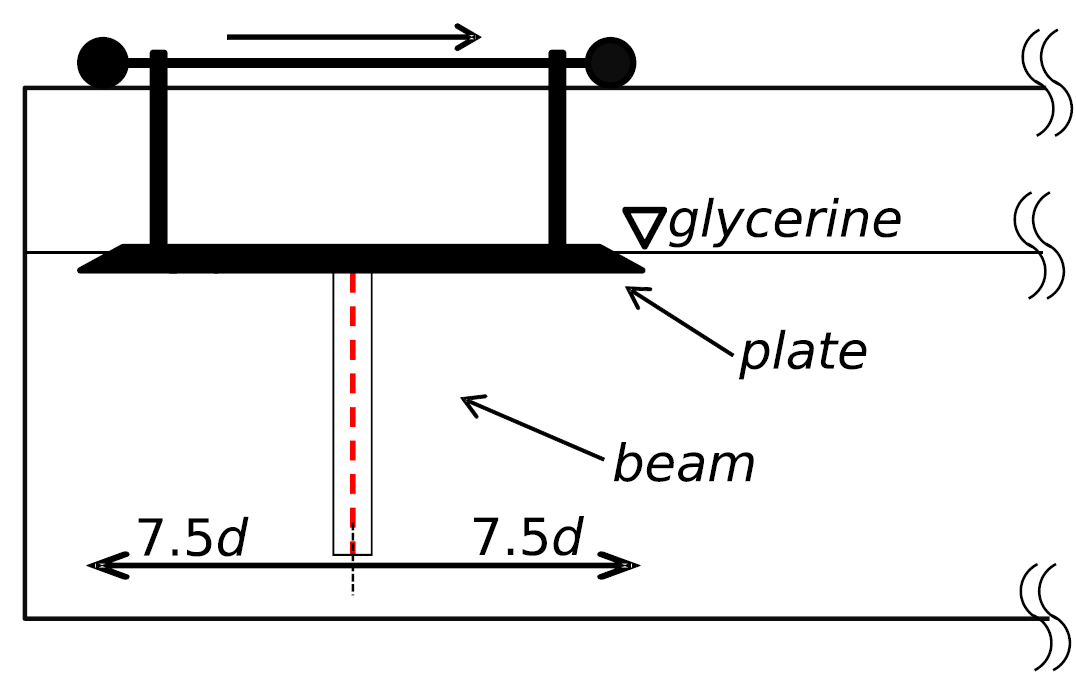}
  \end{minipage}%
  \caption{Experimental setup and schematic from Axtmann et. al\cite{axtmann2016investigation}} \label{fig:Experimental}
\end{figure}

\begin{figure}[!h]
\centering
   \includegraphics[width=1\linewidth]{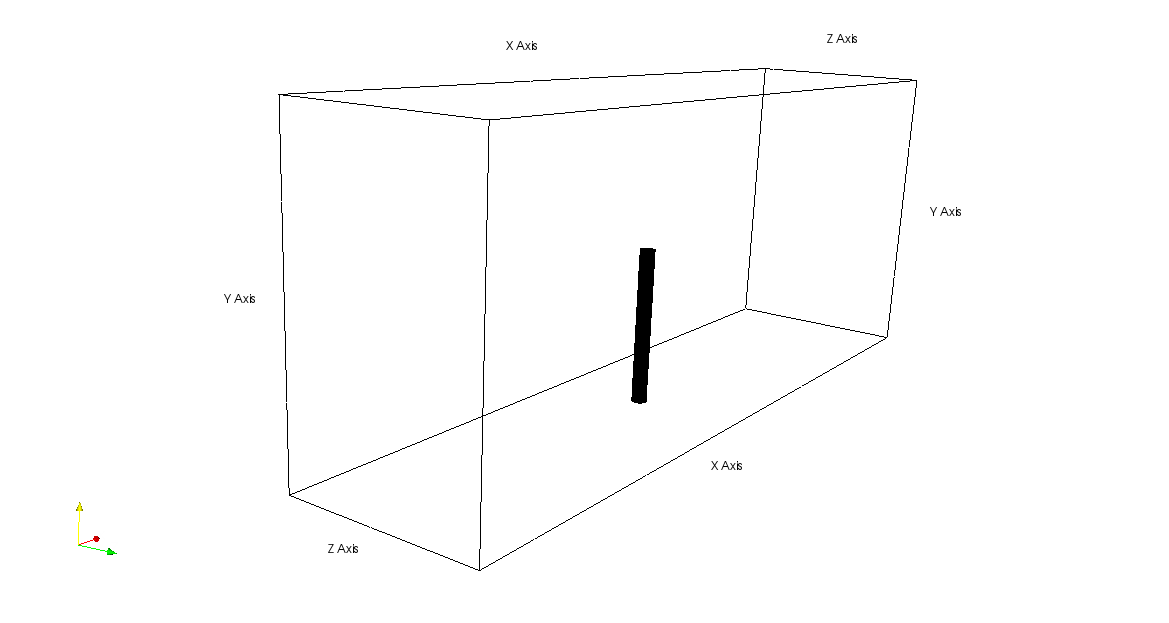}
\caption{Problem set up for bending beam 2}
\label{fig:SetupCylinder}
\end{figure}

\begin{figure}[!h]
\centering
\begingroup  
\setlength{\tabcolsep}{0pt}
\renewcommand{\arraystretch}{0}
\begin{tabular}{c c}
    \begin{subfigure}{0.51\linewidth}
        \includegraphics[width=\linewidth]{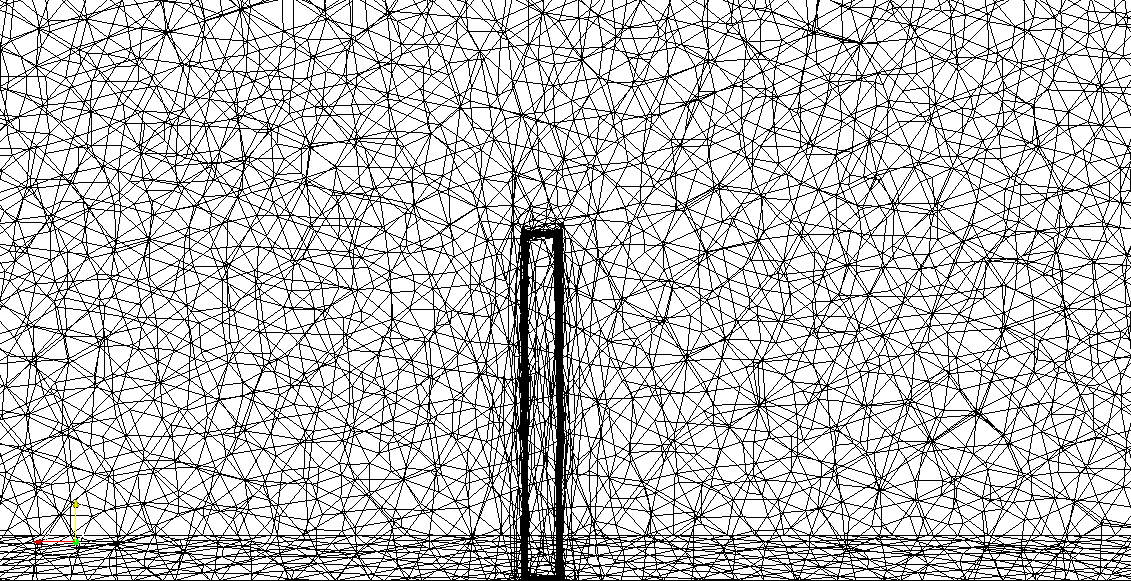}
        \caption{Fluid-Solid Mesh}
        \label{fig:51a}
    \end{subfigure}
    \begin{subfigure}{0.49\linewidth}
        \includegraphics[width=\linewidth]{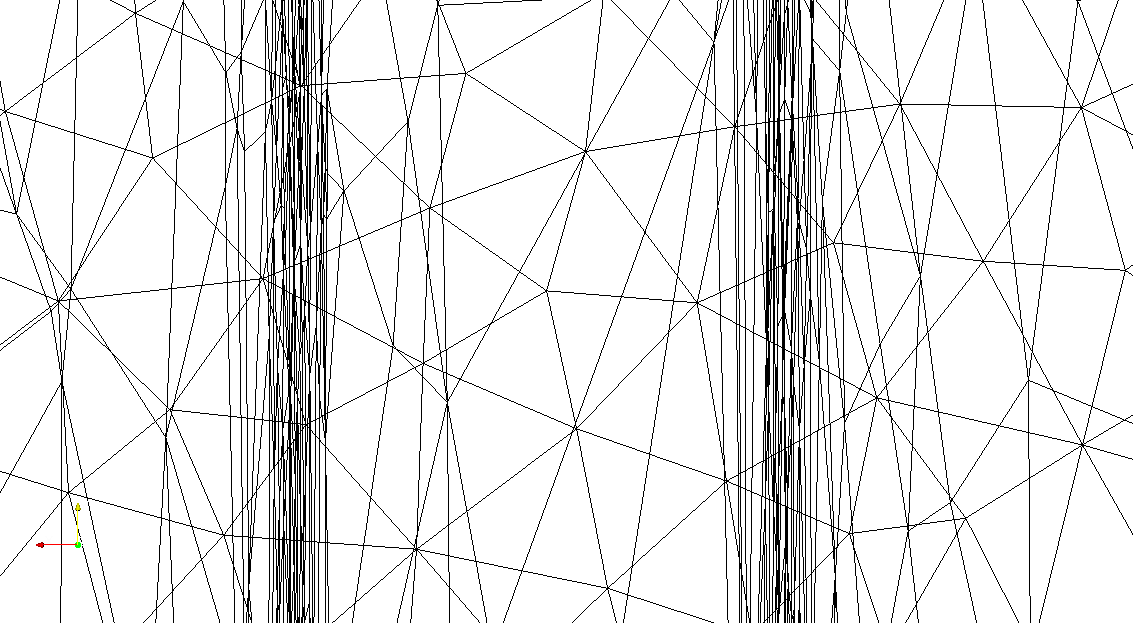}
        \caption{Anisotropic mesh at the interface}
        \label{fig:51b}
    \end{subfigure}
\\
    \begin{subfigure}{0.2\linewidth}
        \includegraphics[width=\linewidth]{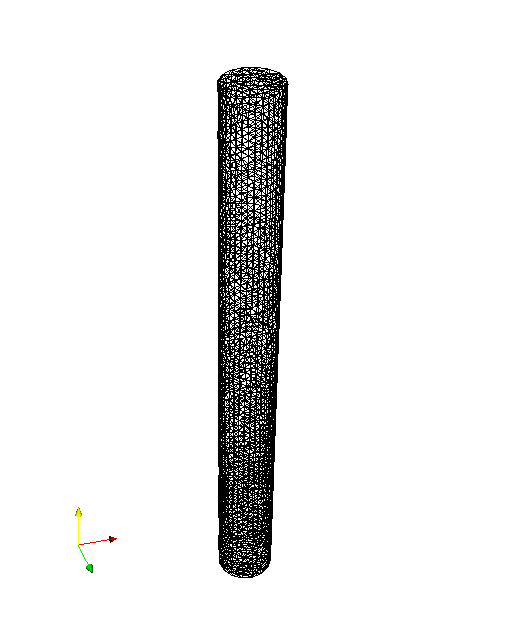}
        \caption{Solid Mesh}
        \label{fig:51c}
    \end{subfigure}
\end{tabular}
\endgroup
\caption{Fluid-solid and solid meshes for the pillar in a Laminar cross flow.}
\label{fig:5ProblemDescription}
\end{figure}

The next step involves configuring our numerical simulation to mirror the experimental setup shown in figure \ref{fig:Experimental}. 
To achieve this, we create a virtual fluid basin, represented as a 3D fluid-solid mesh. 
Additionally, a distinct solid mesh is generated to emulate the polymer structure, which is submerged within the fluid-solid mesh.
In our simulation, we deviate from the experimental procedure where the solid is physically moved. 
Instead, we impose a velocity at the inlet, matching the velocity of the plate where the polymer cylinder is clamped. 
Our simulation setup includes a zero gauge pressure outlet, perfect slip conditions on the virtual domain walls, and zero slip conditions applied to both the cylinder and the plate.
We proceed with the simulation at a Reynolds number of 12, closely mimicking the experimental conditions. 
The geometric layout is visually represented in Figure \ref{fig:SetupCylinder}, offering a clear understanding of the setup. 
For an in-depth perspective, Figures \ref{fig:51a}, \ref{fig:51b}, and \ref{fig:51c} present the fluid-solid mesh, a magnified view highlighting the anisotropic elements near the interface, and the solid structure mesh, respectively.
Our numerical simulation consists of 35,000 solid elements and is capped at 150,000 elements within the fluid-solid mesh. 
We utilize a time step of 0.001 seconds, conducting the simulation for a duration of 2 seconds.

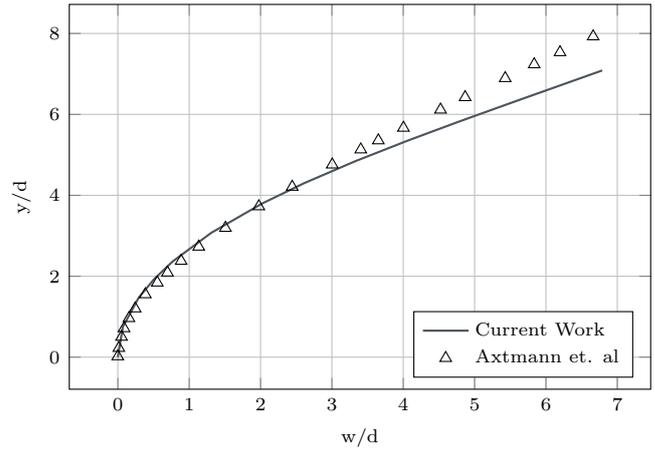
\begin{figure}[p]
        \centering
        \begin{tikzpicture}[trim axis left, trim axis right]
                \begin{axis}[
                scale=0.8,
                transform shape,
                width=.7\textwidth,
                height=.5\textwidth,
                label style={font=\scriptsize},
                tick label style={font=\scriptsize},
                legend style={font=\scriptsize},
                clip=true,
                grid=both,
                legend cell align=left,
                legend pos=south east,
                grid=major,
                xlabel={w/d},
                ylabel={y/d}  
                ]                 
                
                \legend{Current Work, Axtmann et. al}
                       
                \addplot[ 
                mark=none,thick,mygray1
                ]
                table[
                x index=0,
                y index=1
                ]
                {data/BC.dat};

                \addplot[ 
                only marks,
    		    style={solid, fill=gray},
    		    mark=triangle,
   			    mark size=2.5pt
                ]
                 table[
                x index=0,
                y index=1
                ]
                {data/BCB.dat};    
               
                \end{axis}
        \end{tikzpicture}
        \caption{Bending line of the beam at the final time t compared with the work of that of Axtmann et. al\cite{axtmann2016investigation}.}
		\label{fig:BC}
\end{figure}

\begin{figure}[p]
 \centering
        \begin{tikzpicture}[trim axis left, trim axis right]
                \begin{axis}[
                scale=0.8,
                transform shape,
                width=.7\textwidth,
                height=.5\textwidth,
                label style={font=\scriptsize},
                tick label style={font=\scriptsize},
                legend style={font=\scriptsize},
                clip=true,
                grid=both,
                legend cell align=left,
                legend pos=south east,
                grid=major,
                xlabel={time (s)},
                ylabel={y displacement of tip},  
                ]                 
                \addplot[ 
                mark=none,thick,mygray1
                ]
                table[
                x index=0,
                y index=1
                ]
                {data/BCUy.dat}; 
                \end{axis}
        \end{tikzpicture}
\caption{x and y displacement of top center versus time for pillar in a Laminar cross flow }
		\label{fig:BCUy}
\end{figure}
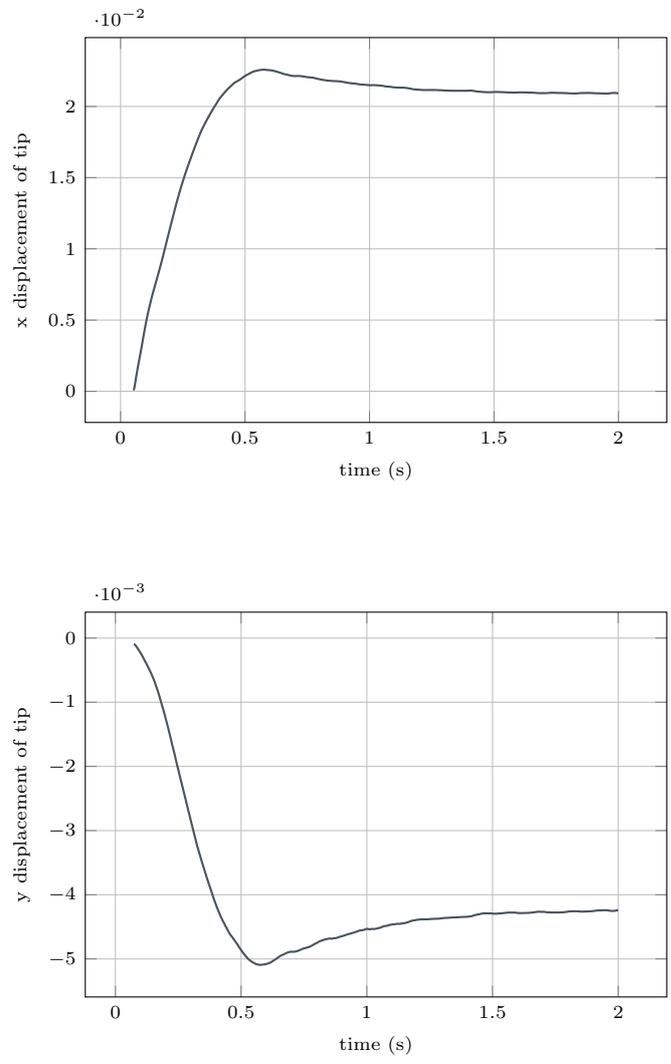

\begin{figure}[p]
\centering
    \begin{subfigure}{0.45\linewidth}
        \includegraphics[width=\linewidth]{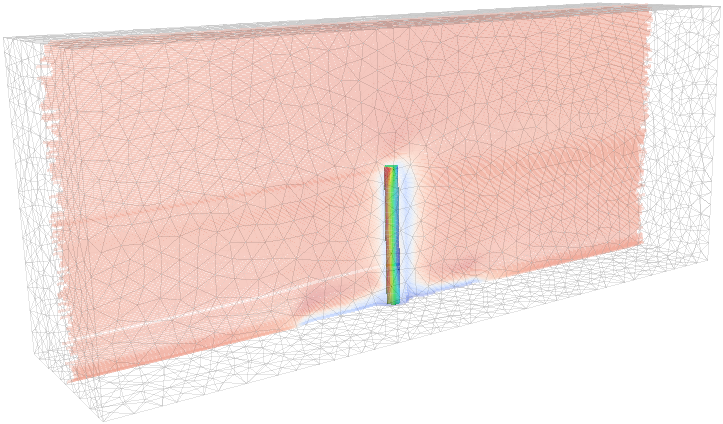}
        \caption{t=0}
        \label{fig:52a}
    \end{subfigure}
    \begin{subfigure}{0.45\linewidth}
        \includegraphics[width=\linewidth]{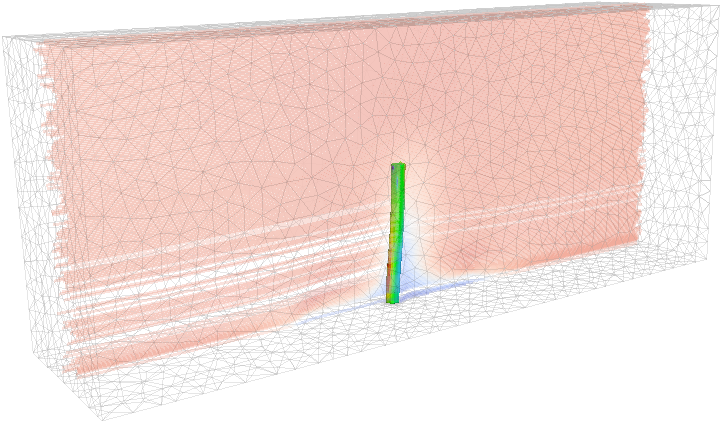}
        \caption{t=0.1}
        \label{fig:52b}
    \end{subfigure}
\\
    \begin{subfigure}{0.45\linewidth}
        \includegraphics[width=\linewidth]{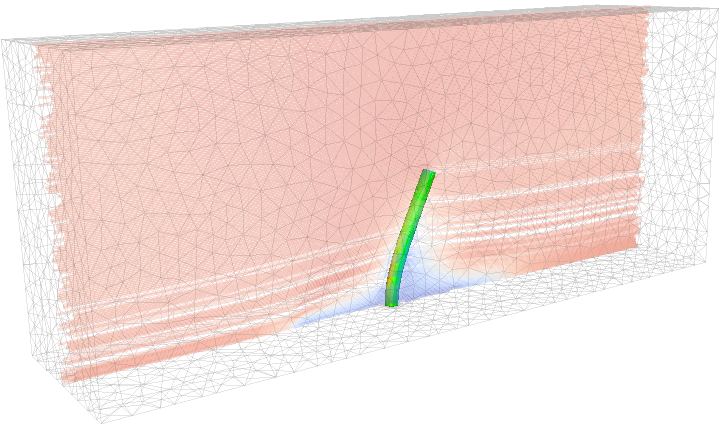}
        \caption{t=0.2}
        \label{fig:52c}
    \end{subfigure}
    \begin{subfigure}{0.45\linewidth}
        \includegraphics[width=\linewidth]{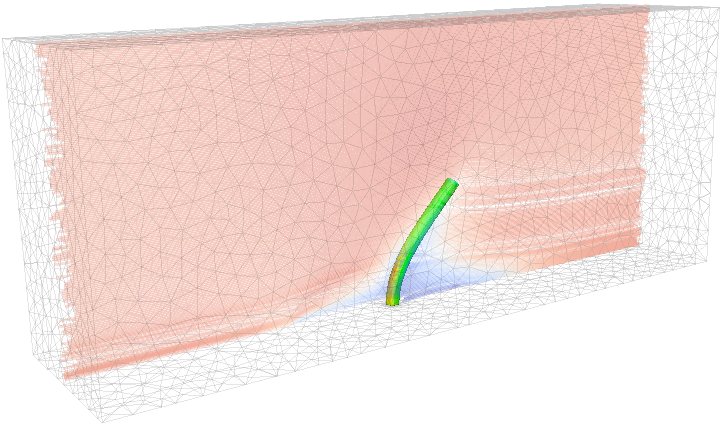}
        \caption{t=0.3}
        \label{fig:52d}
    \end{subfigure}
\\
    \begin{subfigure}{0.45\linewidth}
        \includegraphics[width=\linewidth]{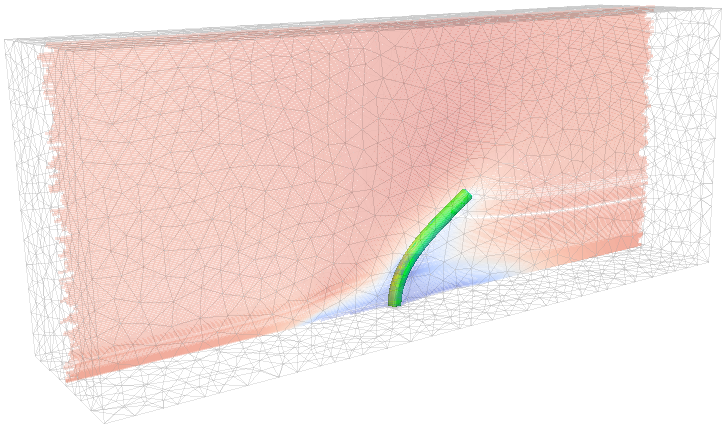}
        \caption{t=0.4}
        \label{fig:52e}
    \end{subfigure}
    \begin{subfigure}{0.45\linewidth}
        \includegraphics[width=\linewidth]{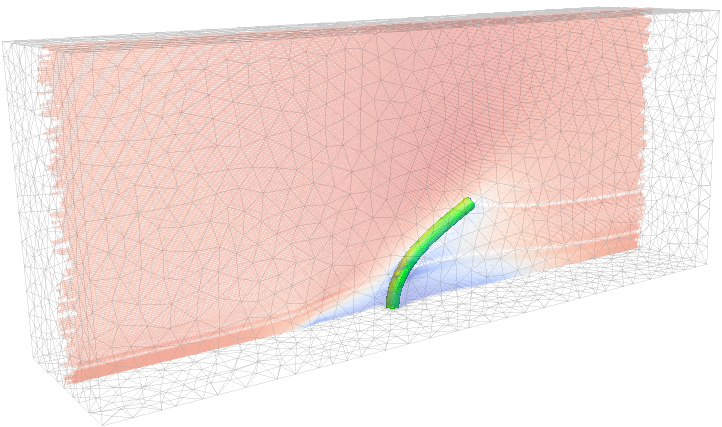}
        \caption{t=0.5}
        \label{fig:52f}
    \end{subfigure}
\\
    \begin{subfigure}{0.45\linewidth}
        \includegraphics[width=\linewidth]{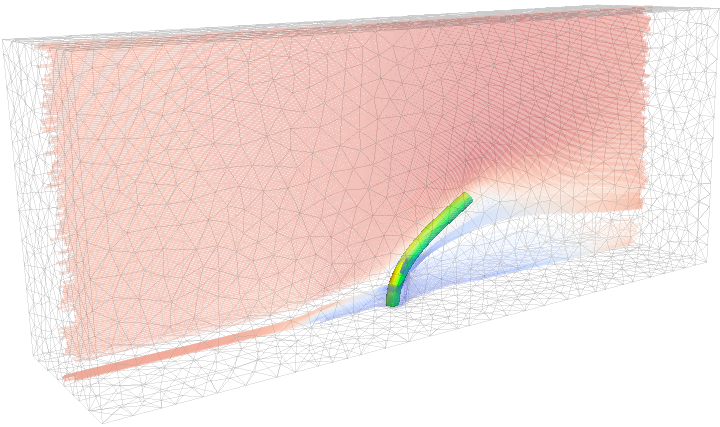}
        \caption{t=1}
        \label{fig:52g}
    \end{subfigure}
    \begin{subfigure}{0.45\linewidth}
        \includegraphics[width=\linewidth]{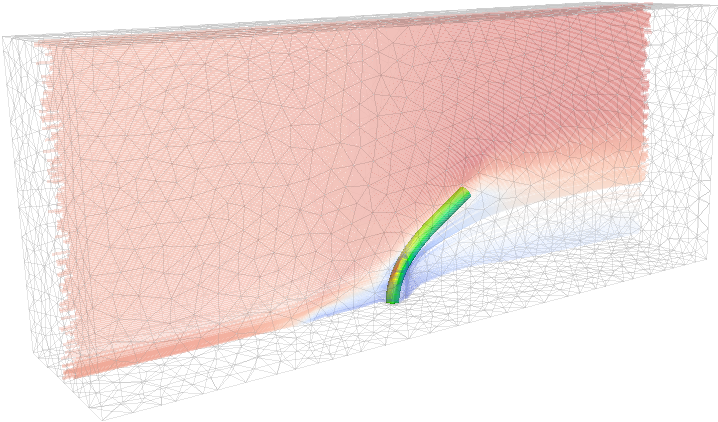}
        \caption{t=2}
        \label{fig:52h}
    \end{subfigure}
\caption{Fluid velocity and pressure at different positions in time.}
\label{fig:VDattimet}
\end{figure}

The normalized bending line is plotted and compared to the experimental results in Figure \ref{fig:BC}. 
A remarkable correlation is observed between the FSI simulation and the experimental data. 
This close agreement underscores the accuracy of our numerical approach in replicating the real-world behavior of the system.
In Figure \ref{fig:VDattimet}, we present the velocity streamlines surrounding the cylindrical structure at different time intervals (t). 
This visual representation vividly illustrates the dynamic flow patterns and displacement field during the simulated scenario.
For a more detailed analysis of the structure's dynamics, we present the x and y components of the displacement plotted against time in Figure \ref{fig:BCUy}. 
These plots provide valuable insights into how the cylinder responds to the fluid forces over time. 
It's evident that our simulation captures the intricate dynamics of the system, enhancing our understanding of its behavior.
These results affirm the suitability of our fully coupled two-way FSI simulation in 3D for this experimental setup, confirming the framework's validity and reliability. 
This numerical approach offers a valuable tool for investigating and understanding the complex interactions between fluid and solid structures, with promising applications in various fields.

\begin{figure}[p]
        \centering
        \begin{tikzpicture}[trim axis left, trim axis right]
                \begin{axis}[
                scale=0.8,
                transform shape,
                width=.7\textwidth,
                height=.5\textwidth,
                label style={font=\scriptsize},
                tick label style={font=\scriptsize},
                legend style={font=\scriptsize},
                clip=true,
                grid=both,
                legend cell align=left,
                legend pos=south east,
                grid=major,
                xlabel={time (s)},
                ylabel={x displacement of tip}  
                ]                 
                \addplot[ 
                mark=none,thick,mygray1
                ]
                table[
                x index=0,
                y index=1
                ]
                {data/BCUx.dat}; 
                \end{axis}
        \end{tikzpicture}
\end{figure}

\subsection{Bending of elastic flaps in a cross flow}

\begin{figure}[H]
\centering
\begin{minipage}{0.33\textwidth}
   \includegraphics[width=1\linewidth]{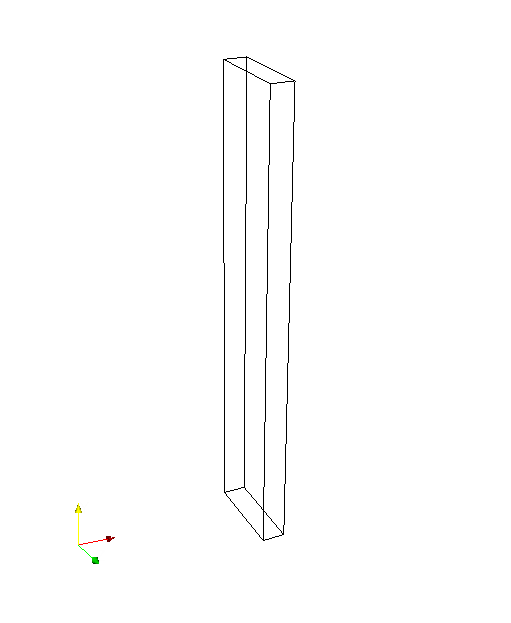}
    \caption{b= 5 mm} \label{fig:53a}
  \end{minipage}%
  \hspace*{\fill}   
  \begin{minipage}{0.33\textwidth}
   \includegraphics[width=1\linewidth]{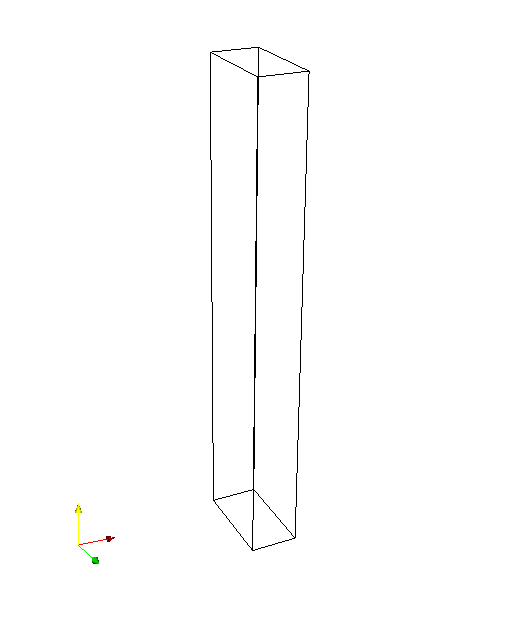}
    \caption{b= 10 mm} \label{fig:53b}
  \end{minipage}%
  \hspace*{\fill}   
  \begin{minipage}{0.33\textwidth}
   \includegraphics[width=1\linewidth]{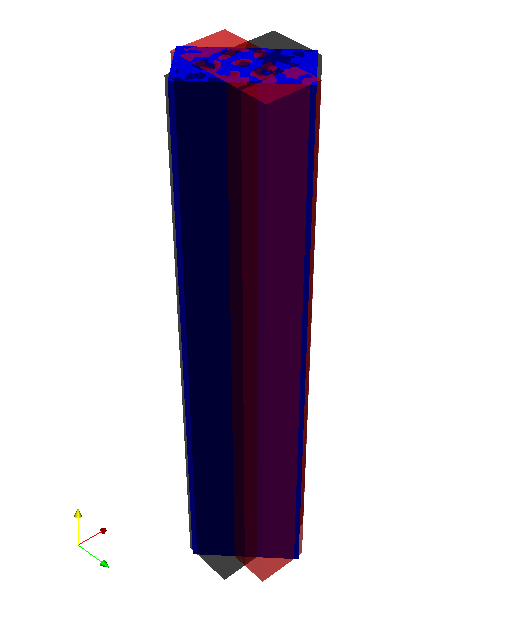}
    \caption{b= 10 mm, for angles 0, 45, and 90 °} \label{fig:53c}
  \end{minipage}%
\caption{Different flap geometries and orientations.}
\label{fig:FlapGeometry}
\end{figure}

\begin{figure}[p]
\centering
    \begin{subfigure}{0.45\linewidth}
        \includegraphics[width=\linewidth]{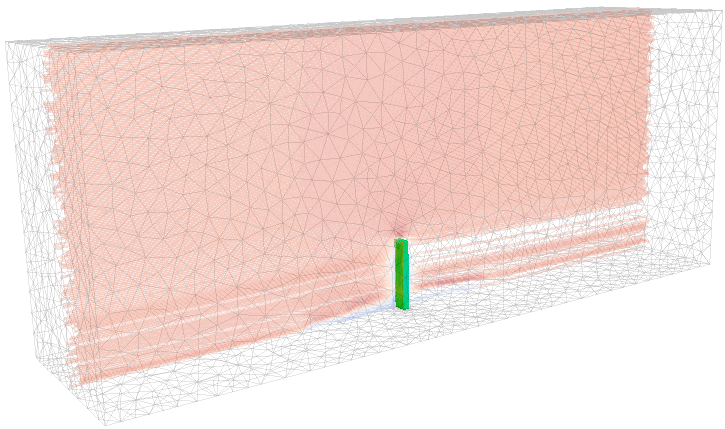}
        \caption{t=0}
        \label{fig:54a}
    \end{subfigure}
    \begin{subfigure}{0.45\linewidth}
        \includegraphics[width=\linewidth]{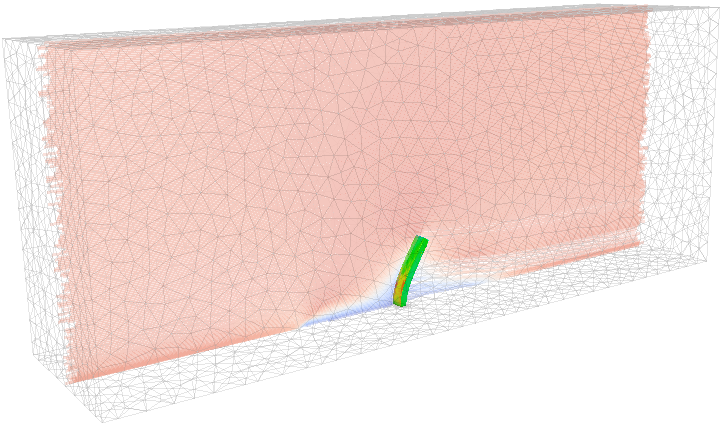}
        \caption{t=0.1}
        \label{fig:54b}
    \end{subfigure}
\\
    \begin{subfigure}{0.45\linewidth}
        \includegraphics[width=\linewidth]{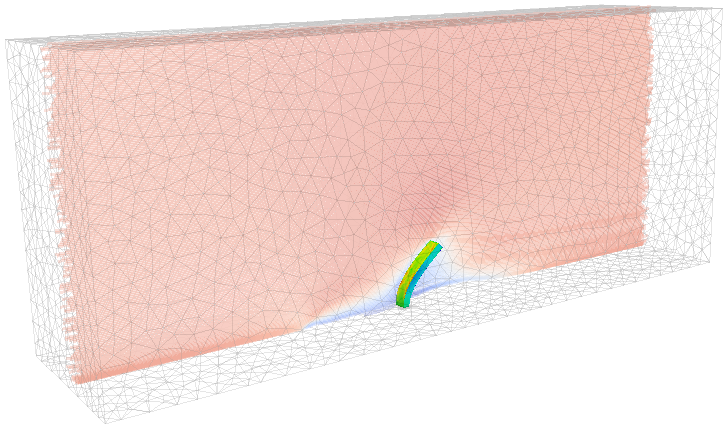}
        \caption{t=0.2}
        \label{fig:54c}
    \end{subfigure}
    \begin{subfigure}{0.45\linewidth}
        \includegraphics[width=\linewidth]{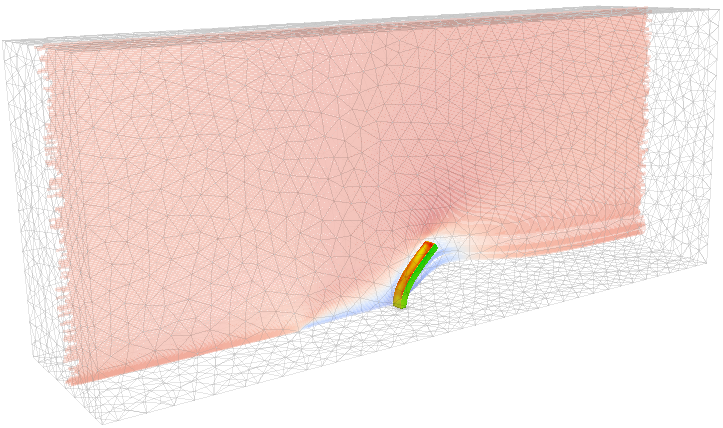}
        \caption{t=0.3}
        \label{fig:54d}
    \end{subfigure}
\\
    \begin{subfigure}{0.45\linewidth}
        \includegraphics[width=\linewidth]{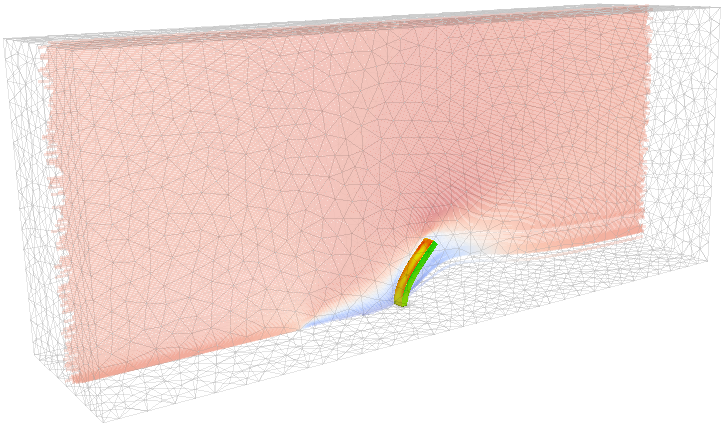}
        \caption{t=0.4}
        \label{fig:54e}
    \end{subfigure}
    \begin{subfigure}{0.45\linewidth}
        \includegraphics[width=\linewidth]{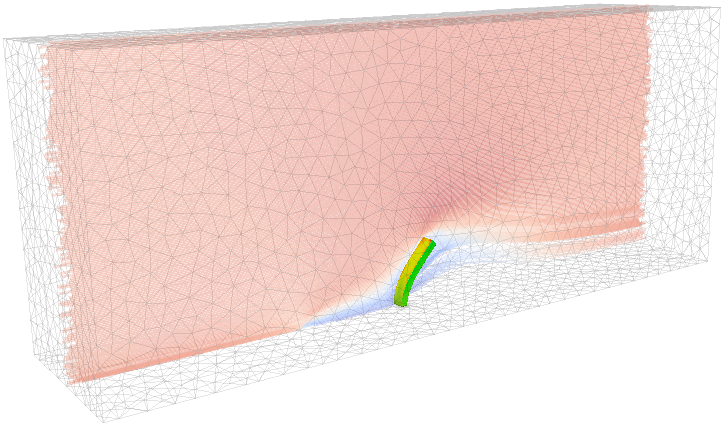}
        \caption{t=0.5}
        \label{fig:54f}
    \end{subfigure}
\\
    \begin{subfigure}{0.45\linewidth}
        \includegraphics[width=\linewidth]{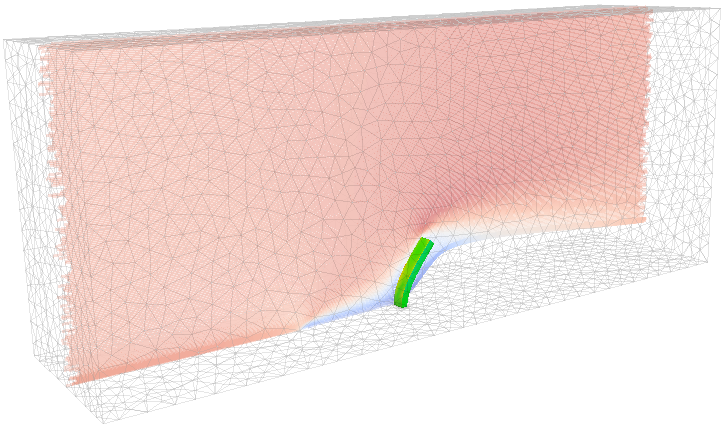}
        \caption{t=1}
        \label{fig:54g}
    \end{subfigure}
    \begin{subfigure}{0.45\linewidth}
        \includegraphics[width=\linewidth]{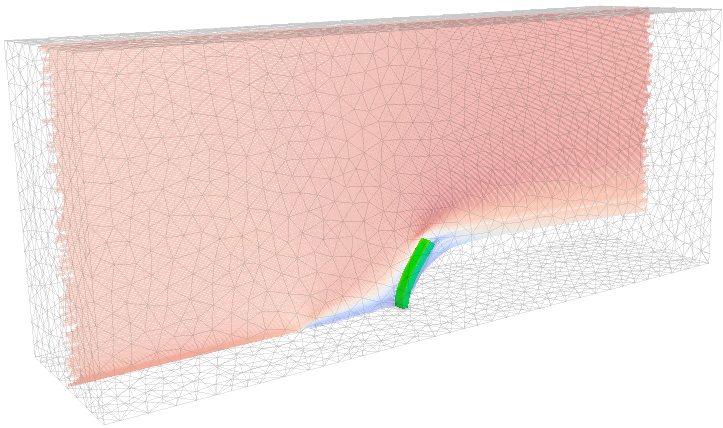}
        \caption{t=2}
        \label{fig:54h}
    \end{subfigure}
\caption{Fluid velocity and pressure at different positions in time for the 10mm flap.}
\label{fig:Flap10}
\end{figure}

\begin{figure}[p]
\centering
    \begin{subfigure}{0.45\linewidth}
        \includegraphics[width=\linewidth]{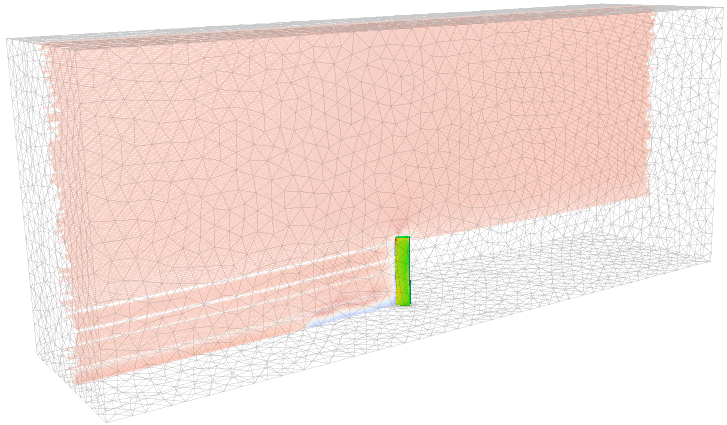}
        \caption{t=0}
        \label{fig:55a}
    \end{subfigure}
    \begin{subfigure}{0.45\linewidth}
        \includegraphics[width=\linewidth]{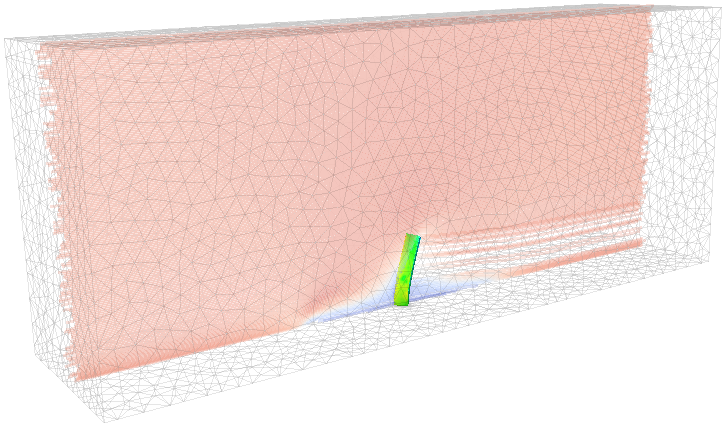}
        \caption{t=0.1}
        \label{fig:55b}
    \end{subfigure}
\\
    \begin{subfigure}{0.45\linewidth}
        \includegraphics[width=\linewidth]{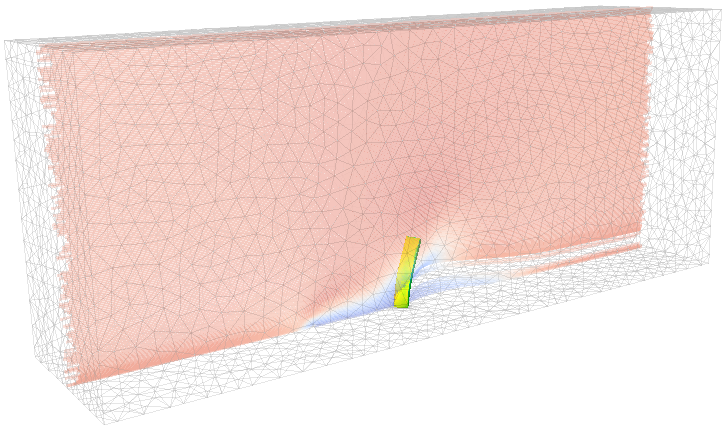}
        \caption{t=0.2}
        \label{fig:55c}
    \end{subfigure}
    \begin{subfigure}{0.45\linewidth}
        \includegraphics[width=\linewidth]{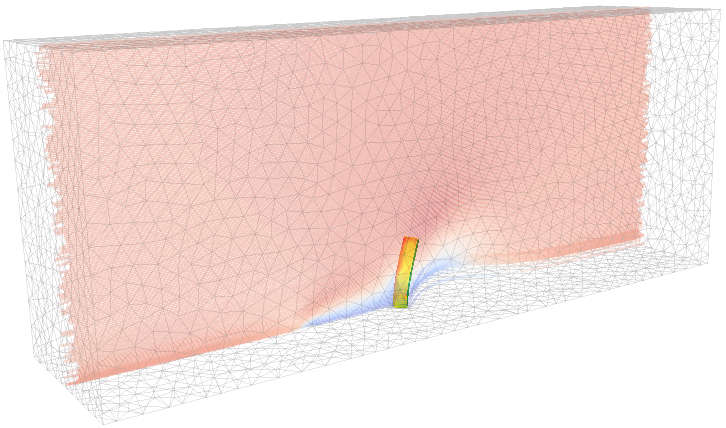}
        \caption{t=0.3}
        \label{fig:55d}
    \end{subfigure}
\\
    \begin{subfigure}{0.45\linewidth}
        \includegraphics[width=\linewidth]{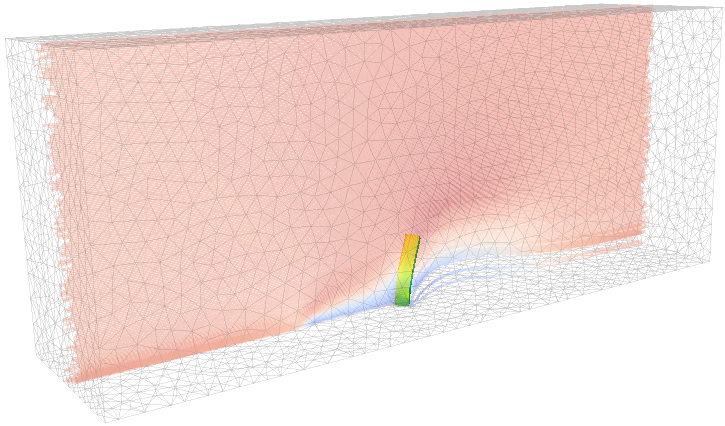}
        \caption{t=0.5}
        \label{fig:55e}
    \end{subfigure}
    \begin{subfigure}{0.45\linewidth}
        \includegraphics[width=\linewidth]{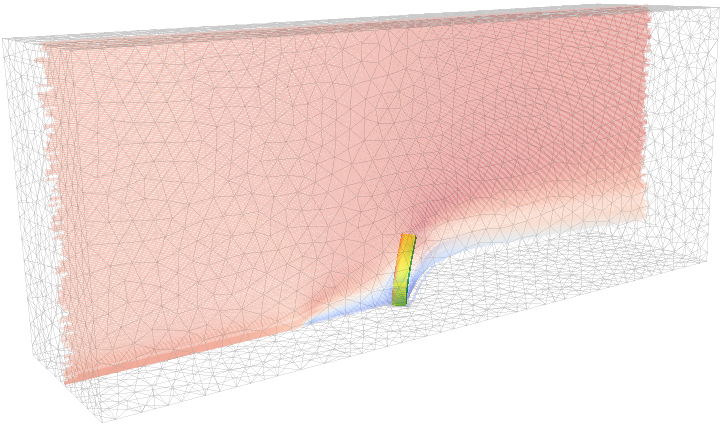}
        \caption{t=2}
        \label{fig:55f}
    \end{subfigure}
\caption{Fluid velocity and pressure at different positions in time for the 45 degrees variation of the 10mm flap.}
\label{fig:Flap1045}
\end{figure}

\begin{figure}[p]
\centering
    \begin{subfigure}{0.45\linewidth}
        \includegraphics[width=\linewidth]{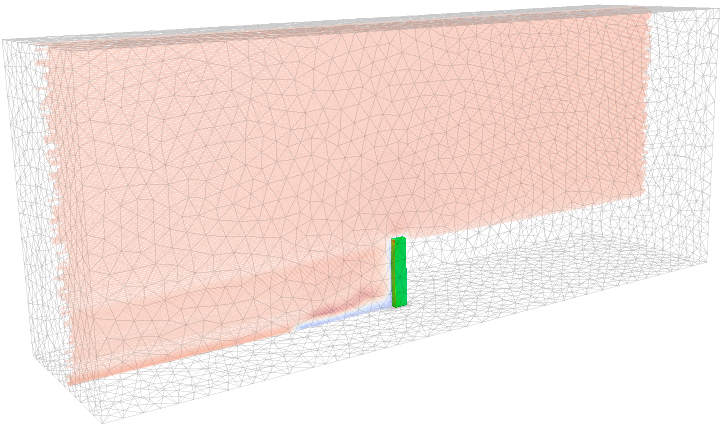}
        \caption{t=0}
        \label{fig:56a}
    \end{subfigure}
    \begin{subfigure}{0.45\linewidth}
        \includegraphics[width=\linewidth]{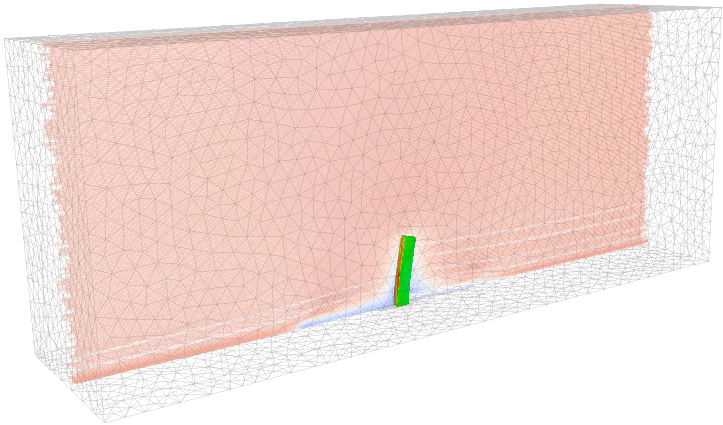}
        \caption{t=0.05}
        \label{fig:56b}
    \end{subfigure}
\\
    \begin{subfigure}{0.45\linewidth}
        \includegraphics[width=\linewidth]{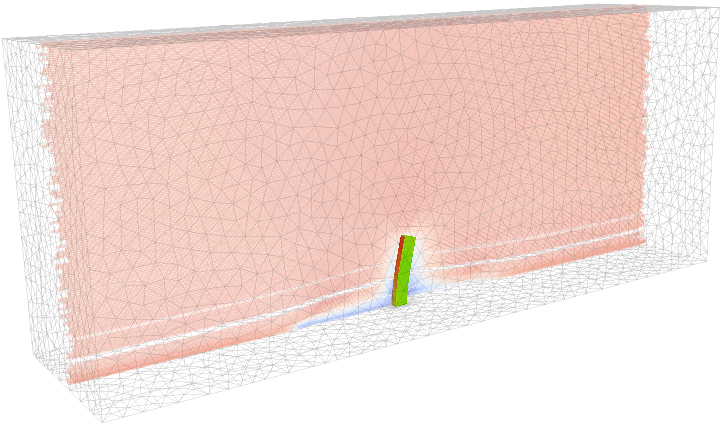}
        \caption{t=0.075}
        \label{fig:56c}
    \end{subfigure}
    \begin{subfigure}{0.45\linewidth}
        \includegraphics[width=\linewidth]{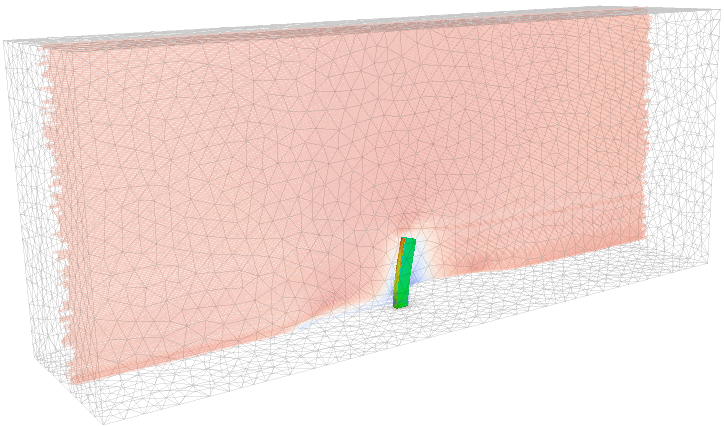}
        \caption{t=0.1}
        \label{fig:56d}
    \end{subfigure}
\\
    \begin{subfigure}{0.45\linewidth}
        \includegraphics[width=\linewidth]{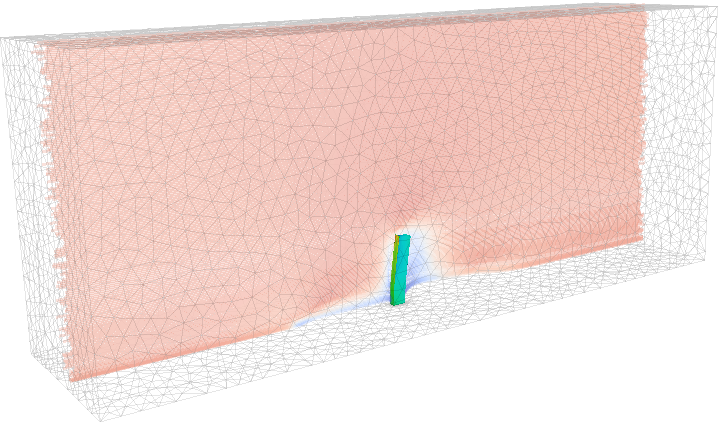}
        \caption{t=0.15}
        \label{fig:56e}
    \end{subfigure}
    \begin{subfigure}{0.45\linewidth}
        \includegraphics[width=\linewidth]{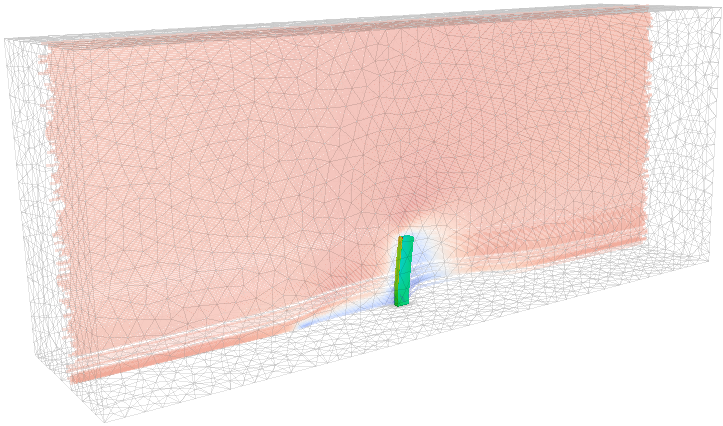}
        \caption{t=0.2}
        \label{fig:56f}
    \end{subfigure}
\caption{Fluid velocity and pressure at different positions in time for the 90 degrees variation of the 10mm flap.}
\label{fig:Flap1090}
\end{figure}

An intriguing variation of the previously discussed case was conducted in \cite{article33}, where flap-like structures are immersed in the fluid instead of the cylindrical beam. 
The fluid and solid properties of this problem are identical to the previous scenario, ensuring consistency in the physical characteristics of the system.
In this configuration, the solid flap structure features a rectangular geometry with a length of $l=100 mm$, a width of $w=20 mm$, and varying thicknesses of $b=5, 10 mm$. 
Figure \ref{fig:53a} illustrates the solid geometry for a thickness of $5 mm$, while Figure \ref{fig:53b} depicts the geometry for a thickness of $10 mm$. 
Furthermore, different orientations were considered for the case of the thick $10mm$ flap, where angles vary between 0, 45, and 90 degrees, as shown in Figure \ref{fig:53c}.
The results of this investigation are summarized in Table \ref{table:52}, which presents the maximum displacement magnitudes and steady-state displacement magnitudes for each flap orientation, along with a comparison to the benchmark maximum transient deformation.
To provide a visual representation of the flow patterns and structural dynamics, we present Figures \ref{fig:Flap10}, \ref{fig:Flap1045}, and \ref{fig:Flap1090}, which show streamlines highlighted by velocity magnitude around the structural flap, as well as the displacement magnitude for the different beam orientations of the $10mm$ thick flap at various time steps. 
These figures offer valuable insights into how different flap orientations and thicknesses impact the interaction between the fluid flow and the solid structure.
This study on the bending of elastic flaps in a cross flow enriches our understanding of complex fluid-structure interactions and provides insights into the behavior of flexible structures immersed in flowing fluids. 

\thispagestyle{empty}

\begin{table}[!h]
\centering
\begin{tabular}{|p{3cm}|p{3cm}|p{3cm}|p{3cm}|}
 \hline
 {Angle(°)} & {Maximum Displacement Magnitude (m)} & {Steady-state Displacement Magnitude (m)} & {Benchmark Maximum Transient Deformation (m)}\\
 \hline
0  & 0.06623 & 0.05082 & 0.0526\\
45 & 0.05581 & 0.04564 & 0.0463\\
90 & 0.01788 & 0.01104 & 0.0196\\
\hline
\end{tabular} 
\caption{Fluid and solid properties for the pillar in a Laminar cross flow problem.}
\label{table:52}
\end{table}

%% file: 06-Perspectives_andconclusion.tex
\section{Perspectives and conclusion}
\indent

This article introduces a novel hybrid Fluid-Structure Interaction (FSI) framework that amalgamates the merits of well-established methods in the existing literature. 
The Adaptive Immersed Mesh Method (AIMM) serves as a bridge between the Monolithic (Eulerian) and Partitioned (Lagrangian) FSI approaches, harnessing the strengths of both to enrich the solid model.
In AIMM, the solid mesh becomes immersed within the fluid-solid mesh at each time step through the level set method. 
This immersion facilitates the natural tracking of the fluid-solid interface on the fluid-solid mesh, where velocity and full stress are communicated at the interface of each grid. 
The framework employs anisotropic mesh adaptation under different criteria, including the creation of stretched elements at the interface based on the level set gradient, which significantly contributes to the precision and accuracy of results.
Furthermore, the Variational Multi-Scale (VMS) Method is applied to both solvers, enabling the use of first-order unstructured finite elements while ensuring compliance with the inf-sup condition. 
This comprehensive methodology is validated through a series of two-dimensional benchmarks, and the obtained results consistently align with the existing literature.
The culmination of our research is exemplified in a diverse array of three-dimensional FSI simulations, showcasing the framework's capability to accurately model flexible, relatively thin structures immersed in a fluid. 
As a testament to our ongoing commitment to advancing this work, we are actively extending the framework to tackle more complex applications, including challenging biomechanical scenarios.
This article presents the hybrid Fluid--Structure Interaction (FSI) framework that combines the advantages of other well-known methods in the literature.

%% file: main.bbl
\begin{thebibliography}{10}
\expandafter\ifx\csname url\endcsname\relax
  \def\url#1{\texttt{#1}}\fi
\expandafter\ifx\csname urlprefix\endcsname\relax\def\urlprefix{URL }\fi
\expandafter\ifx\csname href\endcsname\relax
  \def\href#1#2{#2} \def\path#1{#1}\fi

\bibitem{hamdan1995fluid}
F.~Hamdan, P.~Dowling, Fluid-structure interaction: application to structures in an acoustic fluid medium, part 1: an introduction to numerical treatment, Engineering Computations (1995).

\bibitem{hartwanger20083d}
D.~Hartwanger, A.~Horvat, 3d modelling of a wind turbine using cfd, in: NAFEMS Conference, United Kingdom, 2008.

\bibitem{hou2014coupled}
H.~Hou, Coupled fluid-structure analysis for exhaust system nvh, Tech. rep., SAE Technical Paper (2014).

\bibitem{trimarchi1970fluid}
D.~Trimarchi, S.~Turnock, D.~Chapelle, D.~Taunton, Fluid-structure interactions of anisotropic thin composite materials for application to sail aerodynamics of a yacht in waves (1970).

\bibitem{leung2006fluid}
J.~H. Leung, A.~R. Wright, N.~Cheshire, J.~Crane, S.~A. Thom, A.~D. Hughes, Y.~Xu, Fluid structure interaction of patient specific abdominal aortic aneurysms: a comparison with solid stress models, Biomedical engineering online 5~(1) (2006) 1--15.

\bibitem{deparis2003acceleration}
S.~Deparis, M.~A. Fern{\'a}ndez, L.~Formaggia, Acceleration of a fixed point algorithm for fluid-structure interaction using transpiration conditions, ESAIM: Mathematical Modelling and Numerical Analysis 37~(4) (2003) 601--616.

\bibitem{chouly2006simulation}
F.~Chouly, A.~Van~Hirtum, P.-Y. Lagr{\'e}e, J.-R. Paoli, X.~Pelorson, Y.~Payan, Simulation of the retroglossal fluid-structure interaction during obstructive sleep apnea, in: International Symposium on Biomedical Simulation, Springer, 2006, pp. 48--57.

\bibitem{garg2009numerically}
R.~Garg, C.~Narayanan, S.~Subramaniam, A numerically convergent lagrangian--eulerian simulation method for dispersed two-phase flows, International Journal of Multiphase Flow 35~(4) (2009) 376--388.

\bibitem{41e42816959649909db90a93179d8dbf}
R.~Glowinski, T.~Pan, T.~Hesla, D.~Joseph, A distributed lagrange multiplier/fictitious domain method for particulate flows, International Journal of Multiphase Flow 25~(5) (1999) 755--794.
\newblock \href {https://doi.org/10.1016/S0301-9322(98)00048-2} {\path{doi:10.1016/S0301-9322(98)00048-2}}.

\bibitem{doi:10.1137/040604728}
H.-O. Kreiss, N.~A. Petersson, \href{https://doi.org/10.1137/040604728}{A second order accurate embedded boundary method for the wave equation with dirichlet data}, SIAM Journal on Scientific Computing 27~(4) (2006) 1141--1167.
\newblock \href {http://arxiv.org/abs/https://doi.org/10.1137/040604728} {\path{arXiv:https://doi.org/10.1137/040604728}}, \href {https://doi.org/10.1137/040604728} {\path{doi:10.1137/040604728}}.
\newline\urlprefix\url{https://doi.org/10.1137/040604728}

\bibitem{peskin_2002}
C.~S. Peskin, The immersed boundary method, Acta Numerica 11 (2002) 479–517.
\newblock \href {https://doi.org/10.1017/S0962492902000077} {\path{doi:10.1017/S0962492902000077}}.

\bibitem{wang:hal-00651118}
K.~Wang, A.~Rallu, J.-F. Gerbeau, C.~Farhat, \href{https://hal.inria.fr/hal-00651118}{{Algorithms for Interface Treatment and Load Computation in Embedded Boundary Methods for Fluid and Fluid-Structure Interaction Problems}}, {International Journal for Numerical Methods in Fluids} 67~(9) (2011) 1175--1206.
\newblock \href {https://doi.org/10.1002/fld.2556} {\path{doi:10.1002/fld.2556}}.
\newline\urlprefix\url{https://hal.inria.fr/hal-00651118}

\bibitem{persson2009curved}
P.-O. Persson, J.~Peraire, Curved mesh generation and mesh refinement using lagrangian solid mechanics, in: 47th AIAA Aerospace Sciences Meeting including The New Horizons Forum and Aerospace Exposition, 2009, p. 949.

\bibitem{donea1977lagrangian}
J.~Don{\'e}a, P.~Fasoli-Stella, S.~Giuliani, Lagrangian and eulerian finite element techniques for transient fluid-structure interaction problems (1977).

\bibitem{HUERTA1988277}
A.~Huerta, W.~K. Liu, \href{https://www.sciencedirect.com/science/article/pii/0045782588900448}{Viscous flow with large free surface motion}, Computer Methods in Applied Mechanics and Engineering 69~(3) (1988) 277--324.
\newblock \href {https://doi.org/https://doi.org/10.1016/0045-7825(88)90044-8} {\path{doi:https://doi.org/10.1016/0045-7825(88)90044-8}}.
\newline\urlprefix\url{https://www.sciencedirect.com/science/article/pii/0045782588900448}

\bibitem{hughes1981lagrangian}
T.~J. Hughes, W.~K. Liu, T.~K. Zimmermann, Lagrangian-eulerian finite element formulation for incompressible viscous flows, Computer methods in applied mechanics and engineering 29~(3) (1981) 329--349.

\bibitem{GERSTENBERGER20081699}
A.~Gerstenberger, W.~A. Wall, \href{https://www.sciencedirect.com/science/article/pii/S0045782507002915}{An extended finite element method/lagrange multiplier based approach for fluid–structure interaction}, Computer Methods in Applied Mechanics and Engineering 197~(19) (2008) 1699--1714, computational Methods in Fluid–Structure Interaction.
\newblock \href {https://doi.org/https://doi.org/10.1016/j.cma.2007.07.002} {\path{doi:https://doi.org/10.1016/j.cma.2007.07.002}}.
\newline\urlprefix\url{https://www.sciencedirect.com/science/article/pii/S0045782507002915}

\bibitem{wall2009advances}
W.~Wall, A.~Gerstenberger, U.~Mayer, Advances in fixed-grid fluid structure interaction, in: ECCOMAS Multidisciplinary Jubilee Symposium, Springer, 2009, pp. 235--249.

\bibitem{wall2006large}
W.~A. Wall, A.~Gerstenberger, P.~Gamnitzer, C.~F{\"o}rster, E.~Ramm, Large deformation fluid-structure interaction--advances in ale methods and new fixed grid approaches, in: Fluid-structure interaction, Springer, 2006, pp. 195--232.

\bibitem{legay2006eulerian}
A.~Legay, J.~Chessa, T.~Belytschko, An eulerian--lagrangian method for fluid--structure interaction based on level sets, Computer Methods in Applied Mechanics and Engineering 195~(17-18) (2006) 2070--2087.

\bibitem{article}
A.~Legay, J.~Chessa, T.~Belytschko, An eulerian-lagrangian method for fluid-structure interaction based on level sets, Computer Methods in Applied Mechanics and Engineering 195 (2006) 2070--2087.
\newblock \href {https://doi.org/10.1016/j.cma.2005.02.025} {\path{doi:10.1016/j.cma.2005.02.025}}.

\bibitem{zilian2008enriched}
A.~Zilian, A.~Legay, The enriched space--time finite element method (est) for simultaneous solution of fluid--structure interaction, International Journal for Numerical Methods in Engineering 75~(3) (2008) 305--334.

\bibitem{article2}
T.~Wick, Flapping and contact fsi computations with the fluid-solid interface-tracking/interface-capturing technique and mesh adaptivity, Computational Mechanics 53 (10 2013).
\newblock \href {https://doi.org/10.1007/s00466-013-0890-3} {\path{doi:10.1007/s00466-013-0890-3}}.

\bibitem{NEMER2021113923}
R.~Nemer, A.~Larcher, T.~Coupez, E.~Hachem, \href{https://www.sciencedirect.com/science/article/pii/S0045782521002607}{Stabilized finite element method for incompressible solid dynamics using an updated lagrangian formulation}, Computer Methods in Applied Mechanics and Engineering 384 (2021) 113923.
\newblock \href {https://doi.org/https://doi.org/10.1016/j.cma.2021.113923} {\path{doi:https://doi.org/10.1016/j.cma.2021.113923}}.
\newline\urlprefix\url{https://www.sciencedirect.com/science/article/pii/S0045782521002607}

\bibitem{babuvska1971error}
I.~Babu{\v{s}}ka, Error-bounds for finite element method, Numerische Mathematik 16~(4) (1971) 322--333.

\bibitem{HUGHES198685}
T.~J. Hughes, L.~P. Franca, M.~Balestra, A new finite element formulation for computational fluid dynamics: V. circumventing the babuška-brezzi condition: a stable petrov-galerkin formulation of the stokes problem accommodating equal-order interpolations, Computer Methods in Applied Mechanics and Engineering 59~(1) (1986) 85 -- 99.

\bibitem{HACHEM20108643}
E.~Hachem, B.~Rivaux, T.~Kloczko, H.~Digonnet, T.~Coupez, Stabilized finite element method for incompressible flows with high reynolds number, Journal of Computational Physics 229~(23) (2010) 8643 -- 8665.

\bibitem{HUGHES198785}
T.~J. Hughes, L.~P. Franca, A new finite element formulation for computational fluid dynamics: Vii. the stokes problem with various well-posed boundary conditions: Symmetric formulations that converge for all velocity/pressure spaces, Computer Methods in Applied Mechanics and Engineering 65~(1) (1987) 85 -- 96.

\bibitem{10.1007/BF01395881}
L.~P. Franca, T.~J. Hughes, A.~F. Loula, I.~Miranda, A new family of stable elements for nearly incompressible elasticity based on a mixed petrov-galerkin finite element formulation, Numer. Math. 53~(1–2) (1988) 123–141.

\bibitem{hughes1998variational}
T.~J. Hughes, G.~R. Feij{\'o}o, L.~Mazzei, J.-B. Quincy, The variational multiscale method—a paradigm for computational mechanics, Computer methods in applied mechanics and engineering 166~(1-2) (1998) 3--24.

\bibitem{dubois1999dynamic}
T.~Dubois, Y.~Matras, F.~Jauberteau, R.~Temam, Dynamic multilevel methods and the numerical simulation of turbulence, Cambridge University Press, 1999.

\bibitem{https://doi.org/10.1002/fld.1481}
R.~Codina, J.~Principe, \href{https://onlinelibrary.wiley.com/doi/abs/10.1002/fld.1481}{Dynamic subscales in the finite element approximation of thermally coupled incompressible flows}, International Journal for Numerical Methods in Fluids 54~(6‐8) (2007) 707--730.
\newblock \href {http://arxiv.org/abs/https://onlinelibrary.wiley.com/doi/pdf/10.1002/fld.1481} {\path{arXiv:https://onlinelibrary.wiley.com/doi/pdf/10.1002/fld.1481}}, \href {https://doi.org/https://doi.org/10.1002/fld.1481} {\path{doi:https://doi.org/10.1002/fld.1481}}.
\newline\urlprefix\url{https://onlinelibrary.wiley.com/doi/abs/10.1002/fld.1481}

\bibitem{HACHEM2016238}
E.~Hachem, M.~Khalloufi, J.~Bruchon, R.~Valette, Y.~Mesri, \href{https://www.sciencedirect.com/science/article/pii/S0045782516304236}{Unified adaptive variational multiscale method for two phase compressible–incompressible flows}, Computer Methods in Applied Mechanics and Engineering 308 (2016) 238--255.
\newblock \href {https://doi.org/https://doi.org/10.1016/j.cma.2016.05.022} {\path{doi:https://doi.org/10.1016/j.cma.2016.05.022}}.
\newline\urlprefix\url{https://www.sciencedirect.com/science/article/pii/S0045782516304236}

\bibitem{CODINA20024295}
R.~Codina, \href{https://www.sciencedirect.com/science/article/pii/S0045782502003377}{Stabilized finite element approximation of transient incompressible flows using orthogonal subscales}, Computer Methods in Applied Mechanics and Engineering 191~(39) (2002) 4295--4321.
\newblock \href {https://doi.org/https://doi.org/10.1016/S0045-7825(02)00337-7} {\path{doi:https://doi.org/10.1016/S0045-7825(02)00337-7}}.
\newline\urlprefix\url{https://www.sciencedirect.com/science/article/pii/S0045782502003377}

\bibitem{bruchon:emse-00475556}
J.~Bruchon, H.~Digonnet, T.~Coupez, \href{https://hal-emse.ccsd.cnrs.fr/emse-00475556}{{Using a signed distance function for the simulation of metal forming processes : Formulation of the contact condition and mesh adaptation. From a Lagrangian approach to an Eulerian approach}}, {International Journal for Numerical Methods in Engineering} 78~(8) (2009) 980--1008.
\newblock \href {https://doi.org/10.1002/nme.2519} {\path{doi:10.1002/nme.2519}}.
\newline\urlprefix\url{https://hal-emse.ccsd.cnrs.fr/emse-00475556}

\bibitem{article1}
R.~Codina, O.~Soto, A numerical model to track two‐fluid interfaces based on a stabilized finite element method and the level set technique, International Journal for Numerical Methods in Fluids 40 (2002) 293 -- 301.
\newblock \href {https://doi.org/10.1002/fld.277} {\path{doi:10.1002/fld.277}}.

\bibitem{formaggia2003anisotropic}
L.~Formaggia, S.~Perotto, Anisotropic error estimates for elliptic problems, Numerische Mathematik 94~(1) (2003) 67--92.

\bibitem{hoffman2003adaptive}
J.~Hoffman, C.~Johnson, Adaptive finite element methods for incompressible fluid flow, Error Estimation and Adaptive Discretization in Computational Fluid Dynamics (Ed. T. J. Barth and H. Deconinck), Lecture Notes in Computational Science and Engineering Vol.25, Springer-Verlag Publishing, Heidelberg, 2003.

\bibitem{https://doi.org/10.1002/nme.4481}
E.~Hachem, S.~Feghali, R.~Codina, T.~Coupez, \href{https://onlinelibrary.wiley.com/doi/abs/10.1002/nme.4481}{Immersed stress method for fluid–structure interaction using anisotropic mesh adaptation}, International Journal for Numerical Methods in Engineering 94~(9) (2013) 805--825.
\newblock \href {http://arxiv.org/abs/https://onlinelibrary.wiley.com/doi/pdf/10.1002/nme.4481} {\path{arXiv:https://onlinelibrary.wiley.com/doi/pdf/10.1002/nme.4481}}, \href {https://doi.org/https://doi.org/10.1002/nme.4481} {\path{doi:https://doi.org/10.1002/nme.4481}}.
\newline\urlprefix\url{https://onlinelibrary.wiley.com/doi/abs/10.1002/nme.4481}

\bibitem{Khalloufi201644}
M.~Khalloufi, Y.~Mesri, R.~Valette, E.~Massoni, E.~Hachem, High fidelity anisotropic adaptive variational multiscale method for multiphase flows with surface tension, Computer Methods in Applied Mechanics and Engineering 307 (2016) 44 -- 67.

\bibitem{Legrain2018}
G.~Legrain, N.~Mo\"es, Adaptive anisotropic integration scheme for high-order fictitious domain methods: Application to thin structures, International Journal for Numerical Methods in Engineering 114~(8) (2018) 882--904.

\bibitem{Coupez11}
T.~Coupez, Metric construction by length distribution tensor and edge based error for anisotropic adaptive meshing, Journal of Computational Physics 230~(7) (2011) 2391--2405.

\bibitem{kunert2000edge}
G.~Kunert, R.~Verf{\"u}rth, Edge residuals dominate a posteriori error estimates for linear finite element methods on anisotropic triangular and tetrahedral meshes, Numer. Math. 86~(2) (2000) 283--303.

\bibitem{Baiges}
J.~Baiges, R.~Codina, The fixed-mesh ale approach applied to solid mechanics and fluid-structure interaction problems, International Journal for Numerical Methods in Engineering 81 (2009) 1529 -- 1557.
\newblock \href {https://doi.org/10.1002/nme.2740} {\path{doi:10.1002/nme.2740}}.

\bibitem{https://doi.org/10.1002/nme.6321}
A.~Tello, R.~Codina, J.~Baiges, \href{https://onlinelibrary.wiley.com/doi/abs/10.1002/nme.6321}{Fluid structure interaction by means of variational multiscale reduced order models}, International Journal for Numerical Methods in Engineering 121~(12) (2020) 2601--2625.
\newblock \href {http://arxiv.org/abs/https://onlinelibrary.wiley.com/doi/pdf/10.1002/nme.6321} {\path{arXiv:https://onlinelibrary.wiley.com/doi/pdf/10.1002/nme.6321}}, \href {https://doi.org/https://doi.org/10.1002/nme.6321} {\path{doi:https://doi.org/10.1002/nme.6321}}.
\newline\urlprefix\url{https://onlinelibrary.wiley.com/doi/abs/10.1002/nme.6321}

\bibitem{HAN2020106179}
D.~Han, G.~Liu, S.~Abdallah, \href{https://www.sciencedirect.com/science/article/pii/S0045794919303372}{An eulerian-lagrangian-lagrangian method for 2d fluid-structure interaction problem with a thin flexible structure immersed in fluids}, Computers \& Structures 228 (2020) 106179.
\newblock \href {https://doi.org/https://doi.org/10.1016/j.compstruc.2019.106179} {\path{doi:https://doi.org/10.1016/j.compstruc.2019.106179}}.
\newline\urlprefix\url{https://www.sciencedirect.com/science/article/pii/S0045794919303372}

\bibitem{article3}
Z.-Q. Zhang, G.~Liu, B.~Khoo, Immersed smoothed finite element method for two dimensional fluid–structure interaction problems, International Journal for Numerical Methods in Engineering 90 (06 2012).
\newblock \href {https://doi.org/10.1002/nme.4299} {\path{doi:10.1002/nme.4299}}.

\bibitem{xia2008unstructured}
G.~Xia, C.-L. Lin, An unstructured finite volume approach for structural dynamics in response to fluid motions, Computers \& structures 86~(7-8) (2008) 684--701.

\bibitem{schafer1996benchmark}
M.~Sch{\"a}fer, S.~Turek, F.~Durst, E.~Krause, R.~Rannacher, Benchmark computations of laminar flow around a cylinder, in: Flow simulation with high-performance computers II, Springer, 1996, pp. 547--566.

\bibitem{turek2006proposal}
S.~Turek, J.~Hron, Proposal for numerical benchmarking of fluid-structure interaction between an elastic object and laminar incompressible flow, in: Fluid-structure interaction, Springer, 2006, pp. 371--385.

\bibitem{axtmann2016investigation}
G.~Axtmann, F.~Hegner, C.~Br{\"u}cker, U.~Rist, Investigation and prediction of the bending of single and tandem pillars in a laminar cross flow, Journal of Fluids and Structures 66 (2016) 110--126.

\bibitem{article33}
T.~Bano, F.~Hegner, M.~Heinrich, R.~Schwarze, Investigation of fluid-structure interaction induced bending for elastic flaps in a cross flow, Applied Sciences 10 (2020) 6177.
\newblock \href {https://doi.org/10.3390/app10186177} {\path{doi:10.3390/app10186177}}.

\end{thebibliography}
